\documentclass[12pt,a4paper]{amsart}

\title{Structure-preserving scheme for fractional nonlinear diffusion equations}
\author{Hélène Hivert, Florian Salin}
\address{Univ Rennes, Inria, Géosciences Rennes - UMR 6118, IRMAR - UMR 6625, F-35000 Rennes, France.}
\email{helene.hivert@inria.fr}
\address{Mathematical Institute and Graduate School of Science, Tohoku University, 6-3 Aoba, Aramaki, Aoba-ku, Sendai 980-8578, Japan}
\address{Univ Lyon, \'Ecole centrale de Lyon, CNRS UMR 5208, Institut Camille Jordan, F-69134 \'Ecully, France.}
\email{salin.florian.denis.jacques.r2@dc.tohoku.ac.jp}

\thanks{F.S.~is supported by JSPS KAKENHI Grant Number JP24KJ0381.}
\date{\today}

\usepackage{fullpage}
\usepackage[utf8]{inputenc}
\usepackage{amsmath,amssymb}
\usepackage{amsthm}
\usepackage{latexsym}
\usepackage{enumitem}
\usepackage{comment}
\usepackage{appendix}
\usepackage{mathtools}
\usepackage[giveninits=true]{biblatex}
\usepackage{hyperref}
\usepackage{xcolor}
\usepackage{tikz}
\usepackage{caption}
\usepackage{subcaption}
\usepackage{pgfplots}
\pgfplotsset{compat=1.18}

\newtheorem{thm}{Theorem}
\newtheorem{prop}[thm]{Proposition}
\newtheorem{lem}[thm]{Lemma}

\newtheorem{dfn}[thm]{Definition}
\theoremstyle{definition}

\newtheorem{rem}[thm]{Remark}

\newenvironment{prf}[1]
   {{\noindent \bf Proof of {#1}.}}{\hfill \qed}
\newenvironment{pr}
   {{\noindent \bf Proof.   }}{\hfill \qed}

\renewcommand{\th}{\theta}              
\newcommand{\Z}{\mathbb {Z}}    
\newcommand{\N}{\mathbb{N} }
\newcommand{\R}{\mathbb{R}}

\def\<{\langle }

\renewcommand{\d}{\mathrm{d}}

\newcommand{\X}{\mathcal{X}^\theta_h(I)}
\newcommand{\Xh}{\mathcal{X}^\theta_h(I)}
\renewcommand{\lq}{\ell^q_h(\mathbb{R})}

\newcommand{\CDP}{{\sf (CDP)}}
\newcommand{\bu}{\mathbf{u}}
\newcommand{\bv}{\mathbf{v}}
\newcommand{\bof}{\mathbf{f}}
\newcommand{\bog}{\mathbf{g}}
\newcommand{\bor}{\boldsymbol{\rho}}
\newcommand{\boe}{\boldsymbol{\eta}}

\usepackage{ifthen}
\newboolean{details}
\newcommand{\details}[1]{\ifthenelse{\boolean{details}}{\textcolor{blue}{#1}}{\ignorespaces}}
\setboolean{details}{false}

\begin{document}

\begin{abstract}
    In this paper, we introduce and analyze a numerical scheme for solving the Cauchy-Dirichlet problem associated with fractional nonlinear diffusion equations. These equations generalize the porous medium equation and the fast diffusion equation by incorporating a fractional diffusion term. We provide a rigorous analysis showing that the discretization preserves main properties of the continuous equations, including algebraic decay in the fractional porous medium case and the extinction phenomenon in the fractional fast diffusion case. The study is supported by extensive numerical simulations. In addition, we propose a novel method for accurately computing the extinction time for the fractional fast diffusion equation and illustrate numerically the convergence of rescaled solutions towards asymptotic profiles near the extinction time.
\end{abstract}

\maketitle


\section{Introduction}
The aim of the present paper is to investigate quantitative and qualitative properties of a numerical scheme for the fractional nonlinear diffusion equation
\begin{equation}\label{Eqn:CDPAvantChgVar}
\left\lbrace\begin{aligned}
\partial_t \rho+(-\Delta)^\theta \rho^m &=0 &&\text{in $\Omega\times(0,+\infty)$},\\
\rho&=0																					&&\text{in $(\R^d\setminus\Omega)\times(0,+\infty)$},\\
\rho(\cdot,0)&=\rho^0																&&\text{in $\Omega$},
\end{aligned}\right. 
\end{equation}
posed on a smooth bounded domain $\Omega\subset \R^d$, where $1<m<+\infty$ and $0<\theta<1$. Possibly sign-changing solutions are considered, and throughout this paper $\rho^{m}$ is a shorthand for the odd power $\rho^{m}:=|\rho|^{m-1}\rho$. In \eqref{Eqn:CDPAvantChgVar}, $(-\Delta)^\theta$ is the fractional Laplacian, defined by
\begin{equation}\label{Eqn:DefFracLap}
(-\Delta)^\theta\varphi(x) =\mathcal{F}^{-1}[|\xi|^{2\theta}\mathcal{F}\varphi](x)=C_{d,\theta}\:\text{P.V.}\int_{\R^d}\frac{\varphi(x)-\varphi(y)}{|x-y|^{d+2\theta}}{\rm d}y,
\end{equation}
for $x\in\R$ and e.g.~$\varphi\in C_c^\infty(\R^d)$, and where
$\text{P.V.}$ denotes the principal value,
\begin{equation}
(-\Delta)^\theta\varphi(x)=C_{d,\theta}\:\lim_{\varepsilon\rightarrow 0}\int_{\R^d\setminus B_\varepsilon(x)}\frac{\varphi(x)-\varphi(y)}{|x-y|^{d+2\theta}}{\rm d}y.
\end{equation}
\details{The constant $C_{d,\theta}$ is defined by
$$
C_{d,\theta}=\left(\int_{\R^N}\frac{1-\cos(\xi_1)}{|\xi|^{d+2\theta}}{\rm d}\xi\right)^{-1}=\frac{2\theta 2^{2\theta-1}\Gamma\left(\frac{d+2\theta}{2}\right)}{\pi^{d/2}\Gamma\left(\frac{2-2\theta}{2}\right)},
$$
where $\xi=(\xi_1,\dots,\xi_d)$.}
\noindent
Since the fractional Laplacian is a non-local operator, different definitions exist on bounded domains (e.g.~spectral Laplacian, censored fractional Laplacian \cite{SFL,CFL}). Here, we consider the \emph{restricted fractional Laplacian} which consists in imposing a boundary condition not only on $\partial \Omega$, but on $\R^d\setminus\Omega$ \cite{bonforte_sire_vazquez, akagi_existence}. 

In what follows, we use the change of variable $u(x,t)=(|\rho|^{m-1}\rho)(x,t)$ and $q=1/m+1$, and we let $u^0 = |\rho^0|^{m-1}\rho^0$, so that \eqref{Eqn:CDPAvantChgVar} can be rewritten as
\begin{equation}\label{Eqn:CDP}\tag{{\sf CDP}}
\left\lbrace\begin{aligned}
\partial_t u^{q-1}+(-\Delta)^\theta u&=0 &&\text{in $\Omega\times(0,+\infty)$},\\
u&=0																					&&\text{in $(\R^d\setminus\Omega)\times(0,+\infty)$},\\
u(\cdot,0)&=u^0																&&\text{in $\Omega$}.
\end{aligned}\right. 
\end{equation}
The case $1<q<2$, is usually referred to as the \emph{fractional porous medium equation} and the case $q>2$ is the \emph{fractional fast diffusion equation}. The case $q=2$ corresponds to a \emph{fractional heat equation}. 

In this article, we are interested in developing a numerical scheme for the fractional porous medium equation and the fractional fast diffusion equation that respects standard decay estimates, and behaviour for these two equations. In particular, it is well-known that solutions to \eqref{Eqn:CDP} decays in time algebraically in the porous medium case and extincts in finite time in the fast diffusion case. In this article, we shall write a numerical scheme that respects these two properties, and, in particular, that allows to approximate the extinction time in the fast diffusion case. Moreover, it has also been shown \cite{bonforte_figalli_survey, akagi_existence} that rescaled solutions to \eqref{Eqn:CDP} converge to asymptotic profiles, when the time goes to infinity in the porous medium case, or when it goes to the extinction time in the fast diffusion case. By determining the extinction time, our scheme will also allow us to numerically investigate the behaviour in large time of \eqref{Eqn:CDP}, and in particular the rate of convergence, which has been determined in the non-fractional case \cite{bonforte_figalli_survey,bonforte_figalli,akagi_rate}, but seems to be still an open problem for the fractional case.

More precisely, it can be shown (see e.g.~\cite{bonforte_vazquez, bonforte_sire_vazquez,bonforte_ibarrondo_ispizua,akagi_existence}) that for any $q\in(1,2_\theta^*]\setminus\{2\}$, where $2_\theta^*:=2d/(d-2\theta)_+$, it holds 
\begin{equation}\label{Eqn:Asymp1Continue}
\left(\|u^0\|_{L^q(\Omega)}^{-(2-q)}+\frac{2-q}{q-1}R(0)t\right)_+^{\frac{-1}{2-q}}\leq \|u(t)\|_{L^q(\Omega)}\leq \left(\|u^0\|_{L^q(\Omega)}^{-(2-q)}+\frac{2-q}{q-1}C^{-2}t\right)_+^{\frac{-1}{2-q}},
\end{equation}
where $a_+ :x= \max(a,0)$ for any $a\in\R$, and $R$ is the Rayleigh quotient defined by
\begin{equation}\label{Eqn:DefRayleigh}
    R(t):= [u(t)]_{H^\theta(\mathbb{R}^d)}^2/\|u(t)\|_{L^q(\Omega)}^2.
\end{equation}
The constant $C$ corresponds to any constant satisfying the Poincaré-Sobolev inequality 
\begin{equation}\label{Eqn:PoincareSobolevContinue}
\|u\|_{L^q(\Omega)}\leq C [u]_{H^\theta(\mathbb{R}^d)}\quad \text{for any }u\in H^\theta(\Omega)\text{ with }u\equiv 0 \text{ in }\mathbb{R}^d\setminus\Omega,
\end{equation}
and $[u]_{H^\theta(\mathbb{R}^d)}$ denotes the Gagliardo semi-norm defined by
$$
[u]_{H^\theta(\mathbb{R}^d)}:= \left(\int_{\mathbb{R}^d\times\mathbb{R}^d}\frac{|f(x)-f(y)|^2}{|x-y|^{d+2\theta}}{\rm{d}}x{\rm{d}}y\right)^{1/2}.
$$
Estimate \eqref{Eqn:Asymp1Continue} holds both for the porous medium case and the fast diffusion case. In the porous medium case $1<q<2$, \eqref{Eqn:Asymp1Continue} shows that the decay rate of $\|u(t)\|_{L^q(\mathbb{R}^d)}$ is \emph{exactly} of order $(1+t)^{-1/(2-q)}$. Hence, whenever $q\in(1,2)$, the decay rate is \emph{independent} of the initial data. \details{This is very much related to the nonlinearity of \eqref{Eqn:CDP} in that case. Indeed, when $q=2$, \eqref{Eqn:CDP} is linear and it can be shown by spectral analysis that there exists $K_1,K_2,\lambda\in\R_+$, all depending on the initial data, such that 
$$
K_1 e^{-\lambda t}\leq\|u(t)\|_{L^2(\Omega)}\leq K_2 e^{-\lambda t},\quad\text{for $t\geq 0$}.
$$}
In the fast diffusion case $2< q \leq 2_\theta ^*$, the right-hand side of \eqref{Eqn:Asymp1Continue} vanishes at the time 
\begin{equation}\label{Eqn:DefTetoile}
T_*:= \frac{q-1}{q-2}C^2\|u^0\|_{L^q(\Omega)},
\end{equation}
and the left-hand side vanishes at a time 
$$
S_*:=\frac{q-1}{q-2}R(0)^{-1}\|u^0\|_{L^q(\Omega)}.
$$
As a consequence, when $2<q\leq 2_\theta^*$, the solution $u$ of \eqref{Eqn:CDP} \emph{vanishes in finite time}, i.e.~there exists an \emph{extinction time} $t_*\in\R_+$ such that $u(t)\equiv 0$ for $t\geq t_*$. The estimate $S_*\leq t_*\leq T_*$ holds from \eqref{Eqn:Asymp1Continue}. It is worth noticing that estimate \eqref{Eqn:Asymp1Continue} is a poor estimate near the extinction time $t_*$, since a priori $t_*\neq T_*$ and $t_*\neq S_*$. On the other hand, more precise results hold about the behavior of $u$ near the extinction time. Indeed, if $2<q<2_\theta^*$, one has for any $t>0$ (\cite{bonforte_ibarrondo_ispizua, akagi_existence}),
\begin{equation}\label{Eqn:Asymp2Continue}
\left(\frac{q-2}{q-1}C^{-2}(t_*-t)\right)_+^{1/(q-2)}\leq\|u(t)\|_{L^q(\Omega)}\leq \left(\frac{q-2}{q-1}R(0)(t_*-t)\right)_+^{1/(q-2)}.
\end{equation}
Estimate \eqref{Eqn:Asymp2Continue} shows that the decay rate of $\|u(t)\|_{\lq}$ is \emph{exactly} of order $(t_*-t)^{1/(q-2)}$ when $t\rightarrow t_*$. Les us emphasize on the fact that the extinction time $t_*$ is not analytically known, although all the constants in \eqref{Eqn:Asymp1Continue} are explicit. To the best of our knowledge, the problem of finding an expression of $t_*$ in terms of the initial data and of the parameters of the problem is still an open problem, even when \eqref{Eqn:CDP} is considered with a standard Laplacian instead of a fractional one. Finally, we note also that the Rayleigh quotient \eqref{Eqn:DefRayleigh} is non-increasing (see e.g.~\cite{akagi_existence}), which, together with the Poincaré-Sobolev inequality, allows to obtain estimates similar to \eqref{Eqn:Asymp1Continue} and \eqref{Eqn:Asymp2Continue} for the Gagliardo semi-norm of $u$.

Estimates \eqref{Eqn:Asymp1Continue} and \eqref{Eqn:Asymp2Continue} can be obtained by combining the energy equality
\begin{equation}\label{Eqn:EnergyEqualityContinuous}
\frac{1}{q'}\frac{\d}{\d t}\|u(t)\|_{L^q(\Omega)}^q+[u(t)]_{H^\theta(\R^d)}^2 = 0,
\end{equation}
with the Poincaré-Sobolev inequality \eqref{Eqn:PoincareSobolevContinue} or the decay of the Rayleigh quotient. In \cite{akagi_existence} solutions to \eqref{Eqn:CDP} are constructed using the minimizing movement scheme, which, in the present case, corresponds to the semi-discrete in time implicit Euler scheme.
Moreover, the limiting procedure as the time step goes to zero is performed by using a discrete version of \eqref{Eqn:EnergyEqualityContinuous}. This paves the way to construct a \emph{fully discrete implicit Euler scheme} that is presented in this paper. Moreover, we show that it satisfies a fully discrete version of the identity \eqref{Eqn:EnergyEqualityContinuous}. Then our main result is a discrete version of the estimate \eqref{Eqn:Asymp1Continue}, and of the extinction phenomenon in the fast diffusion case. More precisely, we show that \eqref{Eqn:Asymp1Continue} holds for the fully discrete scheme, up to an error term which converges to zero as the time step converges to zero. Despite numerical solutions cannot vanish for non-zero initial data, in the fast diffusion case we show an \emph{almost extinction phenomenon}. It means that a very fast decay of the discrete solution past the extinction time is proved. Finally, the numerical scheme allows to introduce a method to estimate the extinction time in the fast diffusion case, and to investigate numerically the large time asymptotics of the equation. We will observe numerically the convergence of rescaled solutions to asymptotic profiles in the fast diffusion case, and determine numerically the rate of convergence.

\textit{Related work.~}A very large number of contributions have been made to the non-fractional porous medium equation and fast diffusion equation (see e.g.~\cite{vazquez_pme},\cite{bonforte_figalli_survey} and references therein). Various numerical methods have also been introduced for the non-fractional porous medium equation, including finite differences \cite{graveleau_jamet}, finite volumes \cite{bessemoulin_filbet,eymard_gallouit_herbin_michel}, finite elements \cite{baines_moving_2005,baines_scale-invariant_2006,ngo_huang,zhang_numerical_2009}, interface tracking algorithms \cite{di_benedetto_interface_1984, liu_2018,monsaingeon_explicit_2016}, Lagragian schemes \cite{carrillo_lagrangian_2018, liu_lagrangian_2020}, relaxation schemes \cite{cavalli_high_order_2007, naldi_relaxation_2002}, and schemes based on the gradient flow formulation of the equation with respect to the Wasserstein distance \cite{carrillo_numerical_2017,carrillo_numerical_2016}. 
At the continuous level, a large number of contributions have also been made to the fractional porous medium and fast diffusion equation on $\R^d$ (see e.g.~\cite{de_pablo_fractional_2011,de_pablo_general_2012,vazquez_barenblatt_2014,dePablo_singular, bonforte_vazquez_2014, bonforte_endal}), while the theory on bounded domain is less complete and still growing (see e.g. \cite{bonforte_vazquez, bonforte_sire_vazquez, bonforte_figalli_vazquez, franzina_volzone} for the porous medium case and \cite{bonforte_ibarrondo_ispizua,akagi_existence} for the fast diffusion case). At the discrete level, several methods have been introduced to discretize the fractional Laplacian, and solve the fractional Poisson equation \cite{delia_survey}. For example, the fractional Laplacian on the whole domain, or the restricted fractional Laplacian, have been descretized using a finite-differences quadrature approach \cite{huang_oberman}, finite elements \cite{delia_fractional_2013, acosta_fractional_2017,antil_analysis_2022}, or a semigroup approach \cite{ciaurri_nonlocal_2018}. In this paper we will mainly draw on the finite-differences quadrature methods \cite{huang_oberman} in order to discretize the spatial part of \eqref{Eqn:CDP}, but we will also make use of the discrete functional analysis tools first introduced in \cite{ciaurri_nonlocal_2018}. Numerical methods have also been proposed for the spectral Laplacian on bounded domain, e.g.~\cite{cusimano_discretizations_2018, cusimano_numerical, bonito_borthagaray_nochetto_otarola_salgado, nochetto_otarola_salgado,andrea_pasciak}. Then, several methods have been studied for nonlinear and nonlocal equations on $\R^d$ \cite{del_teso_finite_2014, del_teso_convergent_2023, delteso2013finitedifferencemethodgeneral, delTeso_endal_jakobsen_part1, del_teso_robust_2018}. Notably, the fractional porous medium equation and fractional fast diffusion equation on $\R^d$ can be handled by the methods in \cite{delTeso_endal_jakobsen_part1, del_teso_robust_2018}. However, all this methods deal with the problem on whole domain $\R^d$.  To our knowledge, the only work which deals with the problem on bounded domain is \cite{carrillo_finite_2024}, where a spectral fractional Laplacian, and only the porous medium range is covered. 

\textit{Outline.~}This paper is composed of six sections. In the next section we introduce the numerical scheme and state the main results. In section 3, main discrete functional analysis tools used to study the numerical scheme are introduced. Section 4 is dedicated to the proof of well-posedness and weak convergence results of the scheme, and Section 5 to the proof of the discrete decay estimates and extinction phenomenon. Finally, Section 6 is dedicated to the numerical results, and in particular to the computation of the extinction time and large time asymptotics. Eventually, a technical discrete Poincaré-Sobolev inequality is proved in appendix.

\textbf{Acknowledgement} The authors would like to express their gratitude to Frédéric Lagoutière for the many insightful discussions we had about this project, and to Goro Akagi for his valuable suggestions and advice regarding the near extinction asymptotics of the fast diffusion equation.

\section{Presentation of the schemes and main results}\label{Sec:Results}

Throughout this article, the domain $\Omega$ is an open bounded interval $I=(-L,L)$ (non-symmetric domains can be reduced to this case by translation). It is discretized with a uniform space step $h$, chosen in such a way that there exists $M_x\in \mathbb{N}\setminus\{0\}$ such that $L=(M_x+1)h$, and the space nodes are defined by $(x_i = ih)_{i\in\mathbb{Z}}$. Time is discretized with a regular subdivision $(t_n=n\Delta t)_{n\in\N}$ where $\Delta t>0$ denotes the time step. 

\subsection{Presentation of the schemes and notations}
Consider the following fully discrete scheme for problem \eqref{Eqn:CDP},
\begin{equation}\label{Eqn:FD}\tag{{\sf FD}}
\left\{
\begin{aligned}
\frac{(u_i^{n+1})^{q-1}-(u_i^{n})^{q-1}}{\Delta t}+\left[(-\Delta)^\theta_h \bu^{n+1}\right]_i &= 0, && \text{for $|i|\leq M_x$ and $n\geq 0$},\\
u^n_i &= 0, && \text{for $|i|\geq M_x+1$ and $n\geq 0$},\\
u_i^0 &= (\bu^0)_i && \text{for $|i|\leq M_x$},
\end{aligned}
\right.
\end{equation}
for some initial data $\bu^0\in\R^{\Z}$ such that $(\bu^0)_i = 0$ for $|i|>M_x$.  For $n\geq 0$ and $i\in\Z$, $u^n_i$ stands for an approximation of $u(t_n,x_i)$, the solution of \eqref{Eqn:CDP} evaluated at time $t_n$ and point $x_i$, and we denote $\bu^n:=(u^n_i)_{i\in\Z}$. Throughout this article, bold letters will denote vectors of $\R^\Z$. Consider also the semi-discrete scheme,
\begin{equation}\label{Eqn:SD}\tag{{\sf SD}}
\left\{
\begin{aligned}
\frac{{\rm d}}{{\rm d} t} v_i^{q-1}(t) + \left[(-\Delta)^s_h \bv(t)\right]_i &= 0,&& \text{for $|i|\leq M_x$, and $t>0$},\\
v_i(t)&=0 && \text{for $|i|\geq M_x+1$, and $t\geq 0$},\\
v_i(0)&=(\bu^0)_i && \text{for $|i|\leq M_x$}.\end{aligned}
\right.
\end{equation}
where for any $i\in\Z$, $v_i(t)$ is a function from $\R^+$ to $\R$ such that $v_i^{q-1}\in W^{1,1}(0,T;\R)$ for any $T>0$, and $\bu^0$ is the same initial data as \eqref{Eqn:FD}. Again, $v_i(t)$ is an approximation of $u(t,x_i)$, and $\bv(t)$ denotes $(v_i(t))_{i\in\Z}$.

In \eqref{Eqn:FD} and \eqref{Eqn:SD}, $(-\Delta)^\theta_h$ denotes the discretization of the fractional Laplacian. In what follows, we choose the discretization introduced by Huang and Oberman in \cite{huang_oberman}. It is defined for all $\mathbf{f}:=(f_i)_{i\in\Z}$ decaying fast enough at infinity by,
\begin{equation}\label{Eqn:DefDiscreteFracLap}
\forall i\in\Z,\quad\left[(-\Delta)_h^\theta \mathbf{f}\right]_i := \sum_{j\in \Z} w_j^h(f_i-f_{i-j}),
\end{equation}
with symmetric weights $(w_j^h)_{j\in\Z}$ ($w_j^h = w_{-j}^h$), defined by
\begin{equation}\label{Eqn:DefinitionWeights}
w_j^h=h^{-2\theta}\begin{cases}
    \displaystyle
\frac{C_{1,\theta}}{2-2\theta} - G''(1) - \frac{G'(3)+3G'(1)}{2} + G(3) - G(1), & j=\pm 1,\\
\displaystyle
2[G'(j+1) + G'(j-1) - G(j+1) + G(j-1)], & j=\pm 2, \pm 4,\dots \\
\displaystyle
-\frac{G'(j+2)+6G'(j)+G'(j-2)}{2} + G(j-2) - G(j-2), & j=\pm 3, \pm 5,\dots,
\end{cases},
\end{equation}
where
$$
G(t)=
\begin{cases}
    \displaystyle
\frac{C_{1,\theta}}{(2-2\theta)(2\theta-1)2\theta}|t|^{2-2\theta}, & 2\theta\neq 1,\\
\displaystyle
C_{1,\theta}(t-t\log |t|), & 2\theta=1.
\end{cases}
$$
It is worth noticing that the weights depends on $\theta$. We generally omit this dependence in the notation for the sake of readability, but we may ocasionnaly write $w_j^h(\theta)$ to emphasize it. Two properties of this discrete fractional Laplacian will be used for the analysis of the schemes \eqref{Eqn:FD} and \eqref{Eqn:SD}. The first one is the following estimate on the weights, proved in \cite{ayi_herda_hivert_tristani}, that shows that $(-\Delta)_h^\theta$ has a discrete convolution structure matching the definition of the continuous fractional Laplacian \eqref{Eqn:DefFracLap}.
\begin{prop}\label{Prop:EstimatesWeightsDiscreteFracLap}
There exists positive constants $b_\theta$ and $B_\theta$ depending only on $\theta\in(0,1)$ such that
$$
\frac{b_\theta}{|jh|^{1+2\theta}}h\leq w_j^h(\theta) \leq \frac{B_\theta}{|jh|^{1+2\theta}}h.
$$
\end{prop}
\begin{rem}\label{Rem:SumWeights} The sum of the weights can also be explicitly computed (see \cite{huang_oberman}):
\begin{equation}\label{Eqn:SumWeights}
\sum_{j\in\Z} w_j^h(\theta) = \frac{2^{2s} \Gamma((2s+1)/2)}{\pi^{1/2}\Gamma(2-2s/2)} h^{-2\theta},
\end{equation}
where $\Gamma$ is the Euler function.
\end{rem}
\noindent
The second property of $(-\Delta)^\theta$ to be used is the following error estimate for smooth functions, proved in \cite{huang_oberman}.
\begin{thm}[Error estimate]\label{Thm:ErreurDiscreteFracLap}
Let $f\in C^4(\R)$ such that $f^{(3)}$ and $f^{(4)}$ are bounded. 
Then,
$$
\sup_{i\in\Z}\left|\left[(-\Delta)^\theta_h f\right]_i -(-\Delta )^\theta f(x_i)\right| =\mathcal{O}(h^{3-2\theta}),
$$
where $(-\Delta)^\theta_h f$ is the evaluation of the discrete fractional Laplacian to $(f(x_i))_{i\in\Z}$.
\end{thm}

The analysis of the schemes \eqref{Eqn:FD} and \eqref{Eqn:SD} will also use discrete functional spaces, which play the same role as the functional spaces in the analysis of \eqref{Eqn:CDP} in \cite{akagi_existence}. They are defined hereafter. In the following notation, norms are denoted with double-bars $\|\cdot\|$, and semi-norm are denoted with brackets $[\cdot]$. Some properties of those spaces will be given in Section \ref{Sec:DiscreteFunctionalAnalysis}. 

\begin{dfn}\label{Dfn:DefDiscreteSpaces}
Let $0<\theta<1$, and $1\leq q<\infty$.
\begin{itemize}
\item The discrete Lebesgue space of order $q$ is defined by 
\begin{gather*}
\ell^q_h(\R):=\{\bof=(f_i)_{i\in\Z}\in \R^\Z:\:\|\bof\|_{\lq}<+\infty\},
\text{ with }\|\bof\|_{\ell_h^q(\R)}:=\left(\sum_{i\in\Z} h|f_i|^q\right)^{1/q}.
\end{gather*}
\item The discrete fractional Sobolev space of order $2$ and exponent $\theta$ is defined by
\begin{gather*}
H^\theta_h(\R):=\{\bof\in \ell^2_h(\R):\:[\bof]_{H_h^\theta(\R)}<+\infty\},\text{ with }
[\bof]_{H_h^\theta(\R)}:=\left(\frac{1}{2}\sum_{i,j\in\Z} h w^h_j(\theta)|f_i-f_{i-j}|^2\right)^{1/2}.
\end{gather*}
It is an Hilbert space for the scalar product, defined for any $\bof,\bog\in H^\theta_h(\R)$,
$$
(\bof,\bog)_{H^\theta_h(\R)}:=\sum_{i} h f_i g_i + \sum_{i,j\in\Z} h^2w^h_j(\theta)(f_{i}-f_{i-j})(g_{i}-g_{i-j}).
$$
\item The subspace of $H^\theta_h(\R)$ of sequences vanishing outside the domain $\Omega$ is denoted by
$$
\X := \{\bof = (f_i)_{i\in\mathbb{Z}}\in\mathbb{R}^\mathbb{Z}:\: f_i=0\text{ for }|i|>M_x,\,[\bof]_{H^\theta_h(\mathbb{R})}<+\infty\}.
$$
It will be shown in Section \ref{Sec:DiscreteFunctionalAnalysis} that the semi-norm $[\cdot]_{H^\theta_h(\mathbb{R})}$ defines a proper norm on $\X$. Therefore we will denote $\|\bof\|_{\X}:=[\bof]_{H^\theta_h(\R)}$ and 
$$
(\bof,\bog)_{\X}:=\sum_{i,j\in\Z} h^2w^h_j(\theta)(f_{i}-f_{i-j})(g_{i}-g_{i-j}) \quad\text{for }\bof,\bog\in\X.
$$
\end{itemize}
\end{dfn}

Important tools used in \cite{akagi_existence} to prove estimates \eqref{Eqn:Asymp1Continue} and \eqref{Eqn:Asymp2Continue}, i.e.~decay properties of the solution of \eqref{Eqn:CDP}, are a Green formula for the fractional Laplacian, and the Poincaré-Sobolev inequality \eqref{Eqn:PoincareSobolevContinue}. It is worth noticing that, these two formulas have discrete equivalents using spaces introduced in Definition \ref{Dfn:DefDiscreteSpaces}.

\begin{lem}\label{Lem:DiscretePoincareSobolev}
    For any $q\in (1,2_\theta^*]$ there exists $C>0$ such that for any $h$ and any $\bof\in \X$,
    \begin{equation}\label{Eqn:DiscretePoincareSobolev}
    \|\bof\|_{\lq}\leq C \|\bof\|_{\Xh}.
    \end{equation}
    Moreover, for any $\bof,\bog \in \X$,
    \begin{equation}\label{Eqn:GreenFormulaDiscrete}
    ((-\Delta)_h^\theta \bof,\bog)_{\ell^2_h(\R)}=(\bof,\bog)_{\X}.
    \end{equation}
\end{lem}
Lemma \ref{Lem:DiscretePoincareSobolev} is proved in Section \ref{Sec:DiscreteFunctionalAnalysis}, following discrete Sobolev and Poincaré inequalities, first introduced in \cite{ciaurri_nonlocal_2018}\details{ (a proof is also provided in Appendix \ref{Appendix:DiscretePoincareSobolev})}. It mainly relies on the fact that the discrete fractional Laplacian has a discrete convolution structure (see \eqref{Eqn:DefDiscreteFracLap}) with discrete weights matching the continuous weights (see Proposition \ref{Prop:EstimatesWeightsDiscreteFracLap}). It is the main tool to obtain well-posedness of the schemes \eqref{Eqn:FD} and \eqref{Eqn:SD}, and decay estimates. 

\subsection{Well-posedness of the schemes and convergence}
Our first set of results is about the well-posedness and convergence of the schemes \eqref{Eqn:FD} and \eqref{Eqn:SD}. Proposition \ref{Prop:CvgFDtoSDExistenceStabilityEnergyIneqSD} deals with the scheme \eqref{Eqn:FD}, and Proposition \eqref{Prop:EstimatesSD} deals with the scheme \eqref{Eqn:SD} and the convergence of \eqref{Eqn:FD} to \eqref{Eqn:SD}. 

\begin{prop}[Existence and stability for \eqref{Eqn:FD}]\label{Prop:ExistenceStabilityEnergyIneqFD}
Problem \eqref{Eqn:FD} admits a unique solution $\bu = (u^n_i)_{n\in\N,i\in\Z}$. Moreover, it satisfies the stability estimates
\begin{gather}\label{Eqn:EnergyIneqFD1}
-\|\bu^n\|_{\X}^2\leq \frac{1}{q'}\frac{\|\bu^{n+1}\|_{\ell^q_h(\R)}^q - \|\bu^{n}\|_{\ell^q_h(\R)}^q}{\Delta t}  \leq -\|\bu^{n+1}\|_{\X}^2,\quad\text{for any $n\in\N$,}\\
\label{Eqn:EnergyIneqFD2}
\frac{4}{qq'}\left\| \frac{(\bu^{n+1})^{q/2}-(\bu^{n})^{q/2}}{\Delta t} \right\|_{\ell^2_h(\R)}^2+\frac{1}{2}\frac{\|\bu^{n+1}\|_{\X}^2-\|\bu^{n}\|_{\X}^2}{\Delta t} \leq 0,\quad\text{for any $n\in\N$},
\end{gather}
where $(\bu^{n})^{q/2}:=((u^n_i)^{q/2})_{i\in\Z}\in\R^\Z$.
\end{prop}

\begin{prop}[Existence and stability for \eqref{Eqn:SD}]\label{Prop:CvgFDtoSDExistenceStabilityEnergyIneqSD}
Let $T\in \R_+^*$. Problem \eqref{Eqn:SD} admits a unique solution $\bv: \R_+ \ni t \mapsto (v_i(t))_{i\in\Z}$. Moreover for any $i\in\Z$, $v_i^{q-1}\in W^{1,1}(0,T;\R)$, $v_i^q\in W^{1,1}(0,T;\R)$, $v_i^{q/2}\in W^{1,2}(0,T;\R)$, and it satisfies
\begin{gather}
\label{Eqn:EnergyIneqSD1}
			\frac{1}{q'}\frac{\rm{d}}{{\rm{d}}t}\|\bv(t)\|_{\ell^q_h(\R)}^q+\|\bv(t)\|_{\Xh}^2=0,\quad \text{for a.e.~}t\in(0,\infty),\\
\label{Eqn:EnergyIneqSD2}
			\frac{4}{qq'}\int_s^t \|\partial_\tau \bv^{q/2}(\tau)\|_{\ell^2_h(\R)}^2{\rm d}\tau + \frac{1}{2}\|\bv(t)\|_{\Xh}^2\leq\frac{1}{2}\|\bv(s)\|_{\Xh}^2, \quad\text{for $0\leq s < t <\infty$}.
\end{gather}
Furthermore, let $\bu=(u_i^n)_{n\in\N,i\in\Z}$ be the solution of Problem \eqref{Eqn:FD}, and define for all $i\in\Z$, $v^{\Delta t}_i:\R^+\rightarrow \R$ such that $v^{\Delta t}_i(t)=u_i^{n+1}$ when $t\in(n\Delta t,(n+1)\Delta t]$. Then, $v^{\Delta t}_i\rightarrow v_i$ uniformly on $[0,T]$ and for every $i\in\Z$, as $\Delta t\rightarrow 0$. 
\end{prop}

As a by-product of Proposition \ref{Prop:CvgFDtoSDExistenceStabilityEnergyIneqSD}, the fully discrete scheme \eqref{Eqn:FD} converges to the semi-discrete scheme \eqref{Eqn:SD} as the time-step goes to zero. Convergence of the fully discrete scheme \eqref{Eqn:FD} to the continuous problem \eqref{Eqn:CDP} is proved under Lax-Wendroff type hypothesis: if solutions to \eqref{Eqn:FD} converge almost everywhere to a function $u$ as the time step goes to zero, and if it is bounded in $L^\infty$ norm, then $u$ must be a weak solution to \eqref{Eqn:CDP}.

\begin{prop}[Convergence of \eqref{Eqn:FD}]\label{Prop:LaxWendroff}
Let $\bu = (u_i^n)_{n\in \N, i\in\Z}$ the solution of \eqref{Eqn:FD} with some initial condition $\bu^0:=(u^0_i)_{i\in\Z}$ such that $u^0_i=0$ for $|i|>M_x$. Let $u_{\Delta t,h}:\R^+\times\R\rightarrow \R$ such that $u_{\Delta t,h}(t,x) = u_i^n$ for $(t,x)\in[t_n,t_{n+1})\times [x_i-h/2,x_{i}+h/2)$, and $u^0_h:\R\rightarrow\R$ such that $u^0_{h}(x) = u^0_i$ for $x\in[x_i-h/2,x_{i}+h/2)$. Assume that $\|u_{\Delta t,h}\|_{L^\infty((0,T)\times\R)}$ is bounded uniformly in $\Delta t$ and $h$ for every $T>0$, and that there exists $u\in L^\infty((0,\infty)\times\R)$ such that $u_{\Delta t,h}\rightarrow u$ almost everywhere when $|(\Delta t, h)|\rightarrow 0$. Assume also that there exists $u^0\in L^1(I)$ such that $u^0_{h}\rightarrow u^0 $ in $L^1(I)$ as $h\rightarrow 0$. Then $u$ is a weak solution of \eqref{Eqn:CDP}, in the sense that, for all $T>0$ and for all test function $\varphi \in C_c^\infty([0,T)\times I)$ it satisfies
\begin{multline}\label{Eqn:EqFaible}
-\int_0^T\int_I u^{q-1}(t,x)\partial_t\varphi(t,x){\rm d}x {\rm d}t+\int_0^T\int_I u(t,x)(-\Delta)^\theta \varphi(t,x){\rm d}x {\rm d}t \\ - \int_I (u^0)^{q-1}(x)\varphi(0,x){\rm d}x=0.
\end{multline}
\end{prop}
Proposition \ref{Prop:ExistenceStabilityEnergyIneqFD}, \ref{Prop:CvgFDtoSDExistenceStabilityEnergyIneqSD} and \ref{Prop:LaxWendroff} are proved in section \ref{Sec:WellPosednessConvergence}.

\subsection{Decay estimates of the schemes} After looking at the well-posedness and convergence of the schemes \eqref{Eqn:FD} and \eqref{Eqn:SD}, it is natural to examine their qualitative behaviour. The next two theorems show that \eqref{Eqn:FD} and \eqref{Eqn:SD} satisfy decay estimates close to \eqref{Eqn:Asymp1Continue}. 
Proposition \ref{Prop:EstimatesSD} deals with the semi-discrete case and Proposition \ref{Prop:EstimatesFD} deals with the fully-discrete case. 
\begin{prop}[Decay estimates for \eqref{Eqn:SD}]\label{Prop:EstimatesSD}
Assume $q\in (1,2_\theta^*]$ and $\bu^0\not\equiv 0$. Let $C>0$ be any constant satisfying the discrete Poincaré-Sobolev inequality \eqref{Eqn:DiscretePoincareSobolev}, and let $\bv: \R_+ \ni t \mapsto (v_i(t))_{i\in\Z}$ be the solution to \eqref{Eqn:SD} . Then,
\begin{enumerate}
\item if $q\in(1,2_\theta^*]\setminus\{2\}$, then for any $t>0$,
\begin{equation}\label{Eqn:Asymp1DiscreteSD}
\left(\|\bu^0\|_{\lq}^{-(2-q)}+\frac{2-q}{q-1}R(\bu^0)t\right)_+^{\frac{-1}{2-q}}\leq \|\bv(t)\|_{\lq}\leq \left(\|\bu^0\|_{\lq}^{-(2-q)}+\frac{2-q}{q-1}C^{-2}t\right)_+^{\frac{-1}{2-q}},
\end{equation}
\item if $2< q\leq 2_\theta^*$, then there exists $t_*^h>0$ such that for any $t>0$,
\begin{equation}\label{Eqn:Asymp2DiscreteSD}
\left(\frac{q-2}{q-1}C^{-2}(t_*^h-t)\right)_+^{1/(q-2)}\leq\|\bv(t)\|_{\lq}\leq \left(\frac{q-2}{q-1}R(0)(t_*^h-t)\right)_+^{1/(q-2)},
\end{equation}
\item the discrete Rayleigh quotient 
\begin{equation}\label{Eqn:DecayRayleigh}
R(t):=\frac{[\bv(t)]^2_{H_h^\theta(\R)}}{\|\bv(t)\|^2_{l_h^q(\R)}}\text{ is non-increasing}\left\{
\begin{aligned}
&\text{on $[0,\infty)$ if $q\in (1,2]$},\\ 
&\text{on $[0,t_*^h)$ if $q\in (2,2_\theta^*]$}.
\end{aligned}\right.
\end{equation}
\end{enumerate}
\end{prop}
Estimates \eqref{Eqn:Asymp1DiscreteSD},\eqref{Eqn:Asymp2DiscreteSD} are close to the continuous estimates \eqref{Eqn:Asymp1Continue},\eqref{Eqn:Asymp2Continue}, with only discrete norms replacing continuous norms. The proof of discrete estimates works in the same way as that of continuous estimates. Estimate \eqref{Eqn:Asymp1DiscreteSD} converges to the continuous estimate \eqref{Eqn:Asymp1Continue} as $h\rightarrow 0$, under suitable conditions on the initial data $\bu^0$. However, in the fast diffusion case, the extinction time $t_*^h$ of the semi-discrete scheme may depend on $h$. We could not prove that $t_*^h$ converge to the extinction time $t_*$ of the continuous problem \eqref{Eqn:CDP}, which prevents from passing to the limit in \eqref{Eqn:Asymp2DiscreteSD}.

\begin{prop}[Decay estimates for \eqref{Eqn:FD}]\label{Prop:EstimatesFD}
    Assume $q\in (1,2_\theta^*]$, and $\bu^0=\bu^0_h$ is such that $\bu^0_h\not\equiv 0$ and $\|\bu^0_h\|_{\lq}$ converges as $h\rightarrow 0$. Let $C>0$ be any constant satisfying the discrete Poincaré-Sobolev inequality \eqref{Eqn:DiscretePoincareSobolev}, and let $\bu=(\bu_i^n)_{i\in\Z,n\in\N}$ be the solution to \eqref{Eqn:FD} . Then,
    \begin{enumerate}
    \item  If $q\in (1,2_\theta^*]\setminus\{2\}$, for all $\varepsilon > 0$ there exists $\Delta t_0$ such that for any $h$ and any $\Delta t < \Delta t_0$,
    \begin{equation}\label{Eqn:Asymp1DiscreteFD}
    \|\bu^n\|_{\lq}\leq \left(\|\bu^0\|_{\lq}^{-(2-q)}+\frac{2-q}{q-1}C^{-2}n\Delta t\right)_+^{\frac{-1}{2-q}} + \varepsilon.
    \end{equation}
    In particular, the right-hand side of \eqref{Eqn:Asymp1DiscreteFD} converges to the right-hand side of \eqref{Eqn:Asymp1Continue} as $(\Delta t,h)\rightarrow 0$.
    \item If $q\in(2,2_\theta^*]$, let 
    \begin{equation}\label{Eqn:DefT*h}
    T_*^h:=\frac{q-1}{q-2}C^2\|\bu^0\|_{\lq}^{q-2}.
    \end{equation}
    Then,
    \begin{equation}\label{Eqn:ExtinctionEstimateFD}
    \|\bu^{n}\|_{\lq}\leq \|\bu^0\|_{\lq}\left( \frac{T_*^h}{n\Delta t }\right)^{n/2},\quad \text{for any $n\geq 0$}.
    \end{equation}
    \end{enumerate}
\end{prop}

\noindent
Propositions \ref{Prop:EstimatesFD} and \ref{Prop:EstimatesSD} will be proved in Section \ref{Sec:DecayEstimates}. Estimate \eqref{Eqn:Asymp1DiscreteFD} is similar to the upper bound \eqref{Eqn:Asymp1Continue} up to a correction term $\varepsilon$ that can be computed explicitly. Moreover, it is shown to converge uniformly to $0$ as the time-step converges to $0$, and \emph{independently of $h$}. It will be shown in the proof that the error term $\varepsilon$ is smaller than $\mathcal{O}(\Delta t_0)$ in case $q\in(1,2)$, and smaller than $\mathcal{O}(\Delta t_0^{1/2})$ in case $q\in(2,2_\theta^*]$. The fact that the error terms are independent of $h$ is an important feature of this result. As a consequence, the right-hand side in \eqref{Eqn:Asymp1DiscreteFD} converges uniformly to the continuous upper bound in \eqref{Eqn:Asymp1Continue} as $(\Delta t,h)\rightarrow 0$. On the other hand, combining the convergence result of \eqref{Eqn:FD} to \eqref{Eqn:SD} (Proposition \ref{Prop:CvgFDtoSDExistenceStabilityEnergyIneqSD}) and estimates of Proposition \ref{Prop:EstimatesSD} it is possible to obtain an estimate similar to \eqref{Eqn:Asymp1DiscreteFD} but with an error term depending on $h$ which prevents from passing to the limit as $h\rightarrow 0$. 
\details{
More precisely, let us fix a space step $h>0$, and let $(u^n_i)_{i\in\Z,n\in\N}$ be the solution of the fully discrete scheme \eqref{Eqn:FD}. The uniform convergence result in Proposition \ref{Prop:CvgFDtoSDExistenceStabilityEnergyIneqSD} and estimates \eqref{Eqn:Asymp1DiscreteSD} and \eqref{Eqn:Asymp2DiscreteSD} yield the existence of error terms $(e_n^{\Delta t,h})_{n\geq 0}$, uniformly converging to zero as $\Delta t\rightarrow 0$, such that
\begin{multline}\label{Eqn:MauvaisEst1}
    \left(\|u_0\|_{\lq}^{-(2-q)}+\frac{2-q}{q-1}R(u_0)t\right)_+^{\frac{-1}{2-q}} - e_n^{\Delta t,h}\\
    \leq \|u^n\|_{\lq}\leq\\
    \left(\|u_0\|_{\lq}^{-(2-q)}+\frac{2-q}{q-1}C^{-2}t\right)_+^{\frac{-1}{2-q} } + e_n^{\Delta t,h},
\end{multline}
In the fast diffusion case, denoting $t_*^h$ the extinction time of the semi discrete scheme \eqref{Eqn:SD}, we also have
\begin{equation}\label{Eqn:MauvaisEst2}
    \left(\frac{q-2}{q-1}C^{-2}(t_*^h-t)\right)_+^{1/(q-2)}- e_n^{\Delta t,h}\leq\|u^n\|_{\lq}\leq \left(\frac{q-2}{q-1}R(0)(t_*^h-t)\right)_+^{1/(q-2)} + e_n^{\Delta t,h}.
\end{equation}
However, Proposition \ref{Prop:CvgFDtoSDExistenceStabilityEnergyIneqSD} and Proposition \ref{Prop:EstimatesSD} tell us nothing about the behavior of $e_n^{\Delta t,h}$ and $t_*^h$ with respect to $h$, so it is not possible to pass to the limit in \eqref{Eqn:MauvaisEst1},\eqref{Eqn:MauvaisEst2} as $\Delta t,h\rightarrow 0$ and recover the continuous estimates \eqref{Eqn:Asymp1Continue},\eqref{Eqn:Asymp2Continue}.
}

Finally, it should be noted that although estimate \eqref{Eqn:Asymp1DiscreteFD} allows to recover finite time extinction when $\Delta t$ converges to 0, it does not tell whether the solution $(u^n_i)_{i\in\Z,n\in\N}$ of the fully discrete scheme extinguishes or not when $\Delta t$ is fixed. Actually, if $\bu^{n+1}=0$, then \eqref{Eqn:FD} implies $\bu^n=0$, and by induction we obtain $\bu^0=0$, which means that the only solution to \eqref{Eqn:FD} that vanishes in finite time is the one starting from a zero initial condition. However \eqref{Eqn:ExtinctionEstimateFD} shows that the decay of $\|\bu^n\|_{\lq}$ is very fast (more than geometric) on every time node above the time $T_*^h$. Moreover, $T_*^h$ is the discrete equivalent of the time $T_*$, defined in \eqref{Eqn:DefTetoile}, which is an upper bound for the extinction time in the continuous case. 

\section{Discrete Functional Analysis}\label{Sec:DiscreteFunctionalAnalysis}
The goal of this section is to prove Lemma \ref{Lem:DiscretePoincareSobolev}. This lemma links the discrete functional spaces introduced in Definition \ref{Dfn:DefDiscreteSpaces} and the discrete fractional Laplacian. Moreover, it introduces the discrete Poincaré-Sobolev inequality \eqref{Eqn:DiscretePoincareSobolev}, which will be an important tool in the analysis of the schemes. We start with the proof of the discrete Poincaré-Sobolev inequality. As its name suggests, it is the combination of a discrete Poincaré inequality, and a discrete Sobolev inequality, which are given in the next two lemmas.
\begin{lem}[Discrete Sobolev inequality, see {\cite[Theorem 6.4]{ciaurri_nonlocal_2018}}]\label{Lem:DiscreteSobolevIneq}
Let $\theta\in(0,1)$, and let $2_\theta^*= 2/(1-2\theta)_+ \in(2,\infty]$. There exists a constant $C=C(\theta)>0$, independant of $h$, such that 
\begin{enumerate}
\item if $2\theta<1$, then for any $\bof\in H^{\theta}_h(\R)$, $\|\bof\|_{\ell^{2_\theta^*}_h(\R)}\leq C [\bof]_{H^{\theta}_h(\R)}$,
\item if $2\theta\geq 1$, then for any $q\in[2,\infty)$ and $u\in H^{\theta}_h(\R)$, $\|\bof\|_{\ell^q_h(\R)}\leq C\|\bof\|_{H^{\theta}_h(\R)}$.
\end{enumerate}
\end{lem}
\begin{lem}[Discrete Poincaré inequality, see {\cite[Theorem 6.4]{ciaurri_nonlocal_2018}}]\label{Lem:DiscretePoincareIneq}
Let $h<1, \theta\in(0,1)$, and $I=(-L,L)$ with $L=(M_x+1)h$. There exists a positive constant $C=C(L,\theta)$ depending only on $L$ and $\theta$ but not $h$, such that
\begin{equation}\label{Eqn:DiscretePoincareIneq}
\|\bof\|_{\ell^2_h(\R)}\leq C \|\bof\|_{\Xh},\quad\forall\ \bof\in\Xh.
\end{equation}
\end{lem}
As a consequence of Lemma \ref{Lem:DiscretePoincareIneq}, $\|\cdot\|_{\Xh}:=[\cdot]_{H^\theta_h}(\R)$ defines a norm on $\Xh$, which justifies the notation with a double-bar. Moreover, $\|\cdot\|_{\Xh}$ is equivalent to the discrete fractional Sobolev norm $\|\cdot\|_{H^\theta_h}(\R)$ in $\X$. The proofs of Lemma \ref{Lem:DiscreteSobolevIneq} and \ref{Lem:DiscretePoincareIneq} are given in \cite[Theorem 6.4]{ciaurri_nonlocal_2018} (see also \cite[Theorem 1.1]{ciaurri_nonlocal_2018}), and follow the proof for the continuous Sobolev and Poincaré inequalities of \cite{savin_valdinoci}. \details{These lemmas are also proved in Appendix \ref{Appendix:DiscretePoincareSobolev}, where they are stated for discrete fractional Sobolev spaces with general exponents}. By combining inequality \eqref{Eqn:DiscretePoincareIneq} and Lemma \ref{Lem:DiscreteSobolevIneq}, we obtain the discrete Poincaré-Sobolev inequality \eqref{Eqn:DiscretePoincareSobolev} as announced.
\begin{rem}
Let us denote $\ell^2_h(I):=\{\bof\in \ell^2_h(\mathbb{R}):\:\bof_i=0,\,i>M_x\}$. Since the weights $w_i^h$ are summable (see Remark \ref{Rem:SumWeights}), any sequence vanishing outside $I$ is actually in $\Xh$ and, in fact, in $\ell^2_h(I)$. Note that these two discrete functional spaces coincide for fixed $h>0$. However, the two associated norms are not equivalent \emph{uniformly} with respect to $h$, though they are at fixed $h$. The point of Lemma \ref{Lem:DiscretePoincareIneq} is to show that the injection $\Xh\hookrightarrow\ell_h^2(I)$ is continuous \emph{uniformly} in $h$. Of course, this remark also holds for other discrete functional spaces. We emphasize on the fact that the decay estimates we prove in this paper hold uniformly in $h$, meaning that all the estimates, inequalities, and embeddings we deal with are true uniformly in the discretization parameters.
\end{rem}

The discrete fractional Sobolev space $\Xh$ and the discrete fractional Laplacian $(-\Delta)_h^\theta$ are related by formula \eqref{Eqn:GreenFormulaDiscrete}, which is proved below.

\begin{prf}{Lemma \ref{Lem:DiscretePoincareSobolev}}
Let $\bof,\bog\in\Xh$. Note that $(-\Delta)^\theta_h \bof$ is well-defined since $\bof$ has a bounded support and the weights $(w_j^{h})$ are summable. Moreover, $((-\Delta)^\theta_h \bof,\bog)_{\ell^2_h(\R)}$ has also a meaning since $\bog$ has bounded support. One has, thanks to the symmetry of the weights $w_j^{h}=w_{-j}^h$ (see \eqref{Eqn:DefinitionWeights}),
$$
\begin{aligned}
\sum_{i\in\Z}\sum_{j\in\Z} h w_j^{h}(f_i-g_{i-j})g_i=& \frac{1}{2}\sum_{\substack{i\in\Z\\j\in\Z}} h w_j^{h}(f_i-f_{i-j})g_i + \frac{1}{2}\sum_{\substack{i\in\Z\\j\in\Z}} h w_j^{h}(f_i-f_{i+j})g_i\\
=& \frac{1}{2}\sum_{\substack{i\in\Z\\j\in\Z}} h w_j^{h}(f_i-f_{i-j})g_i + \frac{1}{2}\sum_{\substack{i\in\Z\\j\in\Z}} h w_j^{h}(f_{i-j}-f_{i})g_{i-j}\\
=& \frac{1}{2} \sum_{\substack{i\in\Z\\j\in\Z}} h w_j^{h}(f_i-f_{i-j})(g_i-g_{i-j}),
\end{aligned}
$$
which yields \eqref{Eqn:GreenFormulaDiscrete}.
\end{prf}

\section{Proof of well-posedness of the schemes and convergence results}\label{Sec:WellPosednessConvergence}
In this section, proofs to Proposition \ref{Prop:ExistenceStabilityEnergyIneqFD}, \ref{Prop:CvgFDtoSDExistenceStabilityEnergyIneqSD} and \ref{Prop:LaxWendroff} are given, starting with Proposition \eqref{Prop:ExistenceStabilityEnergyIneqFD}, i.e.~well-posedness and energy inequalities for the scheme \eqref{Eqn:FD}. It is proved using a variational method based on Lemma \ref{Lem:DiscretePoincareSobolev}.

\begin{prf}{Proposition \ref{Prop:ExistenceStabilityEnergyIneqFD}}
The solution of the scheme \eqref{Eqn:FD} is built by induction: given $\bu^n=(u_i^{n})_{i\in\Z}\in\Xh$, show that the equation 
\begin{equation}\label{Eqn:PrfEn}\tag{$E_n$}
\frac{(u^{n+1}_i)^{q-1}-(u^{n}_i)^{q-1}}{\Delta t} +\left[(-\Delta)_h^\theta \bu^{n+1}\right]_i = 0,\quad |i|\leq M_x, \quad \bu^{n+1}=(u_i^{n+1})_{i\in\Z}\in\Xh
\end{equation}
admits a unique solution. To do so, define the functional 
$$
E(\bu):= \frac{1}{q\Delta t}\|\bu\|_{l^q_h(\R)}^q + \frac{1}{2}\|\bu\|_{\X}^2 - \frac{h}{\Delta t}\sum_{i\in\Z} (u^{n}_i)^{q-1}u_i,\quad \bu:=(u_i)_{i\in\Z}\in\X.
$$
The functional $E$ is strictly convex and of class $C^1$ over $\X$, as the first two terms are strictly convex and $C^1$ and the last one is linear. Moreover,
$$
\lim_{\|\bu\|_{\X}\rightarrow +\infty} E(\bu)=+\infty.
$$
Indeed, by Hölder's and Young's inequality, it holds for any $\varepsilon>0$
\begin{equation}\label{Eqn:EstimationAvecHolderYoung}
\begin{aligned}
\frac{1}{\Delta t}\left|\sum_{i\in\Z} h (u^{n}_i)^{q-1}u_i\right|&\leq \frac{1}{\Delta t}\|\bu^{n}\|_{l^q_h(\R)}^{q-1}\|\bu\|_{l^q_h(\R)},\\
&\leq\frac{1}{\varepsilon^{q'} q'\Delta t }\|\bu^{n}\|_{l^q_h(\R)}^q+\frac{\varepsilon^q}{q \Delta t}\|\bu\|_{l^q_h(\R)}^q,
\end{aligned}
\end{equation}
where $q':q/(q-1)$. Taking $\varepsilon$ small enough, we see that there exists a constant $C\in\R$ such that
$$
E(\bu)\geq \frac{1}{2q\Delta t}\|\bu\|_{l^q_h(\R)}^2 + \frac{1}{2}\|\bu\|_{\X}^2 + C \underset{\|\bu\|_{\X}\rightarrow+\infty}{\longrightarrow}+\infty.
$$
Therefore, $E$ admits a unique minimizer $\bu^{n+1}\in\X$, and it is the unique solution of the equation
$$
\nabla_{l^2_h(I)} E(\bu^{n+1}) = 0,
$$
which is \eqref{Eqn:PrfEn}. This proves the well-posedness of \eqref{Eqn:FD}.

 To obtain the upper bound in \eqref{Eqn:EnergyIneqFD1}, test \eqref{Eqn:PrfEn} by $\bu^{n+1}$,
$$
\left( \frac{(\bu^{n+1})^{q-1}-(\bu^{n})^{q-1}}{\Delta t}+(-\Delta)^\theta_h \bu^{n+1},  \bu^{n+1} \right)_{l^2_h(\R)} = 0.
$$
Therefore,
$$
\frac{\|\bu^{n+1}\|_{l^q_h(\R)}^q}{\Delta t} + \|\bu^{n+1}\|_{\X}^2  - \left( \frac{(\bu^n)^{q-1}}{\Delta t} , \bu^{n+1}  \right)_{l^2_h(\R)} = 0.
$$ 
The last term is estimated using Hölder's and Young's inequality, as in \eqref{Eqn:EstimationAvecHolderYoung},
$$
\left|\left(  \frac{(\bu^n)^{q-1}}{\Delta t} , \bu^{n+1}  \right)_{l^2_h(\R)}\right| \leq \frac{1}{q'\Delta t}\|\bu^n\|_{l^q_h(\R)}^{q} + \frac{1}{q\Delta t}\|\bu^{n+1}\|_{l^{q}_h(\R)}^q.
$$
Hence,
$$
\left(1 -\frac{1}{q}\right)\frac{\|\bu^{n+1}\|_{l^q_h(\R)}^q}{\Delta t} - \frac{1}{q'}\frac{\|\bu^{n}\|_{l^{q}_h(\R)}^q}{\Delta t} +\|\bu^{n+1}\|_{\X}^2 \leq 0,
$$
which yields the upper bound of \eqref{Eqn:EnergyIneqFD1}. 

Before showing the lower bound of \eqref{Eqn:EnergyIneqFD1}, we begin with the proof of inequality \eqref{Eqn:EnergyIneqFD2}. Recall the following elementary inequality (see \cite{akagi_existence} for a proof),
\begin{equation}\label{Eqn:PrfIneqElem}
\frac{4}{qq'}(x^{q/2}-y^{q/2})^2\leq (x^{q-1}-y^{q-1})(x-y)\quad\text{for any $x,y\in\R$}.
\end{equation}
Testing \eqref{Eqn:PrfEn} with $(\bu^{n+1}-\bu^n)/\Delta t$ gives
$$
\left(\frac{(\bu^{n+1})^{q-1}-(\bu^n)^{q-1}}{\Delta t}, \frac{\bu^{n+1}-\bu^n}{\Delta t}\right)_{l^2_h(\R)} + \frac{\|\bu^{n+1}\|_{\X}^2 - (\bu^{n+1},\bu^n)_{\X}}{\Delta t} = 0,
$$
where the first term is estimated using \eqref{Eqn:PrfIneqElem},
$$
\frac{4}{qq'}\left\| \frac{(\bu^{n+1})^{q/2}-(\bu^{n})^{q/2}}{\Delta t} \right\|_{l^2_h(\R)}^2 \leq \left(\frac{(\bu^{n+1})^{q-1}-(\bu^n)^{q-1}}{\Delta t}, \frac{\bu^{n+1}-\bu^n}{\Delta t}\right)_{l^2_h(\R)}.
$$
Next, Cauchy-Schwarz and Young's inequalities give
$$
(\bu^{n+1},\bu^n)_{\X}\leq \frac{1}{2}\|\bu^n\|_{\X}^{2} + \frac{1}{2}\|\bu^{n+1}\|_{\X}^2,
$$
which yields \eqref{Eqn:EnergyIneqFD2}.

To prove the lower bound of \eqref{Eqn:EnergyIneqFD1}, \eqref{Eqn:PrfEn} is tested by $\bu^n$. By an argument similar to the proof of the upper bound of \eqref{Eqn:EnergyIneqFD1}, one has
$$
\left|\left(\frac{(\bu^{n+1})^{q-1}-(\bu^n)^{q-1}}{\Delta t},\bu^n\right)_{l_h^2(\R)}\right|\leq \frac{1}{q'}\frac{\|\bu^{n+1}\|_{l^q_h(\R)}^q - \|\bu^{n}\|_{l^q_h(\R)}^q}{\Delta t}.
$$
Moreover, thanks to \eqref{Eqn:EnergyIneqFD2}, $(\|\bu^n\|_{\X})_{n\in\N}$ is non-increasing, so that Cauchy-Schwarz inequality gives
$$
\left((-\Delta)^\theta_h \bu^{n+1}, \bu^n\right)_{l_h^2(\R)}=(\bu^{n+1},\bu^n)_{\X}\leq \|\bu^{n}\|_{\X}^2,
$$
and it yields the lower bound in \eqref{Eqn:EnergyIneqFD1}.
\end{prf}

\details{
\begin{rem}
Inequality \eqref{Eqn:PrfIneqElem} can also be written
$$
\frac{4}{qq'}(x^{1/2}-y^{1/2})^2\leq (x^{1/q}-y^{1/q})(x^{1/q'}-y^{1/q'})\quad\text{for any $x,y\in\R$}.
$$
\end{rem} 
}

After proving well-posedness of the fully discrete scheme \eqref{Eqn:FD}, it is possible to show Proposition \ref{Prop:CvgFDtoSDExistenceStabilityEnergyIneqSD}, that is well-posedness of the semi-discrete scheme \eqref{Eqn:SD}, by passing to the limit in \eqref{Eqn:FD} as $\Delta t\rightarrow 0$. Estimates \eqref{Eqn:EnergyIneqFD1} and \eqref{Eqn:EnergyIneqFD2} provide compactness and allow to obtain convergence of the solutions to the fully discrete scheme. Then, energy estimate \eqref{Eqn:EnergyIneqSD2} is obtained by passing to the limit in the discrete energy estimate \eqref{Eqn:EnergyIneqFD2}. The proof follows the same approach as the continuous case (see \cite{akagi_existence}), but is considerably simplified since $\Xh$ is a finite dimensional space, and thus, \eqref{Eqn:SD} is a nonlinear ODE. However, it is worth noticing that a strategy based on the Cauchy-Lipschitz Theorem would not allow to obtain the energy estimate \eqref{Eqn:EnergyIneqSD2}.

\begin{prf}{Proposition \ref{Prop:CvgFDtoSDExistenceStabilityEnergyIneqSD}}
Let $\bu=(u_i^n)_{n\in\N,i\in\Z}$ be a solution to \eqref{Eqn:FD}, $T>0$, and recall that $\bv^{\Delta t}:\R_+\ni t\mapsto (v_i^{\Delta t}(t))_{i\in\Z}$ is defined by $v_i^{\Delta t}(t)=u_i^{n+1}$ for $t\in(n\Delta t,(n+1)\Delta t]$.
Let us show that $\bv^{\Delta t}$ converges to a solution of \eqref{Eqn:SD} as $\Delta t\rightarrow 0$, at a \emph{fixed space step} $h$.

Introduce $\boldsymbol{\rho}^{\Delta t}:\R_+\ni t \mapsto  (\rho^{\Delta t}_i)_{i\in\Z}$ and $\boe^{\Delta t}:\R_+\ni t\mapsto (\eta^{\Delta t}_i)_{i\in\Z}$, piecewise linear interpolants of $((u^n_i)^{q-1})_{i\in\Z,n\in\N}$ and $((u^n_i)^{q/2})_{i\in\Z,n\in\N}$, defined by
\begin{gather}
\rho_i^{\Delta t}(t):=(u_i^n)^{q-1} + \frac{(u_i^{n+1})^{q-1}-(u_i^{n})^{q-1}}{\Delta t}(t-t_n),\quad\text{for $t\in [t_n,t_{n+1}]$}, \label{Eqn:DefRhoDeltat}\\
\eta_i^{\Delta t}:=(u_i^n)^{q/2} + \frac{(u_i^{n+1})^{q/2}-(u_i^{n})^{q/2}}{\Delta t}(t-t_n),\quad\text{for $t\in [t_n,t_{n+1}]$} \label{Eqn:DefEtaDeltat}.
\end{gather}
Since $\bu$ is the solution of the fully discrete scheme \eqref{Eqn:FD}, 
\begin{align}
\partial_t \rho_i^{\Delta t} (t)+\left[ (-\Delta)^\theta_h \bv^{\Delta t}(t)\right]_i =& 0,\quad{\forall t\neq t_n},\quad|i|\leq M_x, \label{Eqn:EquationFullyDiscRhoVDeltat}\\
\rho_i^{\Delta t}(t) =& 0,\quad \forall t>0,\quad |i|>M_x, \label{Eqn:CondBordRhoDeltat}\\
\rho_i^{\Delta t}(0) =& (u^0_i)^{q-1},\quad \forall i\in\Z.  \label{Eqn:CondInitialeRhoDeltat}
\end{align}
Moreover, \eqref{Eqn:EnergyIneqFD2} implies
\begin{equation}
\label{Eqn:EnergyIneqFD2EtaVDeltat}
\frac{4}{qq'}\int_0^t \|\partial_s \boe^{\Delta t}(s)\|_{l^2_h(\R)}^2{\rm d}s + \frac{\|\bv^{\Delta t}(t)\|_{\Xh}^2}{2}\leq \frac{\|\bu^{0}\|_{\Xh}^2}{2},\quad \forall t>0,
\end{equation}
and it yields the existence of a constant $C$ independent of $\Delta t$ such that
\begin{equation}\label{Eqn:BoundVDeltat}
\|\bv^{\Delta t}(t)\|_{\Xh}<C,\quad \forall t>0.
\end{equation}
Unsing the discrete Green formula \eqref{Eqn:GreenFormulaDiscrete}, one has $\|(-\Delta)^\theta_h \bv^{\Delta t}(t)\|_{\Xh^*} = \|\bv^{\Delta t}(t)\|_{\Xh}$, so that \eqref{Eqn:BoundVDeltat} and \eqref{Eqn:EquationFullyDiscRhoVDeltat} yield
\begin{equation}
\|\partial_t\bor^{\Delta t}(t)\|_{\Xh^*}<C,\quad \forall t>0 \label{Eqn:BoundRhoDeltat}.
\end{equation}
Since $\Xh$ and $\Xh^*$ are finite dimensional spaces, \eqref{Eqn:BoundVDeltat} and \eqref{Eqn:BoundRhoDeltat} yield the existence of a constant $C_h$ independant of $\Delta t$ such that,
\begin{gather}
|v_i^{\Delta t}(t)|<C_h,\quad \forall t>0, \quad\forall i\in \Z \label{Eqn:BoundVDeltatInfinity}\\
|\partial_t\rho_i^{\Delta t}(t)|<C_h,\quad \forall t>0,\quad\forall i\in \Z \label{Eqn:BoundRhoDeltatInfinity}.
\end{gather}
It is worth noticing that $C_h$ may depends on $h$, but it is not an issue since $h$ is fixed. Moreover, by definition of $\bv^{\Delta t}$ and $\bor^{\Delta t}$ (see \eqref{Eqn:DefRhoDeltat}),
\begin{equation}\label{Eqn:ErreurRhoDeltatVDeltat}
\forall n\in\N,\ \forall t\in [t_n,t_{n+1}],\ \forall i\in\Z,\quad|\rho^{\Delta t}_i(t)-(v^{\Delta t}_i)^{q-1}(t)|\leq \int_{t_n}^{t_{n+1}} |\partial_t \rho^{\Delta t}_i(t)|{\rm d}t\leq C\Delta t.
\end{equation}
Equations \eqref{Eqn:BoundVDeltatInfinity}, \eqref{Eqn:BoundRhoDeltatInfinity} and \eqref{Eqn:ErreurRhoDeltatVDeltat} imply
$$
\|\rho_i^{\Delta t}\|_{W^{1,\infty}(0,\infty;\R)}<C,\quad\forall i\in\Z.
$$
Therefore, Rellich-Kondrakov Theorem yields the existence of a function $\bor:\R_+\ni t\mapsto (\rho_i(t))_{i\in\Z}$, with $\rho_i\in C([0,\infty),\R)$ for any $i\in\Z$, and a sequence $\Delta t_k\rightarrow 0$ such that
\begin{gather}
\rho_i^{\Delta t_k} {\rightarrow} \rho_i\quad\text{in $C([0,T],\R)$},\quad\forall i\in\Z,  \label{Eqn:CvgRhoDeltatLInfty}\\
\rho_i^{\Delta t_k} \overset{\star}{\rightharpoonup} \rho_i\quad\text{in $W^{1,\infty}(0,T;\R)$},\quad\forall i\in\Z  \label{Eqn:CvgRhoDeltatFaible},
\end{gather}
as $k\rightarrow \infty$. Finally, \eqref{Eqn:ErreurRhoDeltatVDeltat} together with \eqref{Eqn:CvgRhoDeltatLInfty} implies
\begin{equation}
v_i^{\Delta t_k}\rightarrow v_i:=\rho_i^{1/(q-1)}\quad\text{in $L^\infty(0,T;\R)$},\quad\forall i\in \Z.\label{Eqn:CvgVDeltatLInfty}
\end{equation}
Indeed, $x\mapsto x^{1/(q-1)}$ is uniformly continuous on every bounded intervals of $\R$, and $\rho_i^{\Delta t}$ is bounded on $[0,T]$ uniformly in $\Delta t$ by \eqref{Eqn:BoundVDeltatInfinity}. One can pass to the limit in \eqref{Eqn:EquationFullyDiscRhoVDeltat} using \eqref{Eqn:CvgRhoDeltatFaible} and \eqref{Eqn:CvgVDeltatLInfty}, and obtain the desired ODE,
$$
\frac{{\rm d}}{{\rm d}t} v_i^{q-1}(t) + [(-\Delta)_h^\theta \bv(t)]_i = 0,\quad\text{for $|i|\leq M_x$, and $t<T$},
$$
where $\bv:=t\in\R_+\mapsto (v_i(t))_{i\in\Z}$.
Since $T$ is arbitrary, $\bv$ satisfies the ODE in \eqref{Eqn:SD}. Moreover, the boundary and initial conditions for $\bor^{\Delta t}$ \eqref{Eqn:CondBordRhoDeltat}, \eqref{Eqn:CondInitialeRhoDeltat}, together with the uniform convergence of $v_i^{\Delta t}$ to $v_i=\rho_i^{1/(q-1)}$ \eqref{Eqn:CvgVDeltatLInfty}, show that the boundary and initial conditions of \eqref{Eqn:SD} are satisfied by $\bv$. Therefore, $\bv$ is a solution to \eqref{Eqn:SD}. 

The proof of uniqueness is identical to the continuous case and we refer to \cite{akagi_existence}.
It remains to prove the two energy inequalities \eqref{Eqn:EnergyIneqSD1}, \eqref{Eqn:EnergyIneqSD2}. To prove \eqref{Eqn:EnergyIneqSD1}, it is worth noticing that, since $q'=q/(q-1)>1$, the function $x\in\R\mapsto x^{q'}$ is $C^1(\R)$, and thus $v_i^q = (v_i^{q-1})^{q'} \in W^{1,1}(0,T;\R)$ for any $i\in\Z$ and $T>0$. Moreover, the chain rule applies: $v_i(t)\partial_t v_i^{q-1}(t) =(1/q')\partial_t v_i^q(t)$. Therefore, multiplying the ODE in \eqref{Eqn:SD} by $v_i(t)$ and summing in $i$ yields energy identity \eqref{Eqn:EnergyIneqSD1}.

Finally, the second energy inequality \eqref{Eqn:EnergyIneqSD2} is obtained by passing to the limit in \eqref{Eqn:EnergyIneqFD2EtaVDeltat}. From \eqref{Eqn:EnergyIneqFD2EtaVDeltat}, there exists a constant $C_h$ independant of $\Delta t$ such that
\begin{equation}\label{Eqn:BoundEtaDeltatL2}
    \|\partial_t \eta_i^{\Delta t}\|_{L^2(0,\infty;\R)}<C_h,\quad \forall i\in\Z.
\end{equation}
Moreover, by definition of $\bv^{\Delta t}$ and $\boe^{\Delta t}$ (see \eqref{Eqn:DefEtaDeltat}),
\begin{equation}\label{Eqn:ErreurEtaDeltatVDeltat}
\forall n\in\N,\ \forall t>0,\ \forall i\in\Z,\quad|\eta_i^{\Delta t}(t)- (v_i^{\Delta t})^{q/2}(t)|
\leq \int_{t_n}^{t_{n+1}} |\partial_s \eta^{\Delta t}_i(s)|{\rm d}s\leq C_h\sqrt{\Delta t}.
\end{equation}
Equations \eqref{Eqn:BoundVDeltatInfinity}, \eqref{Eqn:BoundEtaDeltatL2}, \eqref{Eqn:ErreurEtaDeltatVDeltat} imply, up to relabelling the constant $C_h$,
$$
\|\eta_i^{\Delta t}\|_{W^{1,2}(0,T;\R)}<C_h.
$$
Therefore, Rellich-Kondrakov Theorem yields the existence of a function $\boe:\R_+\ni t\mapsto (\eta_i(t))_{i\in\Z}$, with $\eta_i\in C([0,\infty),\R)$ for any $i\in\Z$, and a sequence $\Delta t_k\rightarrow 0$ such that,
\begin{gather}
\eta_i^{\Delta t_k} {\rightarrow} \eta_i\quad\text{in $C([0,T],\R)$},\quad\forall i\in\Z,  \label{Eqn:CvgEtaDeltatLInfty}\\
\eta^{\Delta t_k}_i \overset{*}{\rightharpoonup} \eta_i\quad\text{in $W^{1,2}(0,T;\R)$},\quad\forall i\in\Z. \label{Eqn:CvqEtaDeltatFaible}
\end{gather}
Finally, \eqref{Eqn:ErreurEtaDeltatVDeltat} together with \eqref{Eqn:CvgEtaDeltatLInfty} implies $\eta_i = v_i^{q/2}$ for any $i\in\Z$. Indeed, the function $x\in\R\mapsto x^{q/2}$ is uniformly continuous on every bounded interval of $\R$, and $v_i^{\Delta t}$ is bounded uniformly in $\Delta t$. Therefore, $v_i^{q/2}\in W^{1,2}(0,T;\R)$ for any $i \in \Z$, and passing to the limit in \eqref{Eqn:EnergyIneqFD2EtaVDeltat} (using \eqref{Eqn:CvqEtaDeltatFaible} for the first term and \eqref{Eqn:CvgVDeltatLInfty} for the second term) yields
\begin{equation}\label{Eqn:EnergyIneqSD2Startingt0}
			\frac{4}{qq'}\int_0^t \|\partial_\tau \bv^{q/2}(\tau)\|_{l^2_h}^2{\rm d}\tau + \frac{1}{2}\|\bv(t)\|_{\Xh}^2\leq\frac{1}{2}\|\bu^0\|_{\Xh}^2, \quad \forall t>0.
\end{equation}
Inequality \eqref{Eqn:EnergyIneqSD2Startingt0} is the desired energy inequality \eqref{Eqn:EnergyIneqSD2} but with $s=0$. To obtain \eqref{Eqn:EnergyIneqSD2}, one needs to write \eqref{Eqn:EnergyIneqSD2Startingt0} between two general times $0\leq s< t$. To do so, it is worth noticing that $\tilde{\bv}:=\bv(\cdot-s)$ is the unique solution to \eqref{Eqn:SD} with initial condition $\bv(s)$. Therefore, applying \eqref{Eqn:EnergyIneqSD2Startingt0} to $\tilde{\bv}$ yields \eqref{Eqn:EnergyIneqSD2}. This ends the proof of Proposition \ref{Prop:CvgFDtoSDExistenceStabilityEnergyIneqSD}.
\end{prf}

Proposition \ref{Prop:CvgFDtoSDExistenceStabilityEnergyIneqSD} shows that the fully discrete scheme \eqref{Eqn:FD} converges to the semi-discrete scheme \eqref{Eqn:SD}. Proposition \ref{Prop:LaxWendroff}, showing that \eqref{Eqn:FD} converges to the continuous problem \eqref{Eqn:CDP} under Lax-Wendroff type hypothesis, is proved now. We draw on the argument used in \cite[Proposition 2.10]{ayi_herda_hivert_tristani}.

\begin{prf}{Proposition \ref{Prop:LaxWendroff}}
Let $\bu=(u_i^n)_{n\in\N,i\in\Z}$ be a solution to \eqref{Eqn:FD} with some initial condition $\bu^0=(u^0_i)_{i\in\Z}$. Recall that $u_{\Delta t, h}:\R_+\times\R\rightarrow\R$ is defined by $u_{\Delta t,h}(t,x) = u_i^n$ for $(t,x)\in [t_n,t_{n+1})\times[x_i-h/2,x_i+h/2)$, and $u^0_h:\R\rightarrow \R$ is defined by $u^0_h(x) = u^0_i$ for $x\in [x_i-h/2,x_i+h/2)$. Recall also that $u_{\Delta t,h}$ is assumed to converge almost everywhere to some function $u\in L^\infty((0,\infty)\times\R)$ as $(\Delta t,h)\rightarrow 0$, and to be bounded in $L^\infty((0,T)\times \R)$ norm, uniformly in $\Delta t$ and $h$, for every $T>0$. Finally, recall that $u^0_h$ is assumed to converge to some function $u^0\in L^1(I)$ in $L^1(I)$ norm.

Let $\varphi\in C_c^\infty ([0,\infty)\times I)$, and let $T>0$ be such that $\text{supp}\,\varphi\subset [0,T]\times I$. Introduce the piecewise linear function $\varphi_{\Delta t,h}$, such that for any $t\in [t_n,t_{n+1})$, and $x\in[x_i-h/2,x_i+h/2)$,
$$
\varphi_{\Delta t,h}(x,t)=\varphi(t_n,x_i)+\frac{\varphi(t_n,x_{i+1})-\varphi(t_n,x_{i})}{h}(x-x_i).
$$
For any function $\phi$ bounded, denote also 
$$
(-\Delta)^\theta_h \phi (x):= \sum_{j\in\Z} w_j^{h}(\phi(x)-\phi(x- jh)).
$$
By definition, $\varphi_{\Delta,h}$ satisfies
\begin{equation}
\begin{aligned}\label{Eqn:CalculPhiDelta}
\MoveEqLeft[7]\int_{x_i-h/2}^{x_i+h/2}\left[(-\Delta)_h^\theta \varphi_{\Delta t,h}(t_n,\cdot)\right](x){\rm d}x \\
=&\sum_{j\in\Z} w_j^{h}\left(\int_{x_i-h/2}^{x_i+h/2}\varphi_{\Delta t,h}(t_n,x){\rm d}x-\int_{x_i-h/2}^{x_i+h/2}\varphi_{\Delta t,h}(t_n,x-jh){\rm d}x\right),\\
=& h\sum_{j\in\Z} w_j^{h}\left(\varphi(t_n,x_i)-\varphi(t_n,x_{i-j})\right),\\
=& h\left[(-\Delta)^\theta_h \varphi(t_n,\cdot)\right](x_i).
\end{aligned}
\end{equation}
Moreover, as $\varphi_{\Delta t,h}$ interpolates $\varphi$, one has, for some constant $K_1$, 
\begin{equation}
\|\varphi(t_n,\cdot) - \varphi_{\Delta t,h}(t_n,\cdot)\|_{\infty}\leq K_1 h^2,
\end{equation}
so that 
\begin{equation}\label{Eqn:ErreurPhiPhiDelta}
\begin{aligned}
\|(-\Delta)^\theta_h (\varphi - \varphi_{\Delta t,h})(t_n,\cdot)\|_\infty\leq& 2\|\varphi(t_n,\cdot) - \varphi_{\Delta t,h}(t_n,\cdot)\|_\infty\sum_{j\in\Z}w_j^{h}\\
\leq & K_2 h^{2-2\theta}\rightarrow 0,
\end{aligned}
\end{equation}
for some constant $K_2$.

Multiplying the ODE in \eqref{Eqn:FD} by $\varphi(t_{n+1},x_i)$, and summing in $n\in\N$ and $i\in \Z$ yields a discrete weak formulation. Indeed, denoting $\varphi_i^n=\varphi(t_n,x_i)$, 
$$
\sum_{i\in Z}\sum_{n\in\N} h\Delta t \frac{(u_i^{n+1})^{q-1}-(u_i^{n})^{q-1}}{\Delta t} \varphi_i^{n+1} + \sum_{i\in \Z}\sum_{n\in\N} h\Delta t \left[(-\Delta)^\theta_h u^{n+1}\right]_i\varphi_i^{n+1} = 0.
$$
Recall that $\varphi\in C_c^\infty ([0,\infty)\times I)$, so that both sums are actually finite. A discrete integration by parts yields
\begin{multline}\label{Eqn:DiscreteWeakForm}
-\sum_{i\in Z}\sum_{n\geq 1} h\Delta t (u_i^n)^{q-1}\frac{\varphi_i^{n+1}-\varphi_i^n}{\Delta t} + \sum_{i\in\Z}\sum_{n\in\N} h\Delta t u^{n+1}_i \left[(-\Delta)^\theta_h \varphi(t_{n+1},\cdot)\right](x_i) \\ - \sum_{i\in \Z} h (u^0_i)^{q-1}\varphi_i^1 = 0.
\end{multline}
The rest of the proof consists in proving that the three terms in \eqref{Eqn:DiscreteWeakForm} converge to the three terms in the weak form \eqref{Eqn:EqFaible}. For the second term, using the computation \eqref{Eqn:CalculPhiDelta},
\begin{align*}
\MoveEqLeft[10] \sum_{i\in\Z}\sum_{n\in\N} h\Delta t u^{n+1}_i \left[(-\Delta)^\theta_h \varphi(t_{n+1},\cdot)\right](x_i)\\
=& \sum_{n\geq 0}\Delta t\sum_{i\in\Z}\int_{x_i-h/2}^{x_i+h/2}u_{\Delta t,h}(t_{n+1},x)(-\Delta)^\theta_h \left[\varphi_{\Delta t,h}(t_{n+1},\cdot)\right](x) {\rm d}x\\
=&\sum_{n\geq 0}\Delta t\int_{\R} u_{\Delta t,h}(t_{n+1},x)(-\Delta)^\theta_h \left[\varphi_{\Delta t,h}(t_{n+1},\cdot)\right](x) {\rm d}x\\
=& \int_{\Delta t}^\infty \int_\R u_{\Delta t,h}(t,x)\left[(-\Delta)^\theta\varphi(t,\cdot)\right](x){\rm d}x{\rm d}t \\
&\hspace{99,592pt}+E_1(\Delta t,h)+E_2(\Delta t,h)+E_3(\Delta t,h),
\end{align*}
with
\begin{gather*}
E_1(\Delta t,h) = \int_{\Delta t}^\infty \int_\R u_{\Delta t,h}(t,x)\left(\left[(-\Delta)^\theta_h \varphi(t,\cdot)\right](x) - \left[(-\Delta)^\theta \varphi(t,\cdot)\right](x)\right){\rm d}x,\\
E_2(\Delta t,h)=\sum_{n\geq 1} \int_{\R} u_{\Delta t,h}(t_n,x)\left( \Delta t \left[(-\Delta)^\theta_h \varphi(t_n,\cdot)\right](x)  -  \int_{t_n}^{t_{n+1}}\left[(-\Delta)^\theta_h \varphi(t,\cdot)\right](x){\rm d}t \right) {\rm d}x, \\
E_3(\Delta t,h) = \sum_{n\geq 1} \int_\R u_{\Delta t,h}(t_n,x)\Delta t\left[(-\Delta)^\theta_h(\varphi_{\Delta t,h}(t_n,\cdot) -  \varphi(t_n,\cdot))\right](x){\rm d}x.
\end{gather*}
On the one hand, Theorem \ref{Thm:ErreurDiscreteFracLap} and the boundedness of $u_{\Delta t,h}$ (uniform in $\Delta t$ and $h$) yields $E_1(\Delta t,h)\rightarrow 0$ as $h\rightarrow 0$, uniformly in $\Delta t$. Similarly, using \eqref{Eqn:ErreurPhiPhiDelta}, $E_3(\Delta t,h)\rightarrow 0$ as $h\rightarrow 0$, uniformly in $\Delta t$. Also, $E_2$ is a sum of local quadrature errors:
$$
\left|\Delta t \left[(-\Delta)^\theta_h \varphi(t_n,\cdot)\right](x)  -  \int_{t_n}^{t_{n+1}}\left[(-\Delta)^\theta_h \varphi(t,\cdot)\right](x){\rm d}t\right|\leq \frac{M_1}{2} \Delta t^2,
$$
with $M_1:=\sup_{t>0,x\in I} \partial_t \left[(-\Delta)^\theta_h \varphi(t,\cdot)\right](x)$. 
Hence, $|E_2(\Delta t,h)|\leq K\Delta t$ for some constant $K>0$. Therefore, $E_2(\Delta t,h)\rightarrow 0$ as $\Delta t\rightarrow 0$, uniformly in $h$.

For the first term in \eqref{Eqn:DiscreteWeakForm}, using the definition of $(\varphi_i^n)_{n\in\N,i\in\Z}$ and $u_{\Delta t,h}$,
\begin{align*}
-\sum_{i\in Z}\sum_{n\geq 1} h\Delta t (u_i^n)^{q-1}\frac{\varphi_i^{n+1}-\varphi_i^n}{\Delta t} =& - \sum_{i\in\Z}\sum_{n\geq 1} h(u_i^n)^{q-1}\int_{t_n}^{t_{n+1}}\partial_t \varphi(t,x_i){\rm d}t\\
=& -\int_{\Delta t}^\infty \int_{\R} u_{\Delta t,h}^{q-1}(t,x)\partial_t\varphi(t,x){\rm d}x{\rm d}t + E_4(\Delta t,h),
\end{align*}
with
$$
E_4(\Delta t,h) = \int_{\Delta t}^\infty \sum_{i\in\Z} u_{\Delta t,h}^{q-1}(t,x)\left[\int_{x_i - h/2}^{x_i+h/2}\partial_t\varphi(t,x){\rm d}x-h\partial_t\varphi(t,x_i)\right]{\rm d}t.
$$
Since $E_4(\Delta t,h)$ is a sum of local quadrature errors, as above, $|E_4(\Delta t,h)|\rightarrow 0$ as $h\rightarrow 0$, uniformly in $\Delta t$.

Finally, for the third term in \eqref{Eqn:DiscreteWeakForm}, using the definition of $u^0_h$ and $(\varphi^1_i)_{i\in\Z}$,
$$
-\sum_{i\in\Z} h (u^0_i)^{q-1}\varphi_i^1 = - \int_{\R} (u^0_h)^{q-1}(x)\varphi(\Delta t,x){\rm d}x + E_5(\Delta t,h),
$$
with 
$$
E_5(\Delta t,h)= \sum_{i\in\Z} (u_h^0)^{q-1} (x)\left[\int_{x_i-h/2}^{x_i+h/2}\varphi(\Delta t,x){\rm d}x - h\varphi(\Delta t, x_{i})\right].
$$
Again, since $E_5(\Delta t,h)$ is a sum of local quadrature errors, we see that $E_5(\Delta t,h)\rightarrow 0$ as $h\rightarrow 0$, uniformly in $\Delta t$. Therefore,
\begin{multline*}
-\int_{\Delta t}^T\int_I u_{\Delta t,h}^{q-1}(t,x)\partial_t\varphi(t,x){\rm d}x {\rm d}t+\int_{\Delta t}^T\int_I u_{\Delta t,h}(t,x)(-\Delta)^\theta \varphi(t,x){\rm d}x {\rm d}t \\
- \int_I (u^0_h)^{q-1}(x)\varphi(\Delta t,x){\rm d}x
=-E_1(\Delta t,h)-E_2(\Delta t,h)-E_3(\Delta t,h)-E_4(\Delta t,h)-E_5(\Delta t,h).
\end{multline*}
Passing to the limit as $(\Delta t,h)\rightarrow 0$ using dominated convergence yields \eqref{Eqn:EqFaible}.
\end{prf}

\section{Proof of energy decay estimates}\label{Sec:DecayEstimates}

The goal of this section is to prove Propositions \ref{Prop:EstimatesSD} and \ref{Prop:EstimatesFD}. We start with the proof of Proposition \ref{Prop:EstimatesSD}, i.e.~ the $l^q_h(\R)$-decay estimates for the semi-discrete scheme \eqref{Eqn:SD}. Thanks to the discrete functional analysis tools introduced in Lemma \ref{Lem:DiscreteSobolevIneq}, the proof is almost identical to the one used in \cite{akagi_existence} for proving the continuous estimates \eqref{Eqn:Asymp1Continue}, \eqref{Eqn:Asymp2Continue}. 

\begin{prf}{Proposition \ref{Prop:EstimatesSD}}
    In the case $q\in(2,2_\theta^*)$, it will be shown that the solution $\bv$ of \eqref{Eqn:SD} extincts in finite time. Hence, we shall denote
    $$
    t_*^h:=\sup\{t>0:\:u(t)\not\equiv 0\},
    $$
the \emph{extinction time}. It must be shown that $t_*^h<+\infty$ when $q\in(2,2_\theta^*)$, and $t_*^h=+\infty$ when $q\in (1,2)$. In the case where $t_*^h < +\infty$, the uniqueness to \eqref{Eqn:SD} (Proposition \ref{Prop:CvgFDtoSDExistenceStabilityEnergyIneqSD}) and the continuity of $t\geq 0 \mapsto \|\bv(t)\|_{\lq}^2$ imply $v(t)\equiv 0$ for any $t\geq t_*$. Therefore, the decay of the Rayleigh quotient $R(t) = \|\bv(t)\|_{\Xh}^2/\|\bv(t)\|_{\lq}^2$, as well as estimate \eqref{Eqn:Asymp2DiscreteSD} need to be shown only for $t\in[0,t_*^h)$. 

We shall start with the decay of the Rayleigh quotient. Consider a time $T<t_*^h$. The Rayleigh quotient $R$ may not be absolutely continuous on $[0,T]$, but it is the product of the non-increasing function $[0,T]\ni t \mapsto \|\bv(t)\|_{\Xh}^2$, and the Lipschitz continuous function $[0,T]\ni t \mapsto 1/\|\bv(t)\|_{\lq}^2$. Indeed $\|\bv(t)\|_{\Xh(I)}^2$ is non-increasing thanks to \eqref{Eqn:EnergyIneqSD1} and $\|\bv(t)\|_{\lq}^2$ is uniformly away from $0$ on $[0,T]$ since $T<t_*^h$. Moreover, $\|\bv(t)\|_{\lq}^q$ is absolutely continuous on $[0,T]$, and thanks to \eqref{Eqn:EnergyIneqSD1} and \eqref{Eqn:EnergyIneqSD2}, the time derivative $(\d/\d t)\|\bv(t)\|_{\lq}^q$ is bounded, which yields the Lipschitz continuity of $[0,T]\ni t \mapsto 1/\|\bv(t)\|_{\lq}^2$. Therefore, as shown in \cite{akagi_existence}, proving $({\rm d}/{\rm d}t)R(t)\leq 0$ for almost every $t<t_*$ is enough to obtain the decay of $R$. Roughly speaking, since $R$ is a function of bounded variation, its distributional derivative is a Radon measure $\mu$. This measure has an absolutely continuous part given by the pointwise derivative $(\d/\d t)R(t)$, but it may also have a singular part $\mu_s$ if $R$ is not absolutely continuous. However, since $[0,T]\ni t \mapsto 1/\|\bv(t)\|_{\lq}^2$ is absolutely continuous, the singular part results from the term $\|\bv(t)\|_{\Xh}^2$. Since $[0,T]\ni t \mapsto \|\bv(t)\|_{\Xh}^2$ is non-increasing, it can be shown that $\mu_s \leq 0$, and it is enough to show $({\rm d}/{\rm d}t)R\leq 0$ almost everywhere to obtain $\mu\leq 0$. See \cite{akagi_existence} for a detailed argument.

We show now that $({\rm d}/{\rm d}t)R\leq 0$ almost everywhere on $(0,t_*^h)$. The computations are identical to the continuous case \cite{akagi_existence}, and present them for completeness. For almost all $t\in[0,t_*^h)$,
\begin{equation}\label{Eqn:RayleighQuotientDerivation}
\frac{\d}{\d t}R(t)=\frac{\|\bv(t)\|_{\lq}^2\frac{\d}{\d t}\|\bv(t)\|_{\Xh}^2-\|\bv(t)\|_{\Xh}^2\frac{\d}{\d t}\|\bv(t)\|_{\lq}^2}{\|\bv(t)\|_{\lq}^4}.
\end{equation}
Estimate \eqref{Eqn:EnergyIneqSD1} yields
$$
\|\bv(t)\|_{\Xh}^2=-\frac{1}{q'}\frac{\d}{\d t}\|\bv(t)\|_{\lq}^q,\quad\text{for almost all $t\in(0,t_*^h)$}.
$$
Estimate \eqref{Eqn:EnergyIneqSD2} implies, by multiplying by $1/(t-s)$ and letting $t\rightarrow s$,
$$
\frac{\d}{\d t}\|\bv(t)\|_{\Xh}^2\leq -\frac{8}{qq'}\|\partial_t \bv^{q/2}(t)\|_{l^2_h(\R)}^2\quad\text{for almost all $t\in(0,t_*^h)$}.
$$
Moreover, using the chain rule, 
$$
\frac{\d}{\d t}\|\bv(t)\|_{\lq}^2=\frac{2}{q}\|\bv(t)\|_{\lq}^{2-q}\frac{\d}{\d t}\|\bv(t)\|_{\lq}^q,\quad\text{for almost all $t\in(0,t_*^h)$}.
$$
Hence, replacing in \eqref{Eqn:RayleighQuotientDerivation}, for almost all $t\leq t_*^h$,
$$
\begin{aligned}
\frac{\d}{\d t}R(t)\leq& \frac{-\frac{8}{qq'}\|\partial_t \bv^{q/2}(t)\|_{l^2_h(\R)}^2\|\bv(t)\|_{\lq}^2+\frac{2}{qq'}\|\bv(t)\|_{\lq}^{2-q}\left(\frac{\d}{\d t}\|\bv(t)\|_{\lq}^q\right)^2}{\|\bv(t)\|_{\lq}^4},\\
=&\frac{-\frac{8}{qq'}\|\partial_t \bv^{q/2}(t)\|_{l^2_h(\R)}^2\|\bv(t)\|_{\lq}^q+\frac{2}{qq'}\left(\frac{\d}{\d t}\|\bv(t)\|_{\lq}^q\right)^2}{\|\bv(t)\|_{\lq}^{q+2}}.
\end{aligned}
$$
Finally, we observe
$$
\begin{aligned}
\frac{\d}{\d t}\|\bv(t)\|_{\lq}^q=&\frac{\d}{\d t}\|\bv^{q/2}(t)\|_{l^2_h}^2,\\
=&2(\partial_t \bv^{q/2}(t),\bv^{q/2}(t))_{l^2(\R)},\quad\text{for almost all $t\in(0,t_*^h)$},
\end{aligned}
$$
so Cauchy-Schwartz inequality yields $({\rm d}/{\rm d}t)R(t)\leq 0$ for almost all $t\leq t_*^h$, which proves the decay of $t\in[0,t_*^h)\mapsto R(t)$.

We prove next estimates \eqref{Eqn:Asymp1DiscreteSD}, \eqref{Eqn:Asymp2DiscreteSD}. Assume $q\in(1,2_\theta^*]$. The energy estimates \eqref{Eqn:EnergyIneqSD1} and \eqref{Eqn:EnergyIneqSD2} combined with the discrete Poincaré-Sobolev inequality \eqref{Eqn:DiscretePoincareSobolev} and the decay of the Rayleigh quotient yields
\begin{gather}\label{Eqn:EstimatesSDDerivationSobolevIneq}
\frac{1}{q'}\frac{{\rm d}}{{\rm d}t}\|\bv(t)\|_{\lq}^q+C^{-2}\|\bv(t)\|_{\lq}^2 \leq 0 \quad\text{for almost all $t>0$},\\
\label{Eqn:EstimatesSDDerivationRayleighDecay}
\frac{1}{q'}\frac{{\rm d}}{{\rm d}t}\|\bv(t)\|_{\lq}^q+R(0)\|\bv(t)\|_{\lq}^2 \geq 0 \quad\text{for almost all $t>0$}.
\end{gather}
Next, note that if $q\neq 2$, then
$$
\frac{1}{q'}\frac{{\rm d}}{{\rm d}t}\|\bv(t)\|_{\lq}^q = \frac{q-1}{q-2}\|\bv(t)\|_{\lq}^2\frac{{\rm d}}{{\rm d}t}\|\bv(t)\|^{q-2}_{\lq},\quad\text{for almost all $t<t_*$,}
$$
so that
\begin{gather}\label{Eqn:EstimatesSDDerivationSobolevIneqInter}
\frac{{\rm d}}{{\rm d}t}\|\bv(t)\|_{\lq}^{q-2}  \leq -\frac{q-2}{q-1}C^{-2}\quad\text{for almost all $t<t_*^h$},\\
\label{Eqn:EstimatesSDDerivationRayleighDecayInter}
\frac{{\rm d}}{{\rm d}t}\|\bv(t)\|_{\lq}^{q-2} \geq -\frac{q-2}{q-1}R(0) \quad\text{for almost all $t<t_*^h$}.
\end{gather}
When $q\in(1,2)$, integrating \eqref{Eqn:EstimatesSDDerivationRayleighDecayInter} between $0$ and $t<t_*^h$ yields $t_*^h=+\infty$, and the lower bound in \eqref{Eqn:Asymp1DiscreteSD}. Then integrating \eqref{Eqn:EstimatesSDDerivationSobolevIneqInter} between $0$ and $t$ yields the upper bound in \eqref{Eqn:Asymp1DiscreteSD}. When $q\in(2,2_\theta^*]$, integrating \eqref{Eqn:EstimatesSDDerivationSobolevIneqInter} between $0$ and $t<t_*^h$ yields $t_*^h<+\infty$ and the upper bound in \eqref{Eqn:Asymp1DiscreteSD}, and integrating \eqref{Eqn:EstimatesSDDerivationRayleighDecayInter} between $0$ and $t$ yields the lower bound in \eqref{Eqn:Asymp1DiscreteSD}. Still in the case $q\in (2,2_\theta^*)$, integrating \eqref{Eqn:EstimatesSDDerivationSobolevIneqInter} and \eqref{Eqn:EstimatesSDDerivationRayleighDecayInter} between $t_*^h$ and $t<t_*^h$ yields then \eqref{Eqn:Asymp2DiscreteSD}. 
\end{prf}

We next prove Proposition \ref{Prop:EstimatesFD}. The overall strategy is similar to the one of Proposition \ref{Prop:EstimatesSD}, but since the time derivative is now discrete, most of the computations must be carefully adapted. In particular, the chain rule is not anymore applicable to obtain a discrete version of \eqref{Eqn:EstimatesSDDerivationSobolevIneqInter}, and it is replaced with a convexity inequality, which induces correction terms that must be controlled. More precisely, the following estimation is proved. It can be compared to the upper bound in \eqref{Eqn:Asymp1DiscreteSD}, but involves a corrected coefficient in front of the time.
\begin{lem}\label{Lem:lemest}
Assume $q\in(1,2_\theta^*]\setminus\{2\}$ and $\bu^0\not\equiv 0$. Let $\bu:= (u_i^n)_{i\in\Z,n\in\N}$ be the solution of \eqref{Eqn:FD}. Then for any $m<n$ it holds
\begin{equation}\label{Eqn:Estinter}
\|\bu^n\|_{\lq}\leq \left[\|\bu^m\|_{\lq}^{q-2}-C^{-2}\frac{q-2}{q-1}(n-m)\Delta t\left(\frac{\|\bu^n\|_{\lq}^2}{\|\bu^m\|_{\lq}^2}\right)^{\frac{1}{n-m}}\right]^{\frac{1}{q-2}}.
\end{equation}
\end{lem}
\begin{prf}{Lemma \ref{Lem:lemest}}
Recall that, contrary to the continuous case, if $\bu^0\not\equiv 0$ then $\bu$ never vanishes, even if $q\in (2,2_\theta^*)$ (see the comments after Proposition \ref{Prop:EstimatesFD}). Therefore, the right-hand side of \eqref{Eqn:Estinter} has always a meaning.

The Poincaré-Sobolev inequality \eqref{Eqn:DiscretePoincareSobolev} and the energy inequality \eqref{Eqn:EnergyIneqFD1} yields 
\begin{equation}\label{Eqn:ODE}
\frac{1}{q'}\frac{\|\bu^{n+1}\|_{\lq}^q-\|\bu^n\|_{\lq}^q}{\Delta t}\leq -C^{-2}\|\bu^{n+1}\|_{\lq}^2,
\end{equation}
which is a finite-differences scheme for the differential inequality \eqref{Eqn:EstimatesSDDerivationSobolevIneq}. Once again, since $(\bu^n)_{n\in\N}$ never vanishes, \eqref{Eqn:ODE} implies $\|\bu^{n+1}\|_{\lq}<\|\bu^n\|_{\lq}$. We perform a discrete chain rule to deal with the left-hand side.  Let us write for all $n\geq 0$,
$$
\frac{\|\bu^{n+1}\|_q^q-\|\bu^n\|_q^q}{\Delta t} = \frac{ \left( \|\bu^{n+1}\|_{\lq}^{q-2} \right)^{\frac{q}{q-2}} - \left( \|\bu^{n}\|_{\lq}^{q-2} \right)^{\frac{q}{q-2}} }{ \|\bu^{n+1}\|_{\lq}^{q-2} - \|\bu^{n}\|_{\lq}^{q-2} }\frac{ \|\bu^{n+1}\|_{\lq}^{q-2} - \|\bu^{n}\|_{\lq}^{q-2} }{\Delta t}.
$$
$\bullet$ Case $1<q<2$:  Since $q-2<0$, $\|\bu^{n+1}\|_{\lq}^{q-2}> \|\bu^n\|_{\lq}^{q-2}$ (recall that $(\|\bu^n\|_{\lq})_{n\in\N}$ is decreasing by \eqref{Eqn:ODE}), and $(0,\infty)\ni x \mapsto x^{q/(q-2)}$ is convex. Therefore,
\begin{equation}\label{Eqn:e1}
\begin{aligned}
\frac{q}{q-2}\left(\|\bu^n\|_{\lq}^{q-2}\right)^{\frac{q}{q-2}-1}
\leq& \frac{ \left( \|\bu^{n+1}\|_{\lq}^{q-2} \right)^{\frac{q}{q-2}} - \left( \|\bu^{n}\|_{\lq}^{q-2} \right)^{\frac{q}{q-2}} }{ \|\bu^{n+1}\|_{\lq}^{q-2} - \|\bu^{n}\|_{\lq}^{q-2} } \\
\leq& \frac{q}{q-2}\left(\|\bu^{n+1}\|_{\lq}^{q-2}\right)^{\frac{q}{q-2}-1}.
\end{aligned}
\end{equation}
Hence, for any $n\geq 0$,
\begin{equation}\label{Eqn:e2}
\frac{q}{q'(q-2)}\|\bu^n\|_{\lq}^2\frac{ \|\bu^{n+1}\|_{\lq}^{q-2} - \|\bu^{n}\|_{\lq}^{q-2} }{\Delta t}+C^{-2}\|\bu^{n+1}\|_{\lq}^2\leq 0.
\end{equation}
Dividing by $\|u^n\|_{\lq}^2$ (which is always non-zero) and summing, we find that for all $m<n$,
\begin{equation}\label{Eqn:e3}
\|\bu^n\|_{\lq}^{q-2}\geq \|\bu^m\|_{\lq}^{q-2}+ C^{-2} \frac{2-q}{q-1}\Delta t\sum_{i=m}^{n-1}\frac{\|\bu^{i+1}\|_{\lq}^2}{\|\bu^i\|_{\lq}^2}.
\end{equation}
Moreover, by the arithmetico-geometric inequality,
$$
\frac{1}{n-m}\sum_{i=m}^{n-1}\frac{\|\bu^{i+1}\|_{\lq}^2}{\|\bu^i\|_{\lq}^2} \geq \left( \prod_{i=m}^{n-1}\frac{\|\bu^{i+1}\|_{\lq}^2}{\|\bu^i\|_{\lq}^2} \right)^{\frac{1}{n-m}}=\left(\frac{\|\bu^{n}\|_{\lq}^2}{\|\bu^m\|_{\lq}^2} \right)^{\frac{1}{n-m}}.
$$
Therefore,
$$
\|\bu^n\|_q^{q-2}\geq \|\bu^m\|_q^{q-2}+C^{-2}\frac{2-q}{q-1}(n-m)\Delta t \left(\frac{\|\bu^{n}\|_{\lq}^2}{\|\bu^m\|_{\lq}^2} \right)^{\frac{1}{n-m}},\quad 0\leq m<n<+\infty,
$$
which concludes the proof in case $1<q<2$.\\
\noindent
$\bullet$ Case $2<q\leq 2_\theta^*$: Computations run similarly. Equation \eqref{Eqn:e1} is reversed, but since $( \|\bu^{n+1}\|_{\lq}^{q-2} - \|\bu^{n}\|_{\lq}^{q-2} )/\Delta t$ is positive when $q>2$, \eqref{Eqn:e2} still holds. Then, since $(2-q)/(q-1)<0$, \eqref{Eqn:e3} is reversed, and we can still conclude using the arithmetico-geometric inequality.
\end{prf}

It remains to estimate the correction term $\left(\|\bu^n\|_{\lq}^2/\|\bu^m\|_{\lq}^2\right)^{1/(n-m)}$ in order to obtain Proposition \ref{Prop:EstimatesFD}.

\begin{prf}{Proposition \ref{Prop:EstimatesFD}}
We start with the extinction estimates \eqref{Eqn:ExtinctionEstimateFD}. Let us assume $q\in (2,2_\theta^*]$. From \eqref{Eqn:Estinter},
$$
\begin{aligned}
\|\bu^0\|_{\lq}^{q-2}\geq& \|\bu^m\|_{\lq}^{q-2}+C^{-2}\frac{q-2}{q-1}m\Delta t\left( \frac{\|\bu^m\|_{\lq}^2}{\|\bu^0\|_{\lq}^2} \right)^{1/m}\\
\geq& C^{-2}\frac{q-2}{q-1}m\Delta t\left( \frac{\|\bu^m\|_{\lq}^2}{\|\bu^0\|_{\lq}^2} \right)^{1/m}.
\end{aligned}
$$
Therefore for any $m\geq 0$,
$$
\frac{T_*^h}{m\Delta t}\geq \left( \frac{\|\bu^m\|_{\lq}}{\|\bu^0\|_{\lq}} \right)^{2/m},
$$
which yields the desired result.

We next prove \eqref{Eqn:Asymp1DiscreteFD}. This estimate could be proved by a consistency and stability analysis of the scheme associated to the ODE \eqref{Eqn:ODE}. However, the proof is rather technical since the stability constant of the scheme \eqref{Eqn:ODE} may blow up as we approach the extinction time in case $q\in(2,2_\theta^*)$. Instead, we use a method based on estimate \eqref{Eqn:Estinter}. 
Taking $m=0$ in \eqref{Eqn:Estinter}, with $n\geq 0$ fixed,
$$
\|\bu^n\|_{\lq}^{q-2}\geq \|\bu^0\|_{\lq}^{q-2}+C^{-2}\frac{2-q}{q-1}n\Delta t\left(\frac{\|\bu^n\|_{\lq}^2}{\|\bu^0\|_{\lq}^2}\right)^{1/n}.
$$
The overall strategy is to show that $\|\bu^n\|_{\lq}$ is estimated above by the unique positive solution to the equation
\begin{equation}\label{Eqn:Definition_a_n}
(a_n^{\Delta t,h})^{q-2}=\|\bu^0\|_{\lq}^{q-2}+C^{-2}\frac{2-q}{q-1}n\Delta t\left(\frac{(a_n^{\Delta t,h})^2}{\|\bu^0\|_{\lq}^2}\right)^{1/n}.
\end{equation}
Then, $a_n^{\Delta t,h}$ is compared with the desired expression
\begin{equation}\label{Eqn:Def_bn}
b_n^{\Delta t,h} := \left(\|\bu^0\|_{\lq}^{q-2} + C^{-2}\frac{2-q}{q-1}n\Delta t\right)^{\frac{1}{q-2}}.
\end{equation}

\noindent
$\bullet$ Case $1<q<2$:  Taking $m=0$ in \eqref{Eqn:Estinter}, with $n\geq 0$ fixed,
$$
\|\bu^n\|_{\lq}^{q-2}\geq \|\bu^0\|_{\lq}^{q-2}+C^{-2}\frac{2-q}{q-1}n\Delta t\left(\frac{\|\bu^n\|_{\lq}^2}{\|\bu^0\|_{\lq}^2}\right)^{1/n}.
$$
For any $\varepsilon >0$, if $\|\bu^n\|_{\lq}>\varepsilon$, then it holds
$$
\|\bu^n\|_{\lq}^{q-2}\geq \|\bu^0\|_{\lq}^{q-2}+C^{-2}\frac{2-q}{q-1}n\Delta t\left(\frac{\varepsilon^2}{\|\bu^0\|_{\lq}^2}\right)^{1/n}.
$$
On the other hand, if $\|\bu^n\|_{\lq}\leq\varepsilon$, then $\|\bu^n\|_{\lq}^{q-2}\geq \varepsilon^{q-2}$, since $q<2$. Therefore, for any $\varepsilon>0$,
$$
\|\bu^n\|_{\lq}^{q-2}\geq \min\left(\varepsilon^{q-2},\|\bu^0\|_{\lq}^{q-2}+C^{-2}\frac{2-q}{q-1}n\Delta t\left(\frac{\varepsilon^2}{\|\bu^0\|_{\lq}^2}\right)^{1/n}\right),
$$
hence
\begin{equation}\label{ineqeps}
\|\bu^n\|_{\lq}^{q-2}\geq \sup_{\varepsilon>0}\min\left(\varepsilon^{q-2},\|\bu^0\|_{\lq}^{q-2}+C^{-2}\frac{2-q}{q-1}n\Delta t\left(\frac{\varepsilon^2}{\|\bu^0\|_{\lq}^2}\right)^{1/n}\right).
\end{equation}
\details{On the one hand, $\varepsilon>0\mapsto \varepsilon^{q-2}$ is decreasing, and converging to $+\infty$ as $\varepsilon\rightarrow 0$, and to $0$ as $\varepsilon\rightarrow +\infty$. On the other hand, $\varepsilon>0\mapsto \|\bu^0\|_{\lq}^{q-2}+C^{-2}\frac{2-q}{q-1}n\Delta t(\varepsilon^2/\|\bu^0\|_{\lq}^2)^{1/n}$ is increasing, and it converges to zero as $\varepsilon\rightarrow 0$, and to $+\infty$ as $\varepsilon \rightarrow +\infty$.}
The supremum of the minimum in the right-hand side of \eqref{ineqeps} is achieved when $\varepsilon=a_n^{\Delta t,h}$, where $a_n^{\Delta t,h}$ is the unique positive solution to \eqref{Eqn:Definition_a_n}.
Then, \eqref{ineqeps} reads
\begin{equation}\label{Eqn:ineq_um_am}
    \forall n\in\N,\quad\|\bu^n\|_{\lq}\leq a_n^{\Delta t,h}.
\end{equation}
It remains to estimate the difference between $(a_n^{\Delta t,h})$ and $(b_n^{\Delta t,h})$. By \eqref{Eqn:Definition_a_n}, $a_n^{\Delta t,h}\leq \|\bu^0\|_{\lq}$, which then gives 
$$
(a_n^{\Delta t,h})^{q-2}\leq (b_n^{\Delta t,h})^{q-2}.
$$
Then, using \eqref{Eqn:Definition_a_n} again,
\begin{equation}\label{Eqn:comparaison_an_bn}
\|\bu^0\|_{\lq}^{q-2}+C^{-2}\frac{2-q}{q-1}n\Delta t\left(\frac{(b_n^{\Delta t,h})^2}{\|\bu^0\|_{\lq}^2}\right)^{1/n}\leq (a_n^{\Delta t,h})^{q-2} \leq (b_n^{\Delta t,h})^{q-2}.
\end{equation}
We deduce from \eqref{Eqn:Def_bn} and \eqref{Eqn:comparaison_an_bn}
\details{
$$
\begin{aligned}
\left| (b_n^{\Delta t,h})^{q-2} -(a_n^{\Delta t,h})^{q-2} \right| \leq& C^{-2}\frac{2-q}{q-1}n\Delta t \left( 1 - \left(\frac{(b_n^{\Delta t,h})^2}{\|\bu^0\|_{\lq}^2}\right)^{1/n}\right)\\
\leq& C^{-2}\frac{2-q}{q-1}n\Delta t \left( 1 - \exp\left(\frac{2}{n}\ln\left(\frac{b_n^{\Delta t,h}}{\|u^0\|_{\lq}}\right)
\right)\right)\\
\leq& C^{-2}\frac{2-q}{q-1}n\Delta t \left( 1 - \exp\left(\frac{2\Delta t}{q-2}\frac{\ln\left(1+c^hn\Delta t\right)}{n\Delta t}
\right)\right),
\end{aligned}
$$
}
\begin{equation}\label{Eqn:calcul}
\begin{aligned}
\left| (b_n^{\Delta t,h})^{q-2} -(a_n^{\Delta t,h})^{q-2} \right| \leq& C^{-2}\frac{2-q}{q-1}n\Delta t \left( 1 - \exp\left(\frac{2}{n}\ln\left(\frac{b_n^{\Delta t,h}}{\|\bu^0\|_{\lq}}\right)
\right)\right)\\
=& C^{-2}\frac{2-q}{q-1}n\Delta t \left( 1 - \exp\left(\frac{2\Delta t}{q-2}\frac{\ln\left(1+c^hn\Delta t\right)}{n\Delta t}
\right)\right),
\end{aligned}
\end{equation}
where $c^h:= C^{-2}\frac{2-q}{q-1}\|\bu^0\|_{\lq}^{2-q}$. Since $q-2<0$ and $\ln(1+ct)\leq ct$, we find
$$
\left| (b_n^{\Delta t,h})^{q-2} -(a_n^{\Delta t,h})^{q-2} \right|\leq C^{-2}\frac{2-q}{q-1} T\left(1 - \exp\left(\frac{2c^h\Delta t }{q-2}\right)\right),
$$
for any $n$ such that $n\Delta t\leq T$. Moreover, we assumed the initial data is such that the constant $c^h$ converges as $h\rightarrow 0$ \details{($\|\bu^0\|_{\lq}=\|\bu^0_h\|_{\lq}$ converges as $h\rightarrow 0$)}. Therefore, there exists $c\in\R_+$ such that $c^h\leq c$ for any $h>0$, and
\begin{equation}\label{Eqn:Est_bn_an}
\forall h>0,\quad\left| (b_n^{\Delta t,h})^{q-2} -(a_n^{\Delta t,h})^{q-2} \right|\leq C^{-2}\frac{2-q}{q-1} T\left(1 - \exp\left(\frac{2c\Delta t }{q-2}\right)\right),
\end{equation}
which converges to zero uniformly in $h$ as $\Delta t \rightarrow 0$. Since $x\in\R\mapsto x^{1/(q-2)}$ is Lipschitz continuous on $[\|\bu^0\|_{\lq}^{q-2}, +\infty)$, and $(b_n^{\Delta t,h})^{q-2},(a_n^{\Delta t,h})^{q-2}\in [\|\bu^0\|_{\lq}^{q-2}, +\infty)$ for any $n\in \N$,
$$
\forall h>0,\quad\left| b_n^{\Delta t,h} -a_n^{\Delta t,h} \right|\lessapprox C^{-2}\frac{2-q}{q-1} T\left(1 - \exp\left(\frac{2c\Delta t }{q-2}\right)\right) = \mathcal{O}(\Delta t),
$$
which converges to zero uniformly in $h$ as $\Delta t \rightarrow 0$. This ends the proof in case $q\in (1,2)$.

\noindent
$\bullet$ Case $2<q\leq 2_\theta^*$: Define again $a_n^{\Delta t,h}$ and $b_n^{\Delta t,h}$ by \eqref{Eqn:Definition_a_n} and \eqref{Eqn:Def_bn}. By the same computations as in the case $q\in(1,2)$, \eqref{Eqn:ineq_um_am} still holds. However, \eqref{Eqn:comparaison_an_bn} is reversed since now $q-2>0$:
\begin{equation}\label{Eqn:OrdreAnBn}
(b_n^{\Delta t,h})^{q-2}  \leq (a_n^{\Delta t,h})^{q-2} \leq         \|\bu^0\|_{\lq}^{q-2}-C^{-2}\frac{q-2}{q-1}n\Delta t\left(\frac{(b_n^{\Delta t,h})^2}{\|\bu^0\|_{\lq}^2}\right)^{1/n}.
\end{equation}
Computation \eqref{Eqn:calcul} yields in the present case
$$
\left| (b_n^{\Delta t,h})^{q-2} -(a_n^{\Delta t,h})^{q-2} \right| \leq
C^{-2}\frac{q-2}{q-1}n\Delta t \left( 1 - \exp\left(\frac{2\Delta t}{q-2}\frac{\ln\left(1-n\Delta t/T_*^h\right)}{n\Delta t}
\right)\right),
$$
where $T_*^h$ is defined in \eqref{Eqn:DefT*h}. Thanks to the assumptions on the initial data \details{($\|\bu^0\|_{\lq}=\|\bu^0_h\|_{\lq}$ converges as $h\rightarrow 0$)}, the times $T_*^h$ are bounded by a time $T_0$ independent of $h$. Therefore,
\begin{equation}\label{Eqn:Comparaison_nDt/T*h}
\left| (b_n^{\Delta t,h})^{q-2} -(a_n^{\Delta t,h})^{q-2} \right| \leq
T_0C^{-2}\frac{q-2}{q-1}\frac{n\Delta t}{T_*^h} \left( 1 - \exp\left(T_0\frac{2\Delta t}{q-2}\frac{\ln\left(1-n\Delta t/T_*^h\right)}{n\Delta t/T_*^h}
\right)\right).
\end{equation}
Inequality \eqref{Eqn:Comparaison_nDt/T*h} is the counterpart of \eqref{Eqn:calcul} in the present case. However, in contrast with the case $1<q<2$, the exponent in the exponential is not bounded uniformly in $n$. More precisely $\ln(1-n\Delta t/T_*^h)/(n\Delta t/T_*^h) \rightarrow -\infty$ as $n\Delta t/T_*^h \rightarrow 1$. To get around this difficulty, we obtain the desired estimate only for times $n$ such that the term $\ln(1-n\Delta t/T_*^h)/(n\Delta t/T_*^h)$, remains of finite order, and then use the decay of $(\|\bu^n\|_{\lq}^q)_n$ to conclude that the estimate remains valid for all subsequent times. More precisely, let $\varepsilon >0$ small enough such that there exists $\delta_\varepsilon^h<1$ with
\begin{equation}\label{Eqn:DefDeltaEpsh}
\left(C^{-2}\frac{q-2}{q-1}T_*^h(1-\delta_\varepsilon^h)\right)^{1/(q-2)}=\varepsilon.
\end{equation}
Definition \eqref{Eqn:DefDeltaEpsh} means
\begin{equation}\label{Eqn:bnPlusGrandEps}
\forall n\in\N\text{ s.t.~}n\Delta t<\delta_\varepsilon^h T_*^h,\quad
b_n^{\Delta t,h}\geq \varepsilon.
\end{equation}
Moreover, $\delta_\varepsilon^h<\tilde{\delta_\varepsilon}:=1-\frac{q-1}{C^{-2}(q-2)}\frac{\varepsilon^{q-2}}{T_0}<1$. As the function $\delta\in(0,1)\mapsto \ln(1-\delta)/\delta$) is non-increasing, \eqref{Eqn:Comparaison_nDt/T*h} yields that for all $n\in \N$ such that $n\Delta t<\delta_\varepsilon^h T_*^h$,
\begin{equation}\label{Eqn:Comparaison_a_n_b_n_finale}
\left| (b_n^{\Delta t,h})^{q-2} -(a_n^{\Delta t,h})^{q-2} \right| \leq
T_0C^{-2}\frac{q-2}{q-1}\tilde{\delta_\varepsilon} \left( 1 - \exp\left(T_0\frac{2\Delta t}{q-2}\frac{\ln\left(1-\tilde{\delta_\varepsilon}\right)}{\tilde{\delta_\varepsilon}}
\right)\right).
\end{equation}
Moreover, by the Lipschitz continuity of $x\mapsto x^{1/(q-2)}$ on $[\varepsilon ^{q-2}, \|\bu^0\|_{\lq}^{q-2}]$, and by \eqref{Eqn:OrdreAnBn} and \eqref{Eqn:bnPlusGrandEps}, it holds for all $n\in\N$ such that $n\Delta t<\delta_\varepsilon^h T_*^h$,
$$
\left| b_n^{\Delta t,h} -a_n^{\Delta t,h} \right| \lessapprox T_0C^{-2}\frac{q-2}{q-1}\tilde{\delta_\varepsilon} \left( 1 - \exp\left(T_0\frac{2\Delta t}{q-2}\frac{\ln\left(1-\tilde{\delta_\varepsilon}\right)}{\tilde{\delta_\varepsilon}}
\right)\right).
$$
The right-hand side converges to zero uniformly in $h$ as $\Delta t\rightarrow 0$. Therefore, there exists $\Delta t_0$ independent of $h>0$ such that for any $\Delta t<\Delta t_0$,
\begin{equation}\label{Eqn:Estimation_an_cas_q_plus_grand_que_2}
\forall n\in\N\text{ s.t.~}n\Delta t<\delta_\varepsilon^h T_*^h,\quad \left| b_n^{\Delta t,h} -a_n^{\Delta t,h} \right| \leq 2\varepsilon.
\end{equation}
In particular, for any $\Delta t <\Delta t_0$, \eqref{Eqn:ineq_um_am}, \eqref{Eqn:Estimation_an_cas_q_plus_grand_que_2} and the definition of $b_n^{\Delta t,h}$ \eqref{Eqn:Def_bn} imply
$$
\forall n\in\N \text{ s.t.~}n\Delta t < \delta_\varepsilon^h T_*^h,\quad\|\bu^n\|_{\lq}\leq \left(\|\bu^0\|_{\lq}^{-(2-q)}-\frac{q-2}{q-1}C^{-2}n\Delta t\right)^{\frac{1}{q-2}} + 2\varepsilon.
$$
Moreover, up to taking $\Delta t_0$ smaller, we can always assume that for any $\Delta t<\Delta t_0$ and for any $h>0$,
$$
\left(\|\bu^0\|_{\lq}^{q-2} + C^{-2}\frac{2-q}{q-1}n_\varepsilon^{h,\Delta t}\Delta t \right)^{\frac{1}{q-2}}  \leq 2\varepsilon,
$$
where $n_\varepsilon^{h,\Delta t}\Delta t$ is the largest time node below $\delta_\varepsilon^h T_*^h$. Thus, since $\|\bu^n\|_{\lq}$ is non-increasing, for any $\Delta t<\Delta t_0$,
\details{
$$
\begin{aligned}
\forall n\in\N \text{ s.t.~}n\Delta t \geq  \delta_\varepsilon^h T_*^h,\quad\|\bu^n\|_{\lq} \leq& 4\varepsilon \\
\leq& \left(\|\bu^0\|_{\lq}^{-(2-q)}-\frac{q-2}{q-1}C^{-2}n\Delta t\right)_+^{\frac{1}{q-2}} + 4\varepsilon.
\end{aligned}
$$
}
$$
\forall n\in\N \text{ s.t.~}n\Delta t \geq \delta_\varepsilon^h T_*^h,\quad\|\bu^n\|_{\lq} \leq \left(\|\bu^0\|_{\lq}^{-(2-q)}-\frac{q-2}{q-1}C^{-2}n\Delta t\right)_+^{\frac{1}{q-2}} + 4\varepsilon.
$$
This ends the proof in case $q\in (2,2_\theta^*]$.
\end{prf}

\section{Numerical results}
In this section, numerical tests are conducted for the fully discrete scheme \eqref{Eqn:FD}. Estimate \eqref{Eqn:Asymp1DiscreteFD} and the convergence of the scheme are tested, and the extinction time in the fast diffusion case is investigated numerically, as long as convergence to asymptotic profiles.
\subsection{Implementation}
The assembling of the matrix of the fractional Laplacian is done using formula \eqref{Eqn:DefinitionWeights}. Since the solution $(\bu^n)_{n\in\N}$ satisfies an homogeneous Dirichlet condition on nodes $x_i$ with $|i|>M_x$, the discrete fractional Laplacian is reduced to
$$
\forall i=-M_x,\dots,M_x,\quad \left[(-\Delta)^\theta_h \bu^n\right]_i = \sum_{i-M_x \leq j \leq i+M_x} w_j^h(u^n_i-u^n_{i-j}) + u_i^n\sum_{|j-i|>M_x} w_j^h.
$$
Therefore, the tail of the sum can be explicitly computed using formula \eqref{Eqn:SumWeights}. \details{The discrete fractional Laplacian is assembled into a square $M_x\times M_x$ matrix. The $(i,j)$-coefficient is $-w_{j-i}^h$ if $j\neq i$ and $\sum_{j\in\Z} w_j^h$ if $i=j$.} 
Then, Newton's method is used at each iteration to solve the equation in \eqref{Eqn:FD}. The stopping criterion for Newton's method is $\|\mathbf{v}_{k+1}-\mathbf{v}_{k}\|_\infty< 10^{-8}\|\bv_k\|_\infty$ where $\mathbf{v}_k$ and $\mathbf{v}_{k+1}$ are two consecutive iterations of Newton's method. A relative criterion is needed since an extinction phenomenon will occur when $q>2$ according to Proposition \ref{Prop:EstimatesFD}, hence an absolute criterion would not be enough near the extinction time.

\subsection{Energy decay and Rayleigh quotient}\label{Subsec:EnergyDecay}
\begin{figure}\label{Fig:DecayNorm}
    \begin{subfigure}[b]{0.49\textwidth}
         \centering
          \input{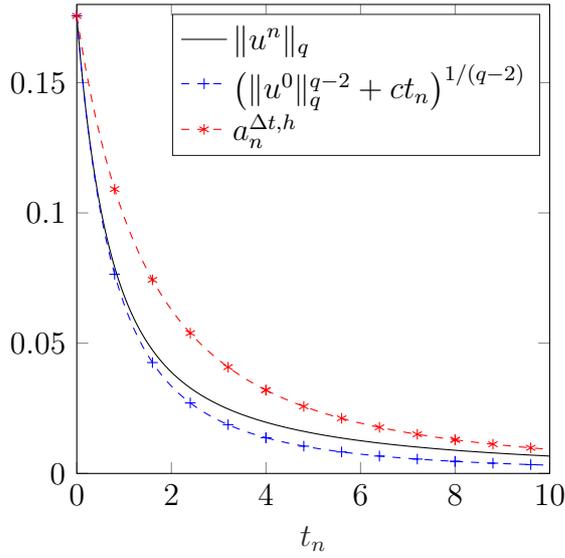}
         \caption{Porous medium case: $q=1.5$}
         \label{Fig:DecayNormPME}
     \end{subfigure}
     \begin{subfigure}[b]{0.49\textwidth}
         \centering
         \input{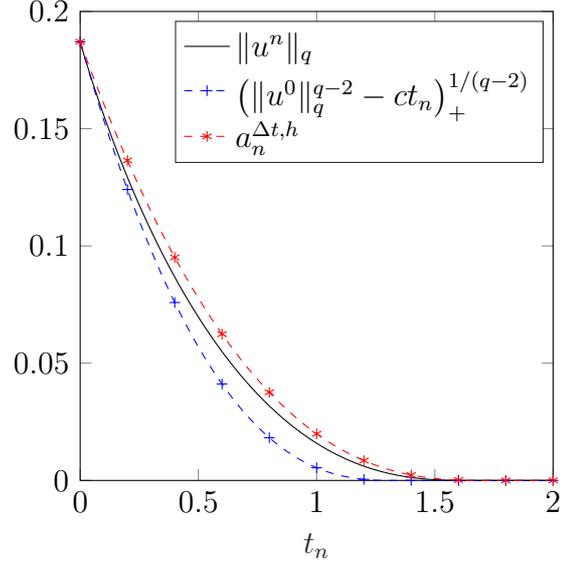}
         \caption{Fast diffusion case: $q=2.4$}
         \label{Fig:DecayNormFDE}
     \end{subfigure}
     \caption{Decay of $\|\bu^n\|_{\lq}$ for $L=1, \Delta t=0.001, h=0.05, \theta=0.25$ and $u^0$ given by \eqref{Eqn:DataIni}.}
\end{figure}

In this test, we illustrate the decay of $\lq$-norm indicated in Proposition \ref{Prop:EstimatesFD}. The domain is $(-L,L)$ with $L=1$, the space step is $h=0.05$, the time step is $\Delta t= 0.001$, the exponent of the fractionnal Laplacian is $\theta = 0.25$, and the initial datum is 
\begin{equation}\label{Eqn:DataIni}
u^0(x) =
\begin{cases}
    (1/2-x)(1/2+x)&\text{if $x\in[-1/2,1/2]$},\\
    0&\text{else}.
\end{cases}
\end{equation}
Figure \ref{Fig:DecayNormPME} deals with the case $q=1.5$ (porous medium case), and Figure \ref{Fig:DecayNormFDE} deals with the case $q=2.4$ (fast diffusion case). The norm $\|\bu^n\|_{\lq}$ and the upper estimate in the form of \ref{Eqn:ineq_um_am} are represented on Figures \ref{Fig:DecayNormPME} and \ref{Fig:DecayNormFDE} as functions of $t_n$. Although we only proved that it is satisfied in the semi-discrete case, the lower estimate in \eqref{Eqn:Asymp1DiscreteSD} is also represented on the figure. Remark that, the optimal (discrete) Poincare-Sobolev constant is required to compute the optimal upper bound \eqref{Eqn:ineq_um_am}. To compute this constant, we solve the following problem using Newton's method
\begin{equation}\label{Eqn:BestConstantEulerLagrange}
    \left\lbrace
    \begin{aligned}
    [(-\Delta)^\theta_h \mathbf{f}]_i&= f_i^{q-1} &&\text{$-M_x\leq i \leq M_x$},\\
    f_i&=0																					&&\text{$i>M_x$},\\
    f_i&>0 &&\text{$-M_x\leq i \leq M_x$},
    \end{aligned}
    \right.
\end{equation} 
and rescale the solution to satisfy $\|\mathbf{f}\|_{\lq}=1$. This yields a solution to the Euler-Lagrange equation for the minimization problem $\inf \{\|\mathbf{f}_{\Xh}:\:\mathbf{f}\in\lq,\, \|\mathbf{f}\|_{\lq}=1\}$, and the optimal Poincare-Sobolev constant is then given by $1/\|\mathbf{f}\|_{\Xh}$. 

According to Figures \ref{Fig:DecayNormPME} and \ref{Fig:DecayNormFDE}, observe that the $\lq$-norm of the solution satisfies the upper estimate \eqref{Eqn:Asymp1DiscreteFD}, as expected, but also the lower estimate in \eqref{Eqn:Asymp1Continue}, not proved in the fully discrete case. Moreover, the extinction phenomenon is almost observed in the fast diffusion case. Indeed, $\|\bu^n\|_{\lq}$ is of order $10^{-9}$ at time $t_n=1.63$, and is less than $10^{-300}$ at $t_n=1.65$ and for all subsequent times. However, as $\bu^n$ does not properly vanish, it is hard to precisely estimate the extinction time. Indeed, in the example of Figure \ref{Fig:DecayNormFDE}, since $\|\bu^n\|_{\lq}$ become of order $\Delta t$ at time $t_n = 1.41$, which leaves a wide range of possible extinction times.

We then investigate numerically the decay of the Rayleigh quotient \eqref{Eqn:DecayRayleigh}. Although it is not proved in the fully discrete case, this property seems to be satisfied. It is plotted as a funtion of time on Fig.~\ref{Fig:DecayRayleighPME} and \ref{Fig:DecayRayleighFDE}, respectivley in the porous medium case and the fast diffusion case. In the latter case, the quotient is well-defined as long as $\bu^n$ is nonzero.
\begin{figure}
    \begin{subfigure}[b]{0.49\textwidth}
         \centering
%
%
\begin{tikzpicture}

\begin{axis}[%
width=\textwidth,
height=\textwidth,
at={(0,0)},
xmin=0,
xmax=10,
xlabel style={font=\color{white!15!black}},
xlabel={$t_n$},
ymin=0.8,
ymax=1.6,
axis background/.style={fill=white},
title style={font=\bfseries},
legend style={legend cell align=left, align=left, draw=white!15!black}
]
\addplot [color=black]
  table[row sep=crcr]{%
0	1.53825824856705\\
0.032	1.53089718623812\\
0.0679999999999996	1.52310209767625\\
0.108000000000001	1.51492465466972\\
0.154999999999999	1.5058040337934\\
0.214	1.49484898977958\\
0.313000000000001	1.47701429110391\\
0.433999999999999	1.4550904299995\\
0.518000000000001	1.43938955154633\\
0.599	1.42378280215346\\
0.685	1.40674064203803\\
0.785	1.38643788638643\\
0.926	1.35729175577803\\
1.157	1.30954402286627\\
1.272	1.28627747960507\\
1.374	1.26610244469246\\
1.469	1.24777247140427\\
1.561	1.23048677856037\\
1.65	1.21422650594272\\
1.738	1.19861153429193\\
1.825	1.18363554347396\\
1.912	1.16912304171842\\
1.999	1.15507498625762\\
2.086	1.1414894983028\\
2.174	1.12821421855256\\
2.262	1.11540188384129\\
2.351	1.10290709496606\\
2.441	1.0907363729141\\
2.532	1.07889480566024\\
2.624	1.06738615947932\\
2.717	1.05621297850945\\
2.812	1.04526381545506\\
2.908	1.03466136680213\\
3.006	1.02430129350634\\
3.105	1.01429425476691\\
3.206	1.00454285243908\\
3.309	0.995056326028692\\
3.414	0.985842517263752\\
3.522	0.976826495987069\\
3.632	0.968102996007667\\
3.745	0.959602249210597\\
3.861	0.951338183473723\\
3.98	0.94332282131289\\
4.102	0.935566292397258\\
4.227	0.928076864020712\\
4.356	0.920806349764163\\
4.489	0.913769573350676\\
4.626	0.906978940423461\\
4.768	0.900399709986365\\
4.915	0.894048764364411\\
5.067	0.887940018664446\\
5.225	0.882048593456723\\
5.39	0.876357699644066\\
5.562	0.870888108600054\\
5.741	0.865656517452971\\
5.928	0.860650065968336\\
6.125	0.855837879704813\\
6.332	0.851244326895189\\
6.55	0.846869085617186\\
6.78	0.842714186951341\\
7.024	0.83876782756238\\
7.284	0.835025665557541\\
7.561	0.831500824485738\\
7.858	0.8281832245976\\
8.178	0.825070925867806\\
8.525	0.822159624806721\\
8.903	0.819452547321191\\
9.318	0.81694596978763\\
9.776	0.814644826737203\\
10	0.813669185801723\\
};
\addlegendentry{$R_n:=\frac{\|u^n\|_{\mathcal{X}_\theta(I)}^2}{\|u^n\|_q^2}$}

\end{axis}

\end{tikzpicture}%
         \caption{Porous medium case: $q=1.5$}
         \label{Fig:DecayRayleighPME}
     \end{subfigure}
     \begin{subfigure}[b]{0.49\textwidth}
         \centering
%
%
\begin{tikzpicture}

\begin{axis}[%
width=\textwidth,
height=\textwidth,
at={(0,0)},
xmin=0,
xmax=1.8,
xlabel style={font=\color{white!15!black}},
xlabel={$t_n$},
ymin=1.05,
ymax=1.4,
axis background/.style={fill=white},
title style={font=\bfseries},
legend style={legend cell align=left, align=left, draw=white!15!black}
]
\addplot [color=black]
  table[row sep=crcr]{%
0	1.35589772030381\\
0.00099999999999989	1.34883289676855\\
0.002	1.34442683246823\\
0.00299999999999989	1.34071520136119\\
0.00499999999999989	1.33438414280063\\
0.0069999999999999	1.32892885809014\\
0.0089999999999999	1.32404166000056\\
0.0109999999999999	1.31956826769163\\
0.014	1.31344192124639\\
0.0169999999999999	1.30785451322769\\
0.02	1.30268893400542\\
0.0229999999999999	1.29786690673657\\
0.0269999999999999	1.29187778939435\\
0.0309999999999999	1.28630697204621\\
0.0349999999999999	1.28108659312836\\
0.0389999999999999	1.27616598438765\\
0.044	1.27037797127608\\
0.0489999999999999	1.26493854127805\\
0.054	1.25980191647225\\
0.0589999999999999	1.25493178718575\\
0.0649999999999999	1.24939835880871\\
0.071	1.24416531309225\\
0.077	1.23919921483138\\
0.083	1.23447249795417\\
0.0900000000000001	1.22922997615178\\
0.097	1.22425134588296\\
0.104	1.21951094759961\\
0.112	1.21435733186184\\
0.12	1.20945908026059\\
0.128	1.20479287037812\\
0.137	1.19979597637559\\
0.146	1.19504313757435\\
0.155	1.19051341378469\\
0.165	1.18572009791056\\
0.175	1.18115835305337\\
0.185	1.17680949115561\\
0.196	1.17225219563603\\
0.207	1.16791376074341\\
0.219	1.1634110377104\\
0.231	1.15913008408438\\
0.243	1.15505429573097\\
0.256	1.15085342205356\\
0.269	1.14685968250827\\
0.283	1.14277369730539\\
0.297	1.1388949269272\\
0.312	1.13495282036113\\
0.327	1.131216539053\\
0.343	1.12744238190697\\
0.359	1.12387152559616\\
0.376	1.12028522568473\\
0.393	1.11689882135385\\
0.411	1.11351672716215\\
0.43	1.11015888819011\\
0.449	1.10700461605702\\
0.469	1.10388982326262\\
0.49	1.1008312664221\\
0.511	1.09797562804171\\
0.533	1.09518737326425\\
0.556	1.09248024808486\\
0.579	1.08997156645892\\
0.603	1.08755130159765\\
0.628	1.08523025355574\\
0.654	1.08301799196134\\
0.681	1.08092281616551\\
0.709	1.07895171648304\\
0.738	1.07711033709993\\
0.768	1.07540294188672\\
0.8	1.07378481058377\\
0.833	1.07231654901248\\
0.868	1.07096095639455\\
0.905	1.0697310202583\\
0.944	1.06863599608955\\
0.986	1.06766022408388\\
1.031	1.06681774323668\\
1.08	1.06610317558782\\
1.134	1.06551804204156\\
1.195	1.06505974683926\\
1.267	1.06472526022286\\
1.357	1.06451705979166\\
1.491	1.06442927722652\\
1.647	1.06442329731579\\
};
\addlegendentry{$R_n:=\frac{\|u^n\|_{\mathcal{X}_\theta(I)}^2}{\|u^n\|_q^2}$}

\end{axis}

\end{tikzpicture}%
         \caption{Fast diffusion case: $q=2.4$}
         \label{Fig:DecayRayleighFDE}
     \end{subfigure}
     \caption{Decay of $\|\bu^n\|_{\Xh}^2/\|\bu^n\|_{\lq}^2$ for $L=1, \Delta t=0.001, h=0.05, \theta=0.25$ and $u^0$ given by \eqref{Eqn:DataIni}.}
\end{figure}

\subsection{Convergence of the scheme}\label{Subsec:ConvergenceTest}
\begin{figure}\label{Fig:ConvergenceTemps}
    \begin{subfigure}[b]{0.49\textwidth}
         \centering
%
%
\definecolor{mycolor1}{rgb}{0.00000,0.44700,0.74100}%
\definecolor{mycolor2}{rgb}{0.85000,0.32500,0.09800}%
\definecolor{mycolor3}{rgb}{0.92900,0.69400,0.12500}%
\begin{tikzpicture}

\begin{axis}[%
width=\textwidth,
height=\textwidth,
xmode=log,
xmin=0.00048828125,
xmax=0.125,
xminorticks=true,
xlabel style={font=\color{white!15!black}},
xlabel={$\Delta t$},
ymode=log,
ymin=1e-05,
ymax=0.0107739414057614,
yminorticks=true,
axis background/.style={fill=white},
legend style={legend cell align=left, align=left, draw=white!15!black, legend pos = south east}
]
\addplot [color=red, dotted, mark=asterisk, mark options={solid}]
  table[row sep=crcr]{%
0.125	0.00334169020933973\\
0.0625	0.00172299868936746\\
0.03125	0.000875405067578189\\
0.015625	0.000440785563254218\\
0.00781249999999999	0.000220637935075487\\
0.00390625	0.000109845144253811\\
0.001953125	5.42693366626979e-05\\
0.000976562499999999	2.64360980183192e-05\\
0.00048828125	1.25081441410177e-05\\
};
\addlegendentry{$\theta=0.25$}

\addplot [color=green, dotted, mark=asterisk, mark options={solid}]
  table[row sep=crcr]{%
0.125	0.00592511208732056\\
0.0625	0.00315551966339615\\
0.03125	0.00162639761382163\\
0.015625	0.000825349851232684\\
0.00781249999999999	0.000414878383584949\\
0.00390625	0.000206993431311327\\
0.001953125	0.000102376542519428\\
0.000976562499999999	4.98976107926459e-05\\
0.00048828125	2.36152773815247e-05\\
};
\addlegendentry{$\theta=0.5$}

\addplot [color=blue, dotted, mark=asterisk, mark options={solid}]
  table[row sep=crcr]{%
0.125	0.0107739414057614\\
0.0625	0.00609052536224106\\
0.03125	0.0032546961259869\\
0.015625	0.00168278224699641\\
0.00781249999999999	0.000854925523608259\\
0.00390625	0.000428799476852325\\
0.001953125	0.000212651129936759\\
0.000976562499999999	0.000103786128800955\\
0.00048828125	4.91528941543712e-05\\
};
\addlegendentry{$\theta=0.75$}

\addplot [color=black, dashdotted]
  table[row sep=crcr]{%
0.125	0.00320208490010053\\
0.00048828125	1.25081441410177e-05\\
};
\addlegendentry{$\mathcal{O}(\Delta t)$}

\end{axis}

\end{tikzpicture}%
         \caption{Porous medium case: $q = 1.5$}
         \label{Fig:ConvergenceTemps_m}
     \end{subfigure}
     \begin{subfigure}[b]{0.49\textwidth}
         \centering
%

\begin{tikzpicture}

\begin{axis}[%
width=\textwidth,
height=\textwidth,
xmode=log,
xmin=0.00048828125,
xmax=0.125,
xminorticks=true,
xlabel style={font=\color{white!15!black}},
xlabel={$\Delta t$},
ymode=log,
ymin=1e-05,
ymax=0.0181977276151774,
yminorticks=true,
axis background/.style={fill=white},
legend style={legend cell align=left, align=left, draw=white!15!black, legend pos = south east}
]
\addplot [color=red, dotted, mark=asterisk, mark options={solid}]
  table[row sep=crcr]{%
0.125	0.00984080802649411\\
0.0625	0.005160076217624\\
0.03125	0.00264246363526175\\
0.015625	0.00133547301863066\\
0.00781249999999999	0.000669522604585152\\
0.00390625	0.000333410792337865\\
0.001953125	0.000164576451820404\\
0.000976562499999999	7.99681435491469e-05\\
0.00048828125	3.76178172781606e-05\\
};
\addlegendentry{$\theta=0.25$}

\addplot [color=green, dotted, mark=asterisk, mark options={solid}]
  table[row sep=crcr]{%
0.125	0.0130532812819356\\
0.0625	0.00754625159876874\\
0.03125	0.00415891994994712\\
0.015625	0.00220719503034087\\
0.00781249999999999	0.00114323233216355\\
0.00390625	0.000581511909919692\\
0.001953125	0.000293108528558414\\
0.000976562499999999	0.000147925993360128\\
0.00048828125	7.16422968446014e-05\\
};
\addlegendentry{$\theta=0.5$}

\addplot [color=blue, dotted, mark=asterisk, mark options={solid}]
  table[row sep=crcr]{%
0.125	0.0181977276151774\\
0.0625	0.0132258519559599\\
0.03125	0.00772532336232267\\
0.015625	0.00436800744464505\\
0.0078125	0.00246542073025575\\
0.00390625	0.0015070029048767\\
0.001953125	0.000923025257886115\\
0.0009765625	0.000515397158456707\\
0.00048828125	0.000275350064920228\\
};
\addlegendentry{$\theta=0.75$}

\addplot [color=black, dashdotted]
  table[row sep=crcr]{%
0.125	0.00963016122320913\\
0.00048828125	3.76178172781606e-05\\
};
\addlegendentry{$\mathcal{O}(\Delta t)$}

\end{axis}

\end{tikzpicture}%
         \caption{Fast diffusion case: $q=2.4$}
         \label{Fig:ConvergenceTemps_q}
     \end{subfigure}
     \caption{$\max_{n} \|u^{n}-u^{n}_{\Delta t,h}\|_{\ell^q_h(I)}$ for $T=2, h=0.05, L=1, \Delta t_{\rm ref}=5e-05$, and $u^0$ given by \eqref{Eqn:DataIni}.}
\end{figure}
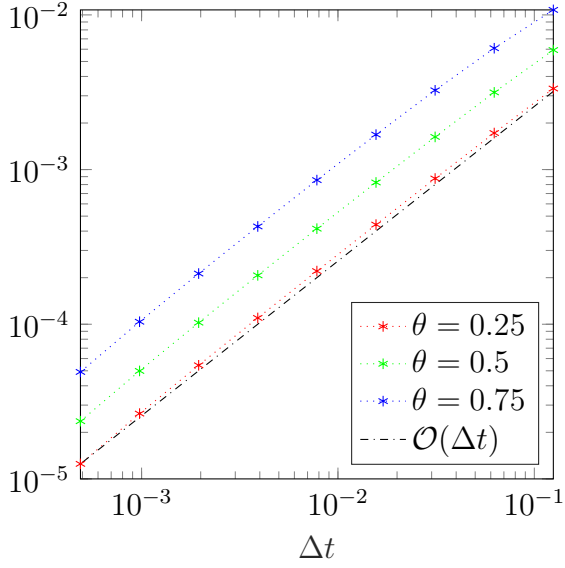
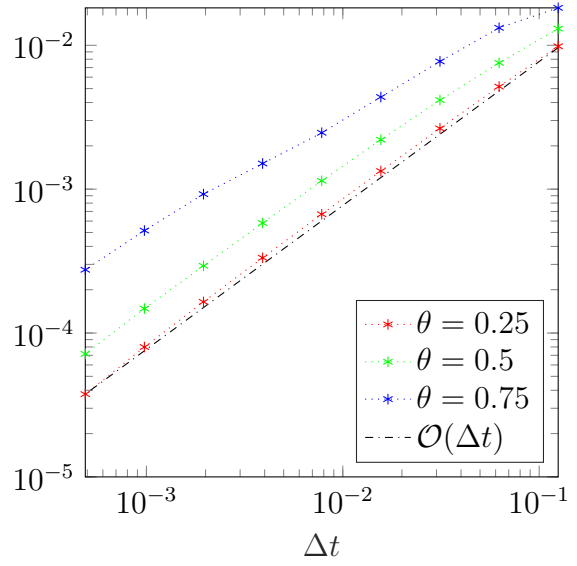

In this section, we test numerically the convergence of scheme \eqref{Eqn:FD}. We consider the convergence in $\Delta t$ and $h$ of the scheme for different values of $q$ and $\theta$. Let us first focus on the convergence in $\Delta t$ with $h$ fixed. For this test, we fix $h=0.05$, and the initial data \eqref{Eqn:DataIni}, and we compute a reference solution with scheme \eqref{Eqn:FD} using $\Delta t_{ref} = 5\cdot 10^{-5}$. The $\Delta t$-dependent error is defined as the $\ell^q_h(\R)$ norm of the difference between the solution given by the scheme and this reference solution. According to Fig.~\ref{Fig:ConvergenceTemps_m} and \ref{Fig:ConvergenceTemps_q}, the convergence is of order $1$ in time, whatever the value of q and $\theta$.

\begin{figure}\label{Fig:ConvergenceEspace}
    \begin{subfigure}[b]{0.49\textwidth}
         \centering
%

\begin{tikzpicture}

\begin{axis}[%
width=\textwidth,
height=\textwidth,
xmode=log,
xmin=0.00390625,
xmax=0.125,
xminorticks=true,
xlabel style={font=\color{white!15!black}},
xlabel={$h$},
ymode=log,
ymin=1e-05,
ymax=0.01,
yminorticks=true,
axis background/.style={fill=white},
legend style={legend cell align=left, align=left, draw=white!15!black, legend pos = south east}
]
\addplot [color=red, dotted, mark=asterisk, mark options={solid}]
  table[row sep=crcr]{%
0.125	0.00427543654321477\\
0.0625	0.00182583150807202\\
0.03125	0.000756110140489236\\
0.015625	0.000305190348907536\\
0.00781249999999999	0.000119637402084722\\
0.00390625	4.48472913362916e-05\\
};
\addlegendentry{$\theta=0.25$}

\addplot [color=green, dotted, mark=asterisk, mark options={solid}]
  table[row sep=crcr]{%
0.125	0.00669118949163915\\
0.0625	0.00258851818797903\\
0.03125	0.000953076841618256\\
0.015625	0.000338657080525933\\
0.00781249999999999	0.000116897337766876\\
0.00390625	3.90655114220723e-05\\
};
\addlegendentry{$\theta=0.5$}

\addplot [color=blue, dotted, mark=asterisk, mark options={solid}]
  table[row sep=crcr]{%
0.125	0.00725669649811189\\
0.0625	0.00254647984409968\\
0.03125	0.000849684876977852\\
0.015625	0.000288496950841345\\
0.00781249999999999	0.000100059076608204\\
0.00390625	3.42323473745418e-05\\
};
\addlegendentry{$\theta=0.75$}

\addplot [color=black, dashdotted]
  table[row sep=crcr]{%
0.125	0.00604888112391389\\
0.00390625	3.34156629789908e-05\\
};
\addlegendentry{$\mathcal{O}(h^{1.5})$}

\end{axis}

\end{tikzpicture}%
         \caption{Porous medium case: $q = 1.5$}
         \label{Fig:ConvergenceEspace_m}
     \end{subfigure}
     \begin{subfigure}[b]{0.49\textwidth}
         \centering
%

\begin{tikzpicture}

\begin{axis}[%
width=\textwidth,
height=\textwidth,
xmode=log,
xmin=0.00390625,
xmax=0.125,
xminorticks=true,
xlabel style={font=\color{white!15!black}},
xlabel={$h$},
ymode=log,
ymin=5.71829564356198e-05,
ymax=0.00224752162010634,
yminorticks=true,
axis background/.style={fill=white},
legend style={legend cell align=left, align=left, draw=white!15!black, legend pos= south east}
]
\addplot [color=red, dotted, mark=asterisk, mark options={solid}]
  table[row sep=crcr]{%
0.125	0.00155307298070456\\
0.0625	0.000721464861288836\\
0.03125	0.000457463651484958\\
0.015625	0.000291976544355227\\
0.0078125	0.000182600795908219\\
0.00390625	0.000109232471788087\\
};
\addlegendentry{$\theta=0.25$}

\addplot [color=green, dotted, mark=asterisk, mark options={solid}]
  table[row sep=crcr]{%
0.125	0.00215641269166382\\
0.0625	0.00114981885501016\\
0.03125	0.000648917596509939\\
0.015625	0.000363489767537407\\
0.0078125	0.000198094121165983\\
0.00390625	0.000102879889880587\\
};
\addlegendentry{$\theta=0.5$}

\addplot [color=blue, dotted, mark=asterisk, mark options={solid}]
  table[row sep=crcr]{%
0.125	0.00224752162010633\\
0.0625	0.000929840004387677\\
0.03125	0.000477635208981251\\
0.015625	0.000245227613696015\\
0.00781249999999999	0.000123073550328991\\
0.00390625	5.90701072861648e-05\\
};
\addlegendentry{$\theta=0.75$}

\addplot [color=black, dashdotted]
  table[row sep=crcr]{%
0.125	0.00182985460593983\\
0.00390625	5.71829564356198e-05\\
};
\addlegendentry{$\mathcal{O}(h)$}

\end{axis}

\end{tikzpicture}%
         \caption{Fast diffusion case: $q=2.4$}
         \label{Fig:ConvergenceEspace_q}
     \end{subfigure}
     \caption{$\max_{n} \|u^{n}-u^{n}_{\Delta t,h}\|_{\ell^q_h(I)}$ for $L=1, \Delta t =0.01, T=2, h_{\rm ref}=2^{-11}$ and $u^0$ given by \eqref{Eqn:DataIni}.}
\end{figure}
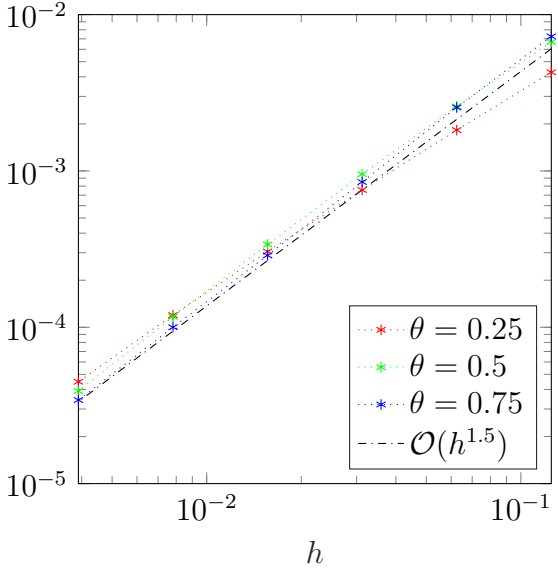
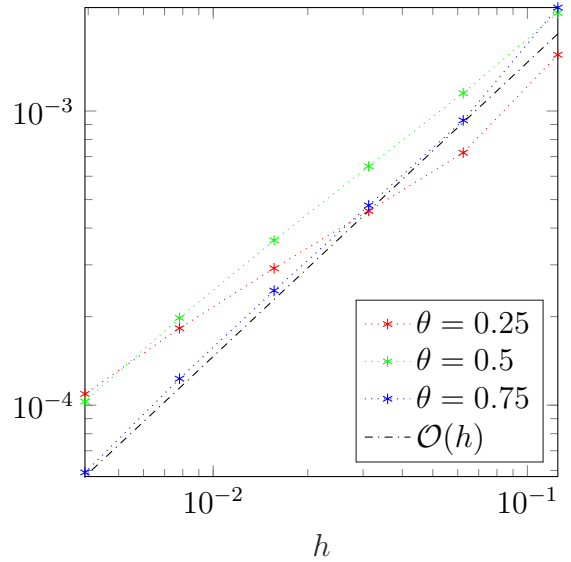

Next we investigate the convergence in space with a fixed time step, here $\Delta t=0.01$. The initial datum is again given by \eqref{Eqn:DataIni} and the reference solution is computed using a space step $h_{\rm ref}=5\cdot 10^{-4}$. Results are shown on figures \ref{Fig:ConvergenceEspace_m} and \eqref{Fig:ConvergenceEspace_q}, and indicates that the rate of convergence in $\lq$-norm is superlinear in the porous medium case, and linear in the fast diffusion case, except for small $\theta$ for which it is slightly slower. For small $\theta$, the rate is slower than the one given by Theorem \ref{Thm:ErreurDiscreteFracLap} for the convergence of the discrete fractional Laplacian, and this is probably due to the fact that the solution $u(x,t)$ of \eqref{Eqn:CDP} is not regular near the boundary of the domain.

\subsection{Numerical computation of the extinction time}\label{Subsec:ExtinctionTime}

As noted in Subsection \ref{Subsec:EnergyDecay}, numerical results for the scheme \eqref{Eqn:FD} in the fast diffusion case does not allow to precisely determine the extinction time even though an extinction phenomenon is observed. In this subsection, we suggest a different approach for the approximation of the extinction time. Consider the solution $(x,t)\in\Omega\times \R_+\mapsto u(x,t)$ of the continuous problem $\CDP$ and denote by $t_*\in\R_+^*$ its extinction time. Now suppose we know some approximation $t_0\in\R_+^*$ of $t_*$, and set 
\begin{equation}\label{Eqn:DefRescaledSolution}
v_{t_0}(x,s)=\left(1-\frac{t}{t_0}\right)_+^{-1/(q-2)}u(x,t),\quad s=t_0\log\left(\frac{t_0}{t_0-t}\right)
\end{equation}
In the new variable $s$, the approximate extinction time $t=t_0$ corresponds to $s=\infty$, and the equation satisfied by $v_{t_0}$ is
\begin{equation}\label{Eqn:RP}\tag{{\sf RP}}
\left\lbrace\begin{aligned}
\partial_t v_{t_0}^{q-1}+(-\Delta)^\theta v_{t_0}&= \frac{\lambda_q}{t_0}v_{t_0}^{q-1}&&\text{in $\Omega\times(0,+\infty)$},\\
v_{t_0}&=0																					&&\text{in $(\R^d\setminus\Omega)\times(0,+\infty)$},\\
v_{t_0}(\cdot,0)&=u_0																&&\text{in $\Omega$},
\end{aligned}\right. 
\end{equation}
where $\lambda_q:=(q-1)/(q-2)$. If $t_0=t_*$, one has by \eqref{Eqn:Asymp2Continue} 
\begin{equation}\label{Eqn:EstimateRescaled}
\forall s\in\R_+^*,\quad (t_0\lambda_q^{-1} C^{-2})^{1/(q-2)}\leq \|v_{t_*}(s)\|_{L^q(\Omega)}\leq \left(t_0\lambda_q^{-1}\frac{\|u_0\|_{H^\theta(\R^d)}^2}{\|u_0\|_{L^q(\Omega)}^2}\right)^{1/(q-2)},
\end{equation}
and in particular $v_{t_*}$ neither vanishes nor blow-up as $s\rightarrow +\infty$.
On the other hand, it follows from the definition of $v_{t_0}$ that, if $t_0>t_*$ then $v_{t_0}$ extinct in finite time again, and if $t_0<t_*$ then $\|v_{t_0}(s)\|_{L^q(\Omega)}\rightarrow\infty$ as $s\rightarrow \infty$. More precisely, if $t_0>t_*$, the upper bound in \eqref{Eqn:EstimateRescaled} is still satisfied, but the lower one is infringed at some time. In case $t_0<t_*$, the lower bound in \eqref{Eqn:EstimateRescaled} is still satisfied, but the upper one is infringed at some time. Consequently, our strategy to compute the extinction time is to introduce a numerical scheme for the rescaled problem $\eqref{Eqn:RP}$ and to perform a dichotomy on $t_0$, according to the above considerations.

The scheme for $\eqref{Eqn:RP}$ is the same as \eqref{Eqn:FD} where we simply add and treat implicitly the reaction term $(\lambda_q/t_0)v_{t_0}^{q-1}$. First the variable $s$ is discretized with a uniform step $\Delta s$ by setting $s_n:=n\Delta s$. Then, the scheme for \eqref{Eqn:RP} is
\begin{equation}\label{Eqn:RescaledScheme}
\left\{
\begin{aligned}
\frac{((v_{t_0})_i^{n+1})^{q-1}-((v_{t_0})_i^{n})^{q-1}}{\Delta s}+\left[(-\Delta)^\theta_h \bv_{t_0}^{n+1}\right]_i &= \frac{\lambda_q}{t_0}((v_{t_0})_i^{n+1})^{q-1}, && \text{$|i|\leq M_x$, $n\leq N_s$},\\
(v_{t_0})^n_i &= 0, && \text{$|i|\geq M_x+1$, $n\leq N_s$},\\
(v_{t_0})_i^0 &= (\bu^0)_i && \text{$|i|\leq M_x$},
\end{aligned}
\right.
\end{equation}
where $(v_{t_0})_i^n$ is an approximation of $v_{t_0}(x_i,s_n)$, $\bv_{t_0}^n=((v_{t_0})_i^n)_{i\in\Z}$, and $S_{t_0}:=N_s\Delta s$ is a final time to be determined later. We expect that there is a \emph{discrete} extinction time $t_*^{\Delta s,h}$ such that if $t_0<t_*^{\Delta s,h}$ then $\|\bv_{t_0}^n\|_{\lq}\rightarrow\infty$, and if $t_0>t_*^{\Delta s,h}$ then $\|\bv_{t_0}^n\|_{\lq}\rightarrow 0$, as $n\rightarrow \infty$. We could not prove this fact but it is observed numerically.

To initialize the dichotomy and determine whether the solution $(\bv^n)_{n\in\N}$ of \eqref{Eqn:RescaledScheme} satisfies $\|\bv^n\|_{\lq}\rightarrow +\infty$ or $\|\bv^n\|_{\lq}\rightarrow 0$, we draw on the continuous case, for which estimate \eqref{Eqn:EstimateRescaled} holds. From \eqref{Eqn:EstimateRescaled}, if $t_0=t_*$ it holds for $h$ small enough
\begin{equation}\label{Eqn:EstimateRescaledApprox}
   \forall s\in\R_+^*,\quad \frac{1}{2}(t_0\lambda_q^{-1} C^{-2})^{1/(q-2)} < \|(v_{t_0}(x_i,s))_{i\in\Z}\|_{\lq}< 2\left(t_0\lambda_q^{-1}\frac{\|\bu_0\|_{\Xh}^2}{\|\bu_0\|_{\lq}^2}\right)^{1/(q-2)}. 
\end{equation}
Taking $s=0$ in \eqref{Eqn:EstimateRescaledApprox} and recalling the definition of $v_{t_*}$ we obtain
\begin{equation}\label{Eqn:InitializationT1T2}
t_1^0:=\lambda_q\frac{\|\bu_0\|_{\lq}^2}{\|\bu_0\|_{\Xh}^2}\frac{(\|\bu^0\|_{\lq})^{q-2}}{2}\leq t_* \leq \lambda_q C^2(2\|\bu^0\|_{\lq})^{q-2}=:t_2^0.
\end{equation}
Therefore, to initialize the dichotomy, we let $t_1=t_1^0$, $t_2=t_2^0$, and we run the scheme \eqref{Eqn:RescaledScheme} on $t_0=(t_1+t_2)/2$. Then, if at some time $n$ we have \begin{equation}\label{Eqn:UpperEstimateDiscrete}
    \|\bv^n\|_{\lq}\geq 2\left(t_0\lambda_q^{-1}\frac{\|\bu_0\|_{\Xh}^2}{\|\bu_0\|_{\lq}^2}\right)^{1/(q-2)}, 
\end{equation}
then we consider the extinction time has been underestimated and we set $t_1\leftarrow t_0$. On the contrary, if at some time $n$ we have 
\begin{equation}\label{Eqn:LowerEstimateDiscrete}
\|\bv^n\|_{\lq}\leq\frac{1}{2}(t_0\lambda_q^{-1} C^{-2})^{1/(q-2)},
\end{equation}
then we consider the extinction time has been overestimated and we set $t_2\leftarrow t_0$. We then repeat the dichotomy with these new times $t_1$ and $t_2$.

The time $s_n$ at which either \eqref{Eqn:UpperEstimateDiscrete} or \eqref{Eqn:LowerEstimateDiscrete} occurs becomes larger as $t_0$ gets closer to the discrete extinction time. Therefore, we shall fix a final time $S_{t_0}$ large enough so that if neither \eqref{Eqn:UpperEstimateDiscrete} nor \eqref{Eqn:LowerEstimateDiscrete} occurs for $n\leq N_s$, this means that $t_0$ is close enough from $t_*^{\Delta s,h}$. To choose this time, we draw again on the continuous case, for which estimate \eqref{Eqn:EstimateRescaledApprox} holds. It is easily shown using the definition of $v_{t_0}$ that, by denoting $\mu_h:=[1 - \|\bu^0\|_{\lq}^2/(4^{q-2}C^2\|\bu^0\|_{\Xh}^2)]^{-1}\in[1,\infty)$, if $t_0\in [t_1^0,t_*)$ then the upper bound in \eqref{Eqn:EstimateRescaledApprox} will be broken at latest at time $s=t_0\log(t_0/((\mu_h-1)(t_*-t_0)))$ (i.e.~$t= t_0 - (\mu_h-1)(t_*-t_0)$, which is non-negative if $t_0\geq t_1^0$), and if $t_0\in (t_*,t_2^0]$ then the lower bound in \eqref{Eqn:EstimateRescaledApprox} will be broken at latest at time $s=t_0\log(t_0/(\mu_h(t_0-t_*)))$ (i.e.~$t=t_0-\mu_h(t_0-t_*)$, which is non-negative if $t_0\leq t_2^0$). Therefore, by setting 
\begin{equation}\label{Eqn:DefinitionS}
    S_{t_0}=t_0\log(t_0/((\mu_h-1)\epsilon))
\end{equation}
(which is non-negative if $\epsilon<t_2^0-t_1^0$ and $t_0\geq t_1^0$), if \eqref{Eqn:EstimateRescaledApprox} holds for $s\in[0,S_{t_0}]$ then it holds $|t_0-t_*|<\epsilon$. Consequently, at each step of the dichotomy we define $S_{t_0}$ by \eqref{Eqn:DefinitionS} and let $N_s=S_{t_0}/\Delta S$ in the scheme \eqref{Eqn:RescaledScheme} to perform the current iteration. If neither \eqref{Eqn:UpperEstimateDiscrete} nor \eqref{Eqn:LowerEstimateDiscrete} occurs by the final time $N_s$, it is expected that $t_0$ is close to $t_*^{\Delta s,h}$, so we stop the dichotomy and return $t_0$ for the extinction time. If this situation does not happen, the dichotomy is terminated when $|t_1-t_2|<\epsilon$ and we then return $t_0=(t_1+t_2)/2$ for the approximated extinction time.

\begin{figure}\label{Fig:ConvergenceExtinctionTime}
        \begin{subfigure}[b]{0.49\textwidth}
         \centering
%

\begin{tikzpicture}

\begin{axis}[%
width=\textwidth,
height=\textwidth,
at={(0,0)},
xmode=log,
xmin=0.00504135701506484,
xmax=0.171467764060357,
xminorticks=true,
ymode=log,
ymin=5.04135701506484e-05,
ymax=0.00183217054747775,
yminorticks=true,
axis background/.style={fill=white},
title style={font=\bfseries},
legend style={at={(0.97,0.03)}, anchor=south east, legend cell align=left, align=left, draw=white!15!black}
]
\addplot [color=red, dotted, mark=asterisk, mark options={solid}]
  table[row sep=crcr]{%
0.171467764060357	0.00183217054747775\\
0.0952598689224204	0.00113669980099607\\
0.0529221494013446	0.000683732397916525\\
0.0294011941118581	0.000400603712189084\\
0.0163339967288101	0.000229204106655611\\
0.0090744426271167	0.000127938845278841\\
0.00504135701506484	6.91884224459649e-05\\
};
\addlegendentry{$\alpha=0.5$}

\addplot [color=green, dotted, mark=asterisk, mark options={solid}]
  table[row sep=crcr]{%
0.171467764060357	0.00168549546055718\\
0.0952598689224204	0.00115143927399908\\
0.0529221494013446	0.000745006487290256\\
0.0294011941118581	0.000460728412486233\\
0.0163339967288101	0.0002742055087952\\
0.0090744426271167	0.00015746883711043\\
0.00504135701506484	8.68887020373243e-05\\
};
\addlegendentry{$\alpha=1$}

\addplot [color=blue, dotted, mark=asterisk, mark options={solid}]
  table[row sep=crcr]{%
0.171467764060357	0.000972497926985172\\
0.0952598689224204	0.000744921606361082\\
0.0529221494013446	0.000531706910793406\\
0.0294011941118581	0.000355606804005926\\
0.0163339967288101	0.000224624893208736\\
0.0090744426271167	0.000134809948465042\\
0.00504135701506484	7.68201054941464e-05\\
};
\addlegendentry{$\alpha=1.5$}

\addplot [color=black, dashdotted]
  table[row sep=crcr]{%
0.171467764060357	0.00171467764060357\\
0.00504135701506484	5.04135701506484e-05\\
};
\addlegendentry{$0.01\Delta s$}

\end{axis}

\end{tikzpicture}%
          \caption{}
         \label{Fig:ConvergenceTimeExtinctionTime}
     \end{subfigure}
     \begin{subfigure}[b]{0.49\textwidth}
         \centering
%

\begin{tikzpicture}

\begin{axis}[%
width=\textwidth,
height=\textwidth,
at={(0,0)},
xmode=log,
xmin=0.0175438596491228,
xmax=0.333333333333333,
xminorticks=true,
xtick distance = 2,
ymode=log,
ymin=0.001,
ymax=0.1,
yminorticks=true,
axis background/.style={fill=white},
title style={font=\bfseries},
legend style={at={(0.97,0.03)}, anchor=south east, legend cell align=left, align=left, draw=white!15!black}
]
\addplot [color=red, dotted, mark=asterisk, mark options={solid}]
  table[row sep=crcr]{%
0.333333333333333	0.0392995012474411\\
0.2	0.0277801965977846\\
0.142857142857143	0.0207988086980095\\
0.0909090909090909	0.0137217760695307\\
0.0588235294117647	0.00902869347305679\\
0.04	0.00616869172825174\\
0.0263157894736842	0.00410245747373405\\
0.0175438596491228	0.00265790933627463\\
};
\addlegendentry{$\alpha=0.5$}

\addplot [color=green, dotted, mark=asterisk, mark options={solid}]
  table[row sep=crcr]{%
0.333333333333333	0.0699061021311204\\
0.2	0.0478587689059693\\
0.142857142857143	0.0356538408234677\\
0.0909090909090909	0.0234041363017226\\
0.0588235294117647	0.0153509452893492\\
0.04	0.0104669163366347\\
0.0263157894736842	0.00690586530880743\\
0.0175438596491228	0.00449524087803321\\
};
\addlegendentry{$\alpha=1$}

\addplot [color=blue, dotted, mark=asterisk, mark options={solid}]
  table[row sep=crcr]{%
0.333333333333333	0.0392324083321323\\
0.2	0.0301507115044324\\
0.142857142857143	0.0231573554496649\\
0.0909090909090909	0.0155069209167432\\
0.0588235294117647	0.0102593841569667\\
0.04	0.0070208014808757\\
0.0263157894736842	0.00464469082464358\\
0.0175438596491228	0.00302320479466589\\
};
\addlegendentry{$\alpha=1.5$}

\addplot [color=black, dashdotted]
  table[row sep=crcr]{%
0.333333333333333	0.0333333333333333\\
0.0175438596491228	0.00175438596491228\\
};
\addlegendentry{$0.1 h$}

\end{axis}
\end{tikzpicture}%
         \caption{}
         \label{Fig:ConvergenceSpaceExtinctionTime}
     \end{subfigure}
         \caption{Convergence of the extinction time computed with the scheme described in subsection \ref{Subsec:ExtinctionTime}, for $q=2.4,\,L=1,\,\epsilon=0.5\cdot 10^{-8}$. {\scriptsize \sc (A)}: Convergence in $\Delta s$ with $h=0.1$ and reference extinction time computed with $\Delta s_{\rm ref}=5\cdot10^{-4}$. {\scriptsize\sc (B)}: Convergence in $h$ with $\Delta s=0.1$ and reference extinction time computed with $h_{\rm ref}=1\cdot 10^{-3}$.}
\end{figure}
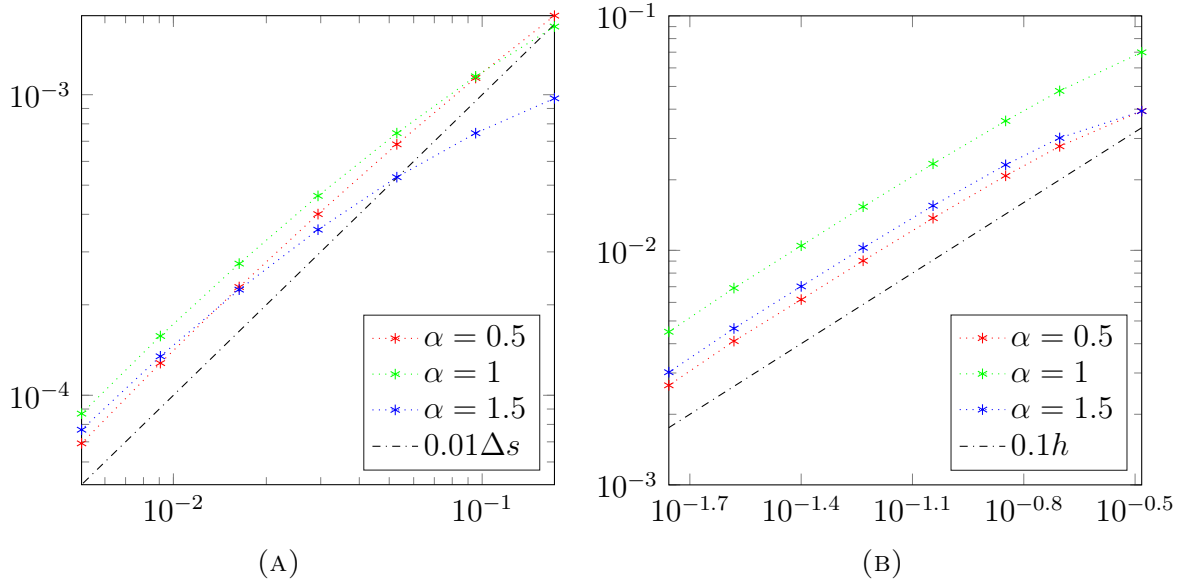

Figure \ref{Fig:ConvergenceTimeExtinctionTime} and \ref{Fig:ConvergenceSpaceExtinctionTime} show the results of convergence tests for the approximation of extinction time. These tests are again performed with the initial datum \eqref{Eqn:DataIni}. On Figure \ref{Fig:ConvergenceTimeExtinctionTime} we fix the space step $h=0.1$ and plot the convergence of the extinction time computed with the above scheme as $\Delta s\rightarrow 0$. The reference extinction time is computed by taking a small $\Delta s_{\rm ref}=0.0005$ and same $h$. The order of convergence is almost $\mathcal{O}(\Delta s)$, and we also point out that the constant is small ($0.01$ on this test), which means that it is enough to take a $\Delta s$ much larger than $\epsilon$ to compute the extinction time with a precision of $\mathcal{O}(\epsilon)$.  On figure \ref{Fig:ConvergenceSpaceExtinctionTime} we fix the time step $\Delta s=0.1$ and we plot the convergence of the extinction time as $h\rightarrow 0$. The reference extinction time is computed by taking a small $h_{\rm ref}=0.005$ and same $\Delta s$. Again the order of convergence is $1$.

\details{
\begin{prf}{the computations}
For any $t_0>0$, \eqref{Eqn:EstimateRescaledApprox} yields, in the original variable $t$.
\begin{multline*}
\frac{1}{2}\left(t_0\lambda_q^{-1}C^{-2}\frac{(t_*-t)_+}{t_0-t}\right)^{1/(q-2)} \leq \left\|\left(v_{t_0}\left(x_i,\lambda_q^{-1}\log\left(\frac{t_0}{t_0-t}\right)\right)\right)_{i\in\Z}\right\|_{\lq} \\
\leq  2\left(t_0\lambda_q^{-1}\frac{\|\bu^0\|^2_{\Xh}}{\|\bu^0\|^2_{\lq}}\frac{(t_*-t)_+}{t_0-t}\right)^{1/(q-2)},
\end{multline*}
for any $t<t_0$. Taking $t=0$ directly yields \eqref{Eqn:InitializationT1T2}. Moreover, if $t_0<t_*$, then $t\in[0,t_0)\mapsto ((t_*-t)/(t_0-t))^{1/(q-2)}$ is non-decreasing, but diverges to $+\infty$ as $t\rightarrow t_0$. Therefore the lower bound in \eqref{Eqn:EstimateRescaledApprox} is satisfied for $t\in[0,t_0)$, but the upper bound in \eqref{Eqn:EstimateRescaledApprox} is broken at latest at $s=t_0\log(t_0/(t_0-t))$ with $t\in(0,t_0)$ the smallest non-negative time such that
$$
2\left(t_0\lambda_q^{-1}\frac{\|\bu_0\|_{\Xh}^2}{\|\bu_0\|_{\lq}^2}\right)^{1/(q-2)} \leq \frac{1}{2}\left(t_0\lambda_q^{-1}C^{-2}\frac{(t_*-t)_+}{t_0-t}\right)^{1/(q-2)},
$$
that is
\begin{gather*}
t_0-t\leq \frac{\|\bu^0\|_{\lq}^2}{4^{q-2}C^2\|\bu^0\|_{\Xh}^2}(t_*-t),\\
t_0 - \frac{\|\bu^0\|_{\lq}^2}{4^{q-2}C^2\|\bu^0\|_{\Xh}^2}t_* \leq \left(1-\frac{\|\bu^0\|_{\lq}^2}{4^{q-2}C^2\|\bu^0\|_{\Xh}^2}\right)t,\\
\left(1-\frac{\|\bu^0\|_{\lq}^2}{4^{q-2}C^2\|\bu^0\|_{\Xh}^2}\right)t_*-(t_*-t_0)\leq \left(1-\frac{\|\bu^0\|_{\lq}^2}{4^{q-2}C^2\|\bu^0\|_{\Xh}^2}\right) t,\\
t_* - \mu_h (t_*-t_0)\leq t,\\
t_0 - (\mu_h-1)(t_*-t_0)\leq t.
\end{gather*}
Moreover, by a direct computation we verify that $\mu_h(t_2^0-t_1^0)=t_2^0$. Therefore, if $t_0\in [t_1^0,t_*)$, then by \eqref{Eqn:InitializationT1T2} it holds $t_0 - (\mu_h-1)(t_*-t_0)\in[0,t_0]$, and the upper bound in \eqref{Eqn:EstimateRescaledApprox} will be broken at latest at $s=t_0\log(t_0/(t_0-t))$ with $t=t_0 - (\mu_h-1)(t_*-t_0)$.\\
Next if $t_0>t_*$, then $t\in[0,t_0)\mapsto ((t_*-t)_+/(t_0-t))^{1/(q-2)}$ is non-increasing, but is identically $0$ from $t=t_*$. Therefore the upper bound in \eqref{Eqn:EstimateRescaledApprox} is satisfied for $t\in[0,t_0)$, but the lower bound in \eqref{Eqn:EstimateRescaledApprox} is broken at latest at $s=t_0\log(t_0/(t_0-t))\geq 0$ with $t$ the smallest non-negative time such that
$$
2\left(t_0\lambda_q^{-1}\frac{\|\bu^0\|^2_{\Xh}}{\|\bu^0\|^2_{\lq}}\frac{(t_*-t)_+}{t_0-t}\right)^{1/(q-2)} \leq \frac{1}{2}(t_0\lambda_q^{-1} C^{-2})^{1/(q-2)} ,
$$
that is
\begin{gather*}
    t_*-t\leq \frac{\|\bu^0\|_{\lq}^2}{4^{q-2}C^2\|\bu^0\|_{\Xh}^2}(t_0-t),\\
    t_* - \frac{\|\bu^0\|_{\lq}^2}{4^{q-2}C^2\|\bu^0\|_{\Xh}^2} t_0 \leq \left(1 - \frac{\|\bu^0\|_{\lq}^2}{4^{q-2}C^2\|\bu^0\|_{\Xh}^2} \right)t,\\
    \left(1 - \frac{\|\bu^0\|_{\lq}^2}{4^{q-2}C^2\|\bu^0\|_{\Xh}^2} \right)t_0 -(t_0-t_*)\leq \left(1 - \frac{\|\bu^0\|_{\lq}^2}{4^{q-2}C^2\|\bu^0\|_{\Xh}^2} \right)t,\\
    t_0 - \mu_h(t_0-t_*)\leq t.
\end{gather*}
Moreover, since $\mu_h(t_2^0-t_1^0)=t_2^0$, we observe that $t_0\leq t_2^0$ implies $t_0-\mu_h(t_0-t_*)=t_*-(\mu_h-1)(t_0-t_*)\geq t_*-(\mu_h-1)(t_2^0-t_1^0)=t_*-t_1^0\geq 0$. Therefore, if $t_0\in(t_*,t_2^0]$, then it holds $t_0-\mu_h(t_0-t_*)\in[0,t_0]$, and the lower bound in \eqref{Eqn:EstimateRescaledApprox} will be broken at latest at $s=t_0\log(t_0/(t_0-t))\geq 0$ with $t=t_0 - \mu_h(t_*-t_0)$.\\
Now let $t_0\in [t_1^0,t_2^0]$, $\epsilon < t_2^0-t_1^0$. If $\epsilon< |t_0-t_*|$, then either the upper bound or the lower bound in \eqref{Eqn:EstimateRescaledApprox} will be broken by $S_{t_0}:=t_0\max\{ \log( t_0/((\mu_h-1)\epsilon) ), \log( t_0/(\mu_h\epsilon))\}= t_0\log(t_0/((\mu_h-1)\epsilon))$. Note that $S_{t_0}>t_0\log(t_0/((\mu_h-1)(t_2^0-t_1^0))=t_0\log(t_0/t_1^0)\geq 0$. If this is not the case then it means $|t_0-t_*|>\epsilon$.
\end{prf}
}


\subsection{Convergence to asymptotic profiles in the fast diffusion case}

We assume in this subsection $u_0\geq 0$, which yields in particular $u(t)\geq 0$ for any $t\geq 0$ \cite{akagi_existence}. When $t_0=t_*$, the rescaled solution $v_{t_*}$ defined in \eqref{Eqn:DefRescaledSolution} is bounded both above and below by positive constants as given in \eqref{Eqn:EstimateRescaled}. This justifies investigating the convergence of $v_{t_*}(s)$ as $s\rightarrow +\infty$. 

In the fractional porous medium case ($\theta \in (0,1)$, $q\in (1,2)$), it has been proved in \cite{bonforte_sire_vazquez, bonforte_figalli_vazquez, franzina_volzone} that $v_{t_*}(s)$ converges towards an asymptotic profile as $s\rightarrow +\infty$, with an exponential rate. The asymptotic profile is a solution to a sublinear fractional Lane-Emden equation \eqref{Eqn:Fractional_Lane_Emden_Fowler}. 

For the standard subcritical fast diffusion case ($\theta=1$, $q\in (2,2d/(d-2)_+)$), we refer to the survey \cite[Section 9.1]{bonforte_figalli_survey}. In this case, quasi-convergence (i.e.~ convergence along a subsequence) was first shown in \cite{berryman_holland}, and then convergence (along the whole sequence) was obtained in \cite{feireisl_simondon}. More precisely, there exists $\phi\in C(\overline{\Omega})$ such that $v_{t_*}(s)\rightarrow \phi$ in $C(\overline{\Omega})$ as $s\rightarrow +\infty$. The asymptotic profile $\phi$ is a solution of the Lane-Emden-Fowler equation
\begin{equation}\label{Eqn:Lane_Emden_Fowler}
    \left\lbrace
    \begin{aligned}
    -\Delta \phi&=\frac{\lambda_q}{t_*}\phi^{q-1} &&\text{on $\Omega$},\\
    \phi&=0																					&&\text{in $\partial\Omega$},\\
    \phi&>0 &&\text{in $\Omega$}.
    \end{aligned}
    \right. 
\end{equation}
It can be shown that $\phi$ is bounded (see \cite{bonforte_figalli_survey} and references therein). Consequently, by elliptic regularity, $\phi$ belongs to $C^\infty(\Omega) \cap C^\alpha(\overline{\Omega})$ for some $\alpha \in (0,1)$. The convergence of $v_{t_*}(s)$ to $\phi$ also holds in relative error, as proved in \cite{bonforte_grillo_vazquez}. 

Then, the question of the rate of convergence of $v_{t_*}$ to $\phi$ was first solved in \cite{bonforte_figalli} under the assumption that $\phi$ is \emph{non-degenerate}. Under this assumption, it is shown that the rate of convergence is given by the linearized problem 
\begin{equation}\label{Eqn:Linearized_rescaled_problem}
    \left\lbrace\begin{aligned}
    (q-1)\phi^{q-2}\partial_s h &= - \mathcal{L}_{\phi} h &&\text{in $\Omega\times(0,+\infty)$},\\
    h&=0																					&&\text{in $\partial\Omega\times(0,+\infty)$},
\end{aligned}\right.  
\end{equation}
where $\mathcal{L}_{\phi}h = -\Delta h -(q-1)(\lambda_q/t_*)\phi^{q-2}h$, and $h(s)\approx v_{t_*}(s)-\phi$. The asymptotic profile $\phi$ is said to be \emph{non-degenerate} if the linear problem
\begin{equation}\label{Eqn:Linearized_Lane_Emden_Fowler}
    \left\lbrace
    \begin{aligned}
    \mathcal{L}_{\phi}h&=0 &&\text{on $\Omega$},\\
    h&=0																					&&\text{in $\partial\Omega$},
    \end{aligned}
    \right. 
\end{equation}
admits no non-trivial solution. In other words, $\phi$ is non-degenerate if $0$ is not an eigenvalue of the weighted eigenvalue problem
\begin{equation}\label{Eqn:weighted_eigenvalue_problem}
    \left\lbrace
    \begin{aligned}
    \mathcal{L}_{\phi}h&=\nu \phi^{q-2} h &&\text{on $\Omega$},\\
    h&=0								    &&\text{in $\partial\Omega$}.
    \end{aligned}
    \right. 
\end{equation}
Since the operator $(1/\phi^{q-2})\mathcal{L}_{\phi}$ is self-adjoint in the weighted $L^2$ space $L^2(\Omega;\phi^{q-2}\d x)$, it can be shown by spectral theory that there is a complete orthonormal system $\{e_j\}_{j=1}^\infty$ of $L^2(\Omega;\phi^{q-2}\d x)$ associated with eigenvalues $\nu_1<\cdots<\nu_k<\nu_{k+1}\rightarrow +\infty$ (see \cite{bonforte_figalli}). Moreover, the first eigenvalue is $\nu_1 = -(q-1)/t_* < 0$, it is simple and the corresponding eigenfunction is $\phi$. If $\phi$ is non-degenerate, then there exists an integer $k_q$ so that $\nu_{k_q}< 0 <\nu_{k_q+1}$. There are no general results that allow to systematically determine whether an asymptotic profile is non-degenerate. For a discussion on this topic, we refer again to the survey \cite[Section 9.1]{bonforte_figalli_survey}. However, we mention a result from \cite[Theorem 4.2]{damascelli_grossi_pacella}, which states that if $\Omega$ is a ball in $\mathbb{R}^d$, then all solutions to the Lane-Emden-Fowler equation are non-degenerate.

Since $v_{t_*}(s)$ converges to $\phi$ as $s\rightarrow+\infty$, we may expect \eqref{Eqn:Linearized_rescaled_problem} to suitably describe the evolution of $v_{t_*}(s)-\phi$ for large $s$. However, since $v_{t_*}(s)$ is bounded according to \eqref{Eqn:EstimateRescaled}, modes corresponding to the negative eigenvalues $\nu_1,\dots,\nu_{k_q}$ should not be active. The result of \cite{bonforte_figalli} is that, if $\phi$ is non-degenerate, then $v_{t_*}(s)$ converges exponentially to $\phi$ as $s\rightarrow +\infty$, with rate at least $\nu_{k_q +1}/(q-1)$. The proof relies on a nonlinear entropy method, and the convergence holds in a relative entropy sense. The topology of convergence can then be improved to relative error and stronger norms using regularity results of \cite{jin_xiong_bubbling}, and the rate of convergence is optimal, as shown in \cite{akagi_maekawa}.

It is worth noticing that the modes corresponding to the negative eigenvalues $\nu_1,\dots,\nu_{k_q}$ are not active only if we choose the \emph{exact} extinction time $t_0=t_*$ in \eqref{Eqn:RP}. If $t_0$ is close, but not equal to $t_*$, we should expect the first mode $\phi$ to become active for large $s$, as it is seen if we take $u_0$ a solution to the Lane-Emden-Fowler equation.

An alternative proof of the optimal rate using an energy method is also given in \cite{akagi_rate}, and it yields convergence in $H^1(\Omega)$ norm with the same rate. More precisely, if $\phi$ is non-degenerate, \cite[Corollary 1.6]{akagi_rate} states that there exists a constant $C>0$ such that
\begin{equation}\label{Eqn:optimalRate}
\left(\int_\Omega |\nabla v_{t_*}(x,s)-\nabla \phi|^2\d x \right)^{1/2}\leq C e^{-\frac{\nu_{k_q + 1}}{q-1}s},\quad \text{for }s\geq 0.
\end{equation}
Moreover, the result of \cite{akagi_rate} also allows for generalization to possibly sign-changing solutions.
For additional results on the long-time behavior of the standard fast diffusion case we also refer to \cite{bonforte_figalli_survey,choi_mcCann_seis,jin_xiong_bubbling}.

In the fractional subcritical fast diffusion case ($\theta \in (0,1)$, $q\in(2,2d/(d-2)_+)$), it is proved in \cite{akagi_existence} that $v_{t_*}(s)$ converges to an asymptotic profile in the fractional Sobolev norm. More precisely there exists $\phi\in H^\theta(\mathbb{R}^d)$ with $\phi\equiv 0$ in $\mathbb{R}^d\setminus \Omega$ such that $v_{t_*}(s)\rightarrow \phi$ in $H^\theta(\mathbb{R}^d)$ as $s\rightarrow +\infty$. The asymptotic profile is a solution to a fractional Lane-Emden-Fowler equation
\begin{equation}\label{Eqn:Fractional_Lane_Emden_Fowler}
    \left\lbrace
    \begin{aligned}
    (-\Delta)^\theta \phi&=\frac{\lambda_q}{t_*}\phi^{q-1} &&\text{on $\Omega$},\\
    \phi&=0																					&&\text{in $\mathbb{R}^d\setminus\Omega$},\\
    \phi&>0 &&\text{in $\Omega$}.
    \end{aligned}
    \right. 
\end{equation}

However, the question of the rate of convergence in the fractional case seems still open. Nevertheless, we may expect a behaviour similar to the standard case $\theta=1$, since the proof of \cite{akagi_rate} relies on a energy method that seems to be adaptable to the fractional case. More precisely, if $\mathcal{L}_{\phi,\theta}h := (-\Delta)^\theta h - (q-1)(\lambda_q/t_*)\phi^{q-2}h$ has a trivial kernel, we conjecture that there exists a constant $C>0$ such that
\begin{equation}\label{Eqn:FractionalOptimalRate}
    \|v_{t_*}(s)-\phi\|_{H^\theta(\mathbb{R}^d)} \leq C e^{-\frac{\nu_0}{q-1} s},\quad \text{for }s\geq 0,
\end{equation}
where $\nu_0$ is the smallest positive eigenvalue to the weighted eigenvalue problem
\begin{equation}\label{Eqn:fractional_weighted_eigenvalue_problem}
    \left\lbrace
    \begin{aligned}
    \mathcal{L}_{\phi,\theta}h&=\nu \phi^{q-2} h &&\text{on $\Omega$},\\
    h&=0								    &&\text{in $\partial\Omega$}.
    \end{aligned}
    \right. 
\end{equation}

In this subsection, we investigate if the convergence result \eqref{Eqn:optimalRate}, and the conjecture \eqref{Eqn:FractionalOptimalRate} are achieved by the scheme \eqref{Eqn:RP}. To do so we first compute an approximation $t_0$ of the extinction time $t_*$ by using the dichotomy method introduced in Subsection \ref{Subsec:ExtinctionTime}. Then we run the scheme \eqref{Eqn:RescaledScheme} with $t_0$ equal to this approximation. To mimic the long time asymptotics of the scheme \eqref{Eqn:RP} we fix a large final time $S=N_s\Delta s$, and we define the \emph{numerical asymptotic profile} $\mathbf{f}=(f_i)_{i\in\Z}$ by rescaling $\bv_{t_*}^{N_s}$ so as $[(-\Delta)^\theta_h \mathbf{f}]_0 = (\lambda_q/t_0)f^{q-1}_0$. The rescaling is needed since, as already mentioned, the first mode of $\mathcal{L}_{\theta,\phi}$ should be active if $t_0$ is not exactly equal to the extinction time.  Numerical simulations then indicate that $\mathbf{f}$ closely approximate a solution to \eqref{Eqn:Fractional_Lane_Emden_Fowler}. Moreover, $\mathbf{f}$ seems to converge to a solution as $\Delta s\rightarrow 0, h\rightarrow 0$ and $S\rightarrow +\infty$. Finally we compute the error $\|\bv_{t_0}^n - \mathbf{f}\|_{\Xh}$ and we compare it to $\exp(-\nu_0 s/(q-1))$ where $\nu_0$ is the smallest positive eigenvalue of the matrix $A D_{\theta,h} - (q-1)(\lambda_q/t_0)I_{M_x}$, where $I_{M_x}$ is the identity matrix of size $M_x\times M_x$, $D_{\theta,h}$ the matrix of the fractional Laplacian, and $A$ a diagonal matrix such that $A_{i,i}=f_i^{2-q}$.

Simulations indeed suggest that \eqref{Eqn:FractionalOptimalRate} holds in the case $\theta\in (0,1)$ and that \eqref{Eqn:optimalRate} holds in the case $\theta =1$.
The decay rate is faster when the initial data is symmetric. However the optimal rates of \eqref{Eqn:FractionalOptimalRate} and \eqref{Eqn:optimalRate} are obtained for non-symmetric  initial data, such as, for example, by shifting the initial data in \eqref{Eqn:DataIni},
\begin{equation}\label{Eqn:DataIniShifted}
u^0(x) =
\begin{cases}
    -x(1/2+x)&\text{if $x\in[-1,0]$},\\
    0&\text{else}.
\end{cases}
\end{equation}

\begin{figure}\label{Fig:OptimalRate}
        \begin{subfigure}[b]{0.49\textwidth}
         \centering
%
%
\definecolor{mycolor1}{rgb}{0.00000,1.00000,1.00000}%
\definecolor{mycolor2}{rgb}{1.00000,0.00000,1.00000}%
\begin{tikzpicture}

\begin{axis}[%
width=\textwidth,
height=\textwidth,
at={(0,0)},
xmin=-1,
xmax=1,
xlabel = {$x$},
ymin=0,
ymax=0.3,
ytick distance = 0.1,
axis background/.style={fill=white},
title style={font=\bfseries},
legend style={legend cell align=left, align=left, draw=white!15!black}
]
\addplot [color=black, dashed]
  table[row sep=crcr]{%
-1	0\\
-0.95	0.0475000000000001\\
-0.9	0.0900000000000001\\
-0.85	0.1275\\
-0.8	0.16\\
-0.75	0.1875\\
-0.7	0.21\\
-0.65	0.2275\\
-0.6	0.24\\
-0.55	0.2475\\
-0.5	0.25\\
-0.45	0.2475\\
-0.4	0.24\\
-0.35	0.2275\\
-0.3	0.21\\
-0.25	0.1875\\
-0.2	0.16\\
-0.15	0.1275\\
-0.0999999999999999	0.0899999999999999\\
-0.05	0.0474999999999999\\
0	-0\\
1	0\\
};
\addlegendentry{$s=0$}

\addplot [color=green, dashed, mark=+, mark options={solid, green}]
  table[row sep=crcr]{%
-1	0\\
-0.95	0.0235303560404121\\
-0.9	0.0418965419005912\\
-0.85	0.0584122969977372\\
-0.8	0.073662396987332\\
-0.75	0.0877656515505374\\
-0.7	0.100686639428078\\
-0.65	0.112333052248141\\
-0.6	0.122599470020677\\
-0.55	0.13139175728214\\
-0.5	0.138641895911223\\
-0.45	0.144316590991902\\
-0.4	0.148421026454847\\
-0.35	0.150998608601709\\
-0.3	0.152127480210007\\
-0.25	0.151914663289255\\
-0.2	0.150488743829158\\
-0.15	0.14799198742661\\
-0.0999999999999999	0.144572662266553\\
-0.05	0.140378167126629\\
0	0.135549351839207\\
0.05	0.13021621135909\\
0.0999999999999999	0.124494959398396\\
0.15	0.118486358586758\\
0.2	0.112275104633164\\
0.25	0.105930026518772\\
0.3	0.0995048628498485\\
0.35	0.0930393940699097\\
0.4	0.0865607398981676\\
0.45	0.0800846614427175\\
0.5	0.0736167296645718\\
0.55	0.0671532278289819\\
0.6	0.0606816332031321\\
0.65	0.0541804489536817\\
0.7	0.0476179759648685\\
0.75	0.0409491821180596\\
0.8	0.034108717598367\\
0.85	0.0269948617895095\\
0.9	0.0194273242598912\\
0.95	0.0110007607409921\\
1	0\\
};
\addlegendentry{$s=0.25$}

\addplot [color=red, dashed, mark=o, mark options={solid, red}]
  table[row sep=crcr]{%
-1	0\\
-0.95	0.0178093932435677\\
-0.9	0.0316267987986349\\
-0.85	0.0440724301420381\\
-0.8	0.0556817634740185\\
-0.75	0.0666319518890093\\
-0.7	0.0769756294636887\\
-0.65	0.0867099985820561\\
-0.6	0.0958050307271676\\
-0.55	0.104217670374511\\
-0.5	0.111900157032442\\
-0.45	0.118805433846437\\
-0.4	0.124890870185373\\
-0.35	0.130120851780307\\
-0.3	0.134468510393146\\
-0.25	0.137916745305843\\
-0.2	0.140458641194755\\
-0.15	0.142097371748392\\
-0.0999999999999999	0.142845676644749\\
-0.05	0.142725001429746\\
0	0.141764390398673\\
0.05	0.139999219487414\\
0.0999999999999999	0.137469848815972\\
0.15	0.134220263284234\\
0.2	0.130296755403036\\
0.25	0.125746688385525\\
0.3	0.12061736027954\\
0.35	0.114954972004321\\
0.4	0.108803683260782\\
0.45	0.102204718994414\\
0.5	0.0951954621104052\\
0.55	0.0878084286631076\\
0.6	0.080069955624426\\
0.65	0.0719983076385864\\
0.7	0.0636006547915084\\
0.75	0.0548677979574854\\
0.8	0.045764060785324\\
0.85	0.0362054830977601\\
0.9	0.0260038457129894\\
0.95	0.0146735646462459\\
1	0\\
};
\addlegendentry{$s=0.5$}

\addplot [color=blue, dashed, mark=x, mark options={solid, blue}]
  table[row sep=crcr]{%
-1	0\\
-0.95	0.0162569600420219\\
-0.9	0.0288422633835808\\
-0.85	0.0401781104714556\\
-0.8	0.0507777635450357\\
-0.75	0.0608261850019369\\
-0.7	0.0703935258193704\\
-0.65	0.0794975009105627\\
-0.6	0.088128201857359\\
-0.55	0.0962600827408144\\
-0.5	0.103858644696108\\
-0.45	0.110884620087055\\
-0.4	0.117296849058329\\
-0.35	0.123054406063044\\
-0.3	0.128118254110258\\
-0.25	0.132452570648343\\
-0.2	0.136025821087333\\
-0.15	0.138811620192323\\
-0.0999999999999999	0.140789402518581\\
-0.05	0.141944913135968\\
0	0.142270525130251\\
0.05	0.141765388533579\\
0.0999999999999999	0.140435415064149\\
0.15	0.138293103438538\\
0.2	0.135357210380997\\
0.25	0.131652272201492\\
0.3	0.127207980326602\\
0.35	0.122058410645421\\
0.4	0.116241099767113\\
0.45	0.109795949226481\\
0.5	0.10276391755994\\
0.55	0.0951854226866722\\
0.6	0.087098308058525\\
0.65	0.0785350923189538\\
0.7	0.069518944641173\\
0.75	0.0600571997237302\\
0.8	0.0501296372892159\\
0.85	0.0396640829016159\\
0.9	0.0284748841535454\\
0.95	0.0160520764883798\\
1	0\\
};
\addlegendentry{$s=1$}

\addplot [color=mycolor1, dashed, mark=diamond, mark options={solid, mycolor1}]
  table[row sep=crcr]{%
-1	0\\
-0.95	0.0161541016674314\\
-0.9	0.0286578227398306\\
-0.85	0.039920052846443\\
-0.8	0.0504524106946471\\
-0.75	0.060440216363008\\
-0.7	0.0699546464170553\\
-0.65	0.0790146807134089\\
-0.6	0.0876117064487392\\
-0.55	0.0957213698600923\\
-0.5	0.103310159954088\\
-0.45	0.110339512356483\\
-0.4	0.116768622745231\\
-0.35	0.122556527623071\\
-0.3	0.127663731423575\\
-0.25	0.132053525109268\\
-0.2	0.135693072919157\\
-0.15	0.138554307243368\\
-0.0999999999999999	0.140614651481086\\
-0.05	0.141857579734807\\
0	0.142273016444816\\
0.05	0.14185757642742\\
0.0999999999999999	0.140614644959571\\
0.15	0.138554297690747\\
0.2	0.135693060601391\\
0.25	0.132053510365475\\
0.3	0.127663714653782\\
0.35	0.122556509274117\\
0.4	0.116768603295527\\
0.45	0.110339492300358\\
0.5	0.103310139786472\\
0.55	0.0957213500622365\\
0.6	0.0876116874756516\\
0.65	0.0790146629836121\\
0.7	0.0699546303052831\\
0.75	0.0604402021965804\\
0.8	0.0504523987547298\\
0.85	0.0399200433769282\\
0.9	0.0286578159719026\\
0.95	0.0161540978930284\\
1	0\\
};
\addlegendentry{$s=3$}

\addplot [color=mycolor2]
  table[row sep=crcr]{%
-1	0\\
-0.95	0.0161541885471426\\
-0.9	0.0286579768308233\\
-0.85	0.0399202674726835\\
-0.8	0.0504526819610738\\
-0.75	0.0604405413992477\\
-0.7	0.0699550227624874\\
-0.65	0.0790151060347837\\
-0.6	0.0876121783891859\\
-0.55	0.0957218859508111\\
-0.5	0.103310717560378\\
-0.45	0.110340108644852\\
-0.4	0.116769254664808\\
-0.35	0.122557191897589\\
-0.3	0.127664424551781\\
-0.25	0.132054243372404\\
-0.2	0.135693812394628\\
-0.15	0.138555063823895\\
-0.0999999999999999	0.140615420898751\\
-0.05	0.141858357589413\\
0	0.142273798235922\\
0.05	0.141858357589353\\
0.0999999999999999	0.140615420898631\\
0.15	0.138555063823719\\
0.2	0.135693812394402\\
0.25	0.132054243372132\\
0.3	0.127664424551472\\
0.35	0.122557191897251\\
0.4	0.11676925466445\\
0.45	0.110340108644483\\
0.5	0.103310717560006\\
0.55	0.0957218859504463\\
0.6	0.0876121783888364\\
0.65	0.079015106034457\\
0.7	0.0699550227621906\\
0.75	0.0604405413989868\\
0.8	0.050452681960854\\
0.85	0.0399202674725092\\
0.9	0.0286579768306987\\
0.95	0.0161541885470731\\
1	0\\
};
\addlegendentry{$s=5$}

\end{axis}
\end{tikzpicture}%
          \caption{}
         \label{Fig:SolutionRenormPlot}
     \end{subfigure}
     \begin{subfigure}[b]{0.49\textwidth}
         \centering
         \input{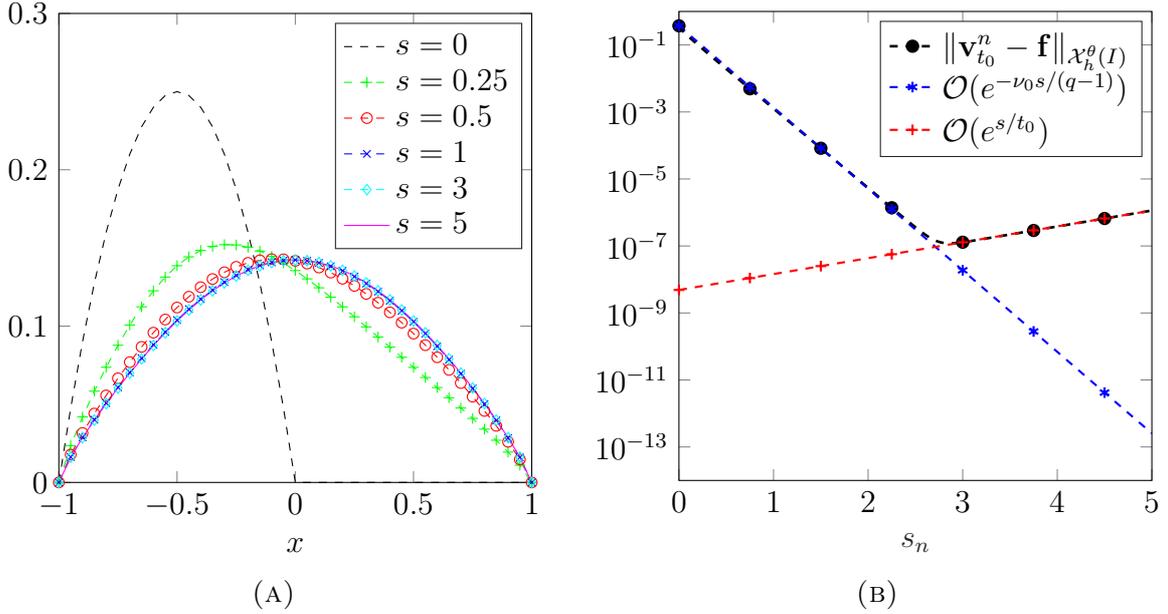}
         \caption{}
         \label{Fig:ErrorAsympProf}
     \end{subfigure}
         \caption{Convergence of the rescaled solution to an asymptotic profile, for parameters $q=2.4$, $\theta=0.75$, $L=1$, $h=0.05$, $\Delta s=0.01$, and the initial data specified in \eqref{Eqn:DataIniShifted}. The extinction time is computed by dichotomy with $\Delta s=0.01$ and $\epsilon=10^{-7}$. The rescaled solution is then computed using the scheme \eqref{Eqn:RescaledScheme}. {\scriptsize \sc (A)}: Rescaled solution at different times, showing convergence to an asymptotic profile $\mathbf{f}$, which is an approximate solution to the fractional Lane-Emden-Fowler equation \eqref{Eqn:Fractional_Lane_Emden_Fowler}. {\scriptsize \sc (B)}: Error between the rescaled solution and the asymptotic profile. Initially, the error decays at a rate corresponding to the first decaying mode of the linearized problem \eqref{Eqn:Linearized_rescaled_problem}. Subsequently, the error increases at a rate corresponding to the first mode of the linearized problem \eqref{Eqn:Linearized_rescaled_problem}, due to the imprecision on the extinction time.}
\end{figure}

On Figure \ref{Fig:SolutionRenormPlot}-\ref{Fig:ErrorAsympProf} we show results for the initial datum \eqref{Eqn:DataIniShifted}, $\theta = 0.75$, $q=2.4$, and for $L=1$, $h = 0.05$, $\Delta s=0.01$. The computed extinction time in this case is $t_*=0.92355$ and it is computed with $\Delta s=0.01$ and $\epsilon = 10^{-7}$ in the dichotomy. On Figure \ref{Fig:SolutionRenormPlot} we plot the solution $\bv^n_{t_*}$ at different times $n$, and on Figure \ref{Fig:SolutionRenormPlot} we plot the error $\|\mathbf{v}^n_{t_*}-\mathbf{f}\|_{\Xh}$ and compare it to the right-hand side of \eqref{Eqn:FractionalOptimalRate}. We observe that $\|\mathbf{v}^n_{t_*}-\mathbf{f}\|_{\Xh}$ decays in $\mathcal{O}\left(e^{-\nu_0 s/(q-1)}\right)$ until $s=2.87$, corresponding to $t \approx 0.95 t_*$. Beyond $s=2.87$, $\|\mathbf{v}^n_{t_*}-\mathbf{f}\|_{\Xh}$ increases as $\mathcal{O}(e^{s/t_0})$. It is noteworthy that $-(q-1)/t_*$ is the first eigenvalue of $\mathcal{L}_{\theta,\phi}$ for the eigenvalue problem \eqref{Eqn:fractional_weighted_eigenvalue_problem} (the eigenvalue is $\phi$). The first mode of $\mathcal{L}_{\theta,\phi}$ is inactive if $t_0$ is the \emph{exact} extinction time $t_*$ in \eqref{Eqn:RP}. However, since our computations use an \emph{approximation} of the extinction time, the error with the asymptotic profile $\|\mathbf{v}^n_{t_*}-\mathbf{f}\|_{\Xh}$ eventually increases at large times $s$. We observe that the rate is determined by the first mode of the linearized problem \eqref{Eqn:Linearized_rescaled_problem}. Other choices of $\theta\in(0,1]$ give results similar to \ref{Fig:SolutionRenormPlot}-\ref{Fig:ErrorAsympProf}, with the rate of \eqref{Eqn:optimalRate} for $\theta=1$ or \eqref{Eqn:FractionalOptimalRate} for $\theta\in(0,1)$.

\subsection{The explicit scheme}
\begin{figure}
    \centering
%
%
\begin{tikzpicture}

\begin{axis}[%
width=0.6\textwidth,
height=0.5\textwidth,
xmin=-1,
xmax=1,
xtick distance = 0.5,
ymin=-0.05,
ymax=0.25,
axis background/.style={fill=white},
legend style={legend cell align=left, align=left, draw=white!15!black}
]
\addplot [color=black, dashed]
  table[row sep=crcr]{%
-1	0\\
-0.52	-0\\
-0.48	0.0196000000000001\\
-0.44	0.0564\\
-0.4	0.0900000000000001\\
-0.36	0.1204\\
-0.32	0.1476\\
-0.28	0.1716\\
-0.24	0.1924\\
-0.2	0.21\\
-0.16	0.2244\\
-0.12	0.2356\\
-0.0800000000000001	0.2436\\
-0.04	0.2484\\
0	0.25\\
0.04	0.2484\\
0.0800000000000001	0.2436\\
0.12	0.2356\\
0.16	0.2244\\
0.2	0.21\\
0.24	0.1924\\
0.28	0.1716\\
0.32	0.1476\\
0.36	0.1204\\
0.4	0.0900000000000001\\
0.44	0.0564\\
0.48	0.0196000000000001\\
0.52	-0\\
1	0\\
};
\addlegendentry{$t=0$}

\addplot [color=green, dashed, mark=+, mark options={solid, green}]
  table[row sep=crcr]{%
-1	0\\
-0.96	0.00457366154049454\\
-0.92	0.00807840397468262\\
-0.88	0.0112879220719009\\
-0.84	0.0142645343673131\\
-0.8	0.0171009805231939\\
-0.76	0.019814068004663\\
-0.72	0.0224251518652088\\
-0.68	0.0249353895050946\\
-0.64	0.0273486171026436\\
-0.6	0.0296603083199229\\
-0.56	0.0318677939145149\\
-0.52	0.0339641615789608\\
-0.48	0.0359439668317356\\
-0.44	0.0377994201160905\\
-0.4	0.0395240417436997\\
-0.36	0.0411100606724542\\
-0.32	0.0425509854096062\\
-0.28	0.0438396977432403\\
-0.24	0.0449703961202663\\
-0.2	0.0459371032693936\\
-0.16	0.0467352239293444\\
-0.12	0.0473602968476214\\
-0.0800000000000001	0.0478093088248344\\
-0.04	0.0480795786808772\\
0	0.0481699207399036\\
0.04	0.0480795786808772\\
0.0800000000000001	0.0478093088248344\\
0.12	0.0473602968476214\\
0.16	0.0467352239293444\\
0.2	0.0459371032693938\\
0.24	0.0449703961202663\\
0.28	0.0438396977432403\\
0.32	0.0425509854096062\\
0.36	0.0411100606724542\\
0.4	0.0395240417436997\\
0.44	0.0377994201160905\\
0.48	0.0359439668317356\\
0.52	0.0339641615789608\\
0.56	0.0318677939145149\\
0.6	0.0296603083199229\\
0.64	0.0273486171026436\\
0.68	0.0249353895050946\\
0.72	0.0224251518652088\\
0.76	0.019814068004663\\
0.8	0.0171009805231939\\
0.84	0.0142645343673131\\
0.88	0.0112879220719009\\
0.92	0.00807840397468262\\
0.96	0.00457366154049454\\
1	0\\
};
\addlegendentry{$t=0.4$}

\addplot [color=red, dashed, mark=o, mark options={solid, red}]
  table[row sep=crcr]{%
-1	0\\
-0.96	0.00908712441690351\\
-0.92	-0.0101302042678926\\
-0.88	0.0137274968525445\\
-0.84	-0.00967901934253623\\
-0.8	0.0143923191651034\\
-0.76	-0.00850490959269012\\
-0.72	0.0148004685075032\\
-0.68	-0.00744173218431388\\
-0.64	0.0150912147287319\\
-0.6	-0.00637786110569416\\
-0.56	0.0148053180396703\\
-0.52	-0.00428799600095164\\
-0.48	0.0115339501000864\\
-0.44	0.00281588632946228\\
-0.4	0.00654075218399908\\
-0.36	0.00662742134689021\\
-0.32	0.00656036420158235\\
-0.28	0.00730778379802199\\
-0.24	0.00723831370549366\\
-0.2	0.00766775794909336\\
-0.16	0.00769010925229097\\
-0.12	0.00794005250168217\\
-0.0800000000000001	0.00793939974604418\\
-0.04	0.00808317338698061\\
0	0.0080201910642117\\
0.04	0.00808286938587188\\
0.0800000000000001	0.00793868359501282\\
0.12	0.00793915349451479\\
0.16	0.00768855665547652\\
0.2	0.00766640728907109\\
0.24	0.00723568124328944\\
0.28	0.00730605314558197\\
0.32	0.00655782889335077\\
0.36	0.00662077144586948\\
0.4	0.00655026255999602\\
0.44	0.00278406082606586\\
0.48	0.0115664770412369\\
0.52	-0.00430490124833005\\
0.56	0.014810284875759\\
0.6	-0.00638105643402054\\
0.64	0.0150907180381974\\
0.68	-0.00744483810992502\\
0.72	0.0147992008963447\\
0.76	-0.00850859801647963\\
0.8	0.0143906278318218\\
0.84	-0.00968209067935621\\
0.88	0.013725930585367\\
0.92	-0.0101319674097504\\
0.96	0.00908640261925164\\
1	0\\
};
\addlegendentry{$t=0.7$}

\addplot [color=blue, dashed, mark=x, mark options={solid, blue}]
  table[row sep=crcr]{%
-1	0\\
-0.96	0.00846178004397236\\
-0.92	-0.0113978955234055\\
-0.88	0.0122669804428506\\
-0.84	-0.0120092790973034\\
-0.8	0.0123575576246482\\
-0.76	-0.0118575074138654\\
-0.72	0.0123662854689923\\
-0.68	-0.011759537769201\\
-0.64	0.0124043528037843\\
-0.6	-0.0116903626828349\\
-0.56	0.0124477920405885\\
-0.52	-0.0116359784194351\\
-0.48	0.012490581688787\\
-0.44	-0.0115945501953822\\
-0.4	0.0125332585991123\\
-0.36	-0.0115705752095019\\
-0.32	0.0125852488486984\\
-0.28	-0.0115860471686009\\
-0.24	0.0126823950970945\\
-0.2	-0.0116645034218514\\
-0.16	0.0125900831984487\\
-0.12	-0.0102764574137322\\
-0.0800000000000001	0.00608421409732007\\
-0.04	0.00734576214601024\\
0	-0.00990040463159958\\
0.04	0.00734576222623828\\
0.0800000000000001	0.00608421420957894\\
0.12	-0.0102764571769767\\
0.16	0.0125900834590107\\
0.2	-0.0116645030546403\\
0.24	0.0126823954749427\\
0.28	-0.0115860466840971\\
0.32	0.0125852493202294\\
0.36	-0.0115705746360242\\
0.4	0.0125332591384895\\
0.44	-0.0115945495669736\\
0.48	0.0124905822661046\\
0.52	-0.0116359777732256\\
0.56	0.0124477926231648\\
0.6	-0.0116903620570146\\
0.64	0.0124043533577058\\
0.68	-0.0117595372010126\\
0.72	0.0123662859605902\\
0.76	-0.0118575069378939\\
0.8	0.0123575580215971\\
0.84	-0.0120092787445276\\
0.88	0.012266980713842\\
0.92	-0.0113978953227778\\
0.96	0.00846178015385224\\
1	0\\
};
\addlegendentry{$t=1$}

\end{axis}

\end{tikzpicture}%
    \caption{Solution of the explicit scheme \eqref{Eqn:E} with $q=2.4$, $\theta=0.5$, $h=0.04$, $\Delta t=0.001$, and $u^0$ defined by \eqref{Eqn:DataIni}. Instability occurs when the solution approaches $0$.}
    \label{Fig:ExplicitSchemeFDE}
\end{figure}
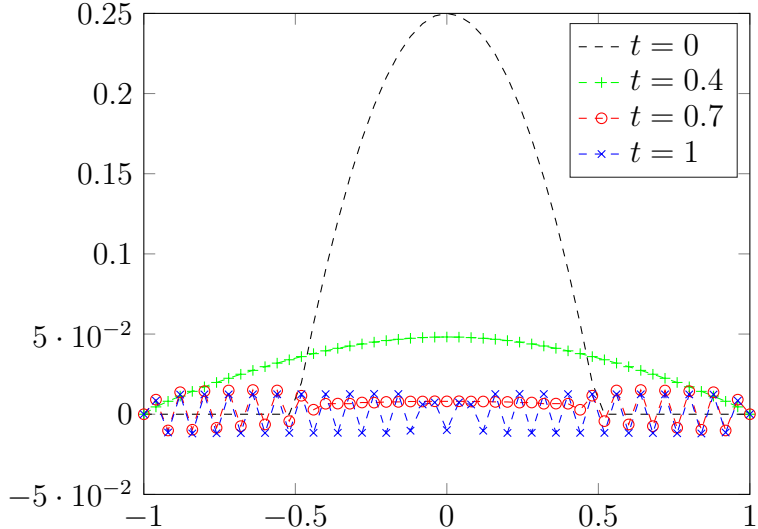

As it would seem a natural choice for the approximation of \eqref{Eqn:CDP}, in this subsection we test the behaviour of an explicit scheme, especially regarding its stability properties. The scheme we consider is the following,
\begin{equation}\label{Eqn:E}\tag{{\sf E}}
\left\{
\begin{aligned}
\frac{(u_i^{n+1})^{q-1}-(u_i^{n})^{q-1}}{\Delta t}+\left[(-\Delta)^\theta_h \bu^{n}\right]_i &= 0, && \text{for $|i|\leq M_x$ and $n\geq 0$},\\
u^n_i &= 0, && \text{for $|i|\geq M_x+1$ and $n\geq 0$},\\
u_i^0 &= (\bu^0)_i && \text{for $|i|\leq M_x$},
\end{aligned}
\right.
\end{equation}
In the fast diffusion case $q>2$ the diffusion coefficient is singular when $u\rightarrow 0$, and therefore \eqref{Eqn:E} is always unstable near the extinction time whatever the values of $\Delta t$ and $h$. This leads to a bad behavior of the scheme close to extinction time and justifies the need of an implicit scheme. On figure \ref{Fig:ExplicitSchemeFDE} we plot the solution of \eqref{Eqn:E} at different time nodes, with the parameters are $q=2.4$, $\theta=0.75$, $h=0.04$, $\Delta t=0.001$ and $u^0$ is defined by \eqref{Eqn:DataIni}, to emphasize on these instabilities.

In the porous medium case $q\in(1,2)$ the scheme \eqref{Eqn:E} is stable under an explicit stability condition,
\begin{equation}\label{Eqn:CFL}
\Delta t\sum_{j=1}^{M_x} w_j^h\leq \frac{q-1}{2}\|\bu^0\|_{\infty}^{q-2}.
\end{equation}
In particular in the non-fractional case $\theta=1$, this yields
\begin{equation}\label{Eqn:CFLlocal}
\frac{\Delta t}{h^2}\leq \frac{q-1}{2}\|\bu^0\|_{\infty}^{q-2},
\end{equation}
while the following relation holds for the fractional case, as the sum of the $w^h_j$ is known (see \eqref{Eqn:SumWeights})
\begin{equation}\label{Eqn:CFLfinal}
\frac{\Delta t}{h^{2\theta}}\leq \frac{\pi^{1/2}\Gamma(2-\theta)}{2^{2\theta}\Gamma((2\theta+1)/2)}(q-1)\|\bu^0\|_{\infty}^{q-2}.
\end{equation}
Note that \eqref{Eqn:CFLfinal} is also consistent with \eqref{Eqn:CFLlocal} when $\theta\rightarrow 1$. It is worth noticing that $q-2<0$, so the larger $\|\bu^0\|_{\infty}$ is, the more restrictive this condition is. This reflects the degeneracy of the diffusion coefficient when $u\rightarrow +\infty$ in the porous medium case. 

When \eqref{Eqn:CFLfinal} is satisfied, numerical results show that the behavior of the explicit scheme is close to that of the implicit scheme: in all the tests performed, the $\lq$-norm of the solution decays, the upper bound \eqref{Eqn:ineq_um_am} is satisfied (but the lower bound in \eqref{Eqn:Asymp2DiscreteSD} is not always satisfied), and the Rayleigh quotient decays. However, we emphasize on the fact that the stability condition \eqref{Eqn:CFLfinal} makes the scheme costly, even when its performance is compared to the implicit scheme, in particular when $\theta\rightarrow 1$ and $\|\bu^0\|_\infty$ is large.

\appendix

\section{Proof of the discrete Poincaré-Sobolev inequality}\label{Appendix:DiscretePoincareSobolev}

We prove in this section the discrete Sobolev and Poincaré inequalities stated in Lemma \ref{Lem:DiscreteSobolevIneq} and Lemma \ref{Lem:DiscretePoincareIneq}. Our approach is a straightforward adaptation of the proofs in \cite[Appendix]{savin_valdinoci} or \cite[Section 6]{di_nezza_palatucci_valdinoci} of the continuous case. It is very similar to the one of \cite{ciaurri_nonlocal_2018}, but we provide it for the sake of completeness. Moreover, it is also valid for discrete fractional Sobolev spaces of any exponent with no additional difficulty. Hence, define for any $p\in [1,\infty)$, $\theta\in (0,1)$ and $\bu\in \mathbb{Z}^\mathbb{N}$,
\begin{equation}\label{Eqn:DefFracSobolevNorm_p}
[\bu]_{W^{\theta,p}_h(\R)}:=\left(\sum_{i\in\Z}\sum_{\substack{j\in\Z\\j\neq i}}h^2\frac{|u_i-u_j|^p}{|h(i-j)|^{1+p\theta}}\right)^{1/p},
\end{equation}
and $W^{\theta,p}_h(\R):=\{\bu\in l^p_h(\R):\:[\bu]_{W^{\theta,p}_h(\R)}<+\infty\}$, endowed with the norm $\|\bu\|_{W^{\theta,p}_h(\R)} :=\|\bu\|_{\ell^p_h(\R)} + [\bu]_{W^{\theta,p}_h(\R)}$. Then it holds $H^\theta_h(\R) = W^{\theta,2}_h(\R)$, and the semi-norm $[\cdot]_{W^{\theta,2}_h(\R)}$ is equivalent to $[\cdot]_{H^\theta(\R)}=\|\cdot\|_{\Xh}$ thanks to Proposition \ref{Prop:EstimatesWeightsDiscreteFracLap}, with constants independant of $h$. Lemma \ref{Lem:DiscreteSobolevIneq} is then a straightforward consequence of the following proposition with $p=2$.

\begin{prop}[Sobolev inequality]\label{Prop:discreteSobolev}
Let $\theta\in(0,1)$, and $p\in[1,\infty)$. Let $p_\theta^*= p/(1-p\theta)_+ \in(p,\infty]$. There exists a positive constant $C=C(p,\theta)>0$ such that 
\begin{enumerate}
\item if $p\theta<1$, then for any $\bu\in W^{\theta,p}_h(\R)$, $\|\bu\|_{l^{p_\theta^*}_h(\R)}\leq C [\bu]_{W^{\theta,p}_h(\R)}$,
\item if $p\theta\geq 1$, then for any $q\in[p,\infty)$ and $\bu\in W^{\theta,p}_h(\R)$, $\|\bu\|_{l^q_h(\R)}\leq C\|\bu\|_{W^{\theta,p}_h(\R)}$.
\end{enumerate}
\end{prop}
The first ingredient is the following estimate on the discrete kernel (see \cite[Lemma 6.1]{di_nezza_palatucci_valdinoci}, \cite[Lemma A.1]{savin_valdinoci_2} for the continuous case).
\begin{lem}\label{Lem:EstimateKernel}
Let $p>1$, $\theta\in(0,1)$, let $A\subset \Z$ be finite and nonempty, and fix $i\in A$. There exists a positive constant $C=C(p,\theta)$ depending only on $p$ and $\theta$ (and not $i$) such that 
\begin{equation}
\label{Eqn:EstimateKernelConst}
\sum_{j\in \Z\setminus A}h\frac{1}{|h(i-j)|^{1+p\theta}}\geq C(p,\theta)(|A|h)^{-p\theta},
\end{equation}
where $|A|$ is the cardinal of $A$.
\end{lem}
The discrete Poincaré inequality stated in Lemma \ref{Lem:DiscretePoincareIneq} is a straightforward consequence of Lemma \ref{Lem:EstimateKernel}:

\begin{prf}{Lemma \ref{Lem:DiscretePoincareIneq}}
Let $\bu\in\Xh$, then
$$
\begin{aligned}
[\bu]_{\mathcal{X}^\theta(\R)}=&\frac{1}{2}\sum_{|i|,|j|\leq M_x}hw_j^h|u_i-u_{i-j}|^2 + \frac{1}{2}\sum_{|i|\leq M_x,|j|\geq M_x+1} hw_j^h|u_i|^2\\
\geq& b_\theta\sum_{|i|\leq M_x}h |u_i|^2\sum_{|j|\geq M_x+1} h\frac{1}{|h(i-j)|^{1+2\theta}},
\end{aligned}
$$
where we used Proposition \ref{Prop:EstimatesWeightsDiscreteFracLap}. We conclude using Lemma \ref{Lem:EstimateKernel} with $A=\{|j|\leq M_x\}$.
\end{prf}

\begin{prf}{Lemma \ref{Lem:EstimateKernel}}
Let us define
$$
\rho:=\begin{cases}
(|A|-1)/2 &\text{if $|A|$ is odd},\\
|A|/2 &\text{if $|A|$ is even},
\end{cases}
$$
and $B_\rho(i)=\{i-\rho,\dots,i+\rho$\}. Note that $|B_\rho(i)|=2\rho+1$, that is $|A|$ if $|A|$ is odd and $|A|+1$ if $|A|$ is even, and in any case $|B_\rho(i)|\geq|A|$. Therefore, thanks to the decay of $x>0\mapsto 1/x^{1+p\theta}$,
$$
\sum_{j\in \Z\setminus A}h\frac{1}{|h(i-j)|^{1+p\theta}}\geq \sum_{j\in \Z\setminus B_\rho(i)}h\frac{1}{|h(i-j)|^{1+p\theta}} = \sum_{|j|\geq\rho+1}h\frac{1}{|hj|^{1+p\theta}}.
$$
\details{More precisely, since $|(\Z\setminus A)\cap B_\rho(i)|=|B_\rho(i)|-|A\cap B_\rho(i)|\geq |A|-|A\cap B_\rho(i)|=|A\cap (\Z\setminus B_\rho(i))|$, the following estimate holds true,
$$
\begin{aligned}
\sum_{j\in \Z\setminus A}h\frac{1}{|h(i-j)|^{1+p\theta}}=&\sum_{j\in (\Z\setminus A)\cap B_\rho(i)}h\frac{1}{|h(i-j)|^{1+p\theta}} + \sum_{j\in (\Z\setminus A)\cap (\Z\setminus B_\rho(i))}h\frac{1}{|h(i-j)|^{1+p\theta}}\\
\geq&\sum_{j\in (\Z\setminus A)\cap B_\rho(i)}h\frac{1}{|h\rho|^{1+p\theta}} + \sum_{j\in (\Z\setminus A)\cap (\Z\setminus B_\rho(i))}h\frac{1}{|h(i-j)|^{1+p\theta}}\\
=&\frac{|(\Z\setminus A)\cap B_\rho(i)|}{|h\rho|^{1+p\theta}} + \sum_{j\in (\Z\setminus A)\cap (\Z\setminus B_\rho(i))}h\frac{1}{|h(i-j)|^{1+p\theta}}\\
\geq& \frac{|A\cap (\Z\setminus B_\rho(i))|}{|h\rho|^{1+p\theta}} + \sum_{j\in (\Z\setminus A)\cap (\Z\setminus B_\rho(i))}h\frac{1}{|h(i-j)|^{1+p\theta}}\\
\geq& \sum_{j\in A\cap (\Z\setminus B_\rho(i))}h\frac{1}{|h(i-j)|^{1+p\theta}} + \sum_{j\in (\Z\setminus A)\cap (\Z\setminus B_\rho(i))}h\frac{1}{|h(i-j)|^{1+p\theta}}\\
=&\sum_{j\in \Z\setminus B_\rho(i)}h\frac{1}{|h(i-j)|^{1+p\theta}}\\
=& \sum_{|j|\geq\rho+1}h\frac{1}{|hj|^{1+p\theta}}.
\end{aligned}
$$}
To conclude, the series of the right-hand side is compared to an integral, by using again the decay of $x>0\mapsto 1/x^{1+p\theta}$,
$$
\begin{aligned}
\sum_{|j|\geq\rho+1}h\frac{1}{|hj|^{1+p\theta}}\geq& 2\int_{(\rho+1)h}^{+\infty} \frac{1}{x^{1+p\theta}}\d x
=\frac{2}{p\theta}(\rho+1)^{-p\theta}h^{-p\theta}\\
\geq& \frac{2}{p\theta}(1+|A|/2)^{-p\theta}h^{-p\theta}
\geq \frac{2}{p\theta}\left(\frac{3}{2}\right)^{-p\theta}(|A|h)^{-p\theta}.
\end{aligned}
$$
\end{prf}

We now turn to the proof of Proposition \ref{Prop:discreteSobolev}. We start by recalling a general lemma (see \cite[Lemma 6.2]{di_nezza_palatucci_valdinoci},\cite[Lemma 5]{savin_valdinoci}) that will be useful to estimate from below $[\bu]_{W^{\theta,p}_h(\R)}$.  

\begin{lem}\label{algebre}
Let $s\in(0,1)$, and $T>1$. Let $N\in\Z$ and $(a_k)_{k\in\Z}$ be a bounded, nonnegative, non-increasing sequence with $a_k=0$ for any $k\geq N$. Then,
$$
\sum_{k\in\Z} a_k^{1-s}T^k\leq C \sum_{k\in\Z, a_k\neq 0} a_{k+1}a_k^{-s}T^k,
$$
for a suitable constant $C=C(s,T)>0$, independant of $N$.
\end{lem}
\details{
\begin{pr}
With a change of index in the sum, one has
$$
\begin{aligned}
\frac{1}{T}\sum_{k\in\Z, a_k\neq 0} a_k^{1-s}T^k=& \frac{1}{T}\sum_{k\in\Z, a_{k+1}\neq 0} a_{k+1}^{1-\th}T^{k+1}\\
=& \sum_{k\in\Z, a_{k+1}\neq 0} a_{k+1}^{1-s}T^{k}\\
=&\sum_{k\in\Z,a_{k+1}\neq 0} a_k^{s(1-s)}a_{k+1}^{1-s}a_k^{-s(1-s)}T^k\\
\leq&\left(\sum_{k\in\Z, a_{k+1\neq 0}} a_k^{1-s}\right)^{s}\left(\sum_{k\in\Z, a_{k+1\neq 0}} a_{k+1}a_k^{-s}T^k\right)^{1-s}\\
\leq&\left(\sum_{k\in\Z, a_{k\neq 0}} a_k^{1-s}\right)^{s}\left(\sum_{k\in\Z, a_{k\neq 0}} a_{k+1}a_k^{-s}T^k\right)^{1-s},
\end{aligned}
$$
where the inequality in the fourth line has been obtained using Holder inequality. Therefore,
$$
\sum_{k\in\Z, a_k\neq 0} a_k^{1-s}T^k\leq T^{1/(1-s)}\sum_{k\in\Z, a_{k\neq 0}} a_{k+1}a_k^{-s}T^k.
$$
\end{pr}
}

\begin{prf}{Proposition \ref{Prop:discreteSobolev}} Let $p\in[0,+\infty)$, $\theta\in(0,1)$ and $\bu\in W^{\theta,p}_h(\R)$. Using a cutoff argument similar to the continuous case, it can be shown that compactly supported sequences are dense in $W^{\theta,p}_h(\R)$. Hence, we assume without loss of generality that $u$ vanishes outside a finite number of indices, and thus, $\bu$ is also bounded. \details{In the continuous case, one would normally use a smooth cutoff function to perform the approximation argument. Here, since we are working with $h$ fixed, and there is no need for uniform estimates in $h$ at this point, a coarse cutoff function is enough. Indeed, let $N\geq 1$ and $\bu^N=(u^N_i)_{i\in\Z}$ with
$$
u_i^N=\begin{cases}
u_i &\text{if $|i|\leq N$},\\
0 &\text{if $|i|> N$},
\end{cases}
$$
Then,
$$
\begin{aligned}
[\bu-\bu^N]_{W^{\theta,p}_h(\R)}^p &= \sum_{i,j\in\Z, i\neq j}h^2 \frac{|u_i-u_i^N - (u_j-u_j^N)|^p}{|h(i-j)|^{1+p\theta}}\\
&= 2\sum_{|i|>N, |j|\leq N}h^2 \frac{|u_i|^p}{|h(i-j)|^{1+p\theta}} + \sum_{|i|>N,|j|>N}h^2\frac{|u_i-u_j|^p}{|h(i-j)|^{1+p\theta}}.
\end{aligned}
$$
The second term on the right-hand side converges to $0$ by dominated convergence. The first one is treated with crude, non-uniform in $h$ estimates
$$
\sum_{|i|>N, |j|\leq N}h^2 \frac{|u_i|^p}{|h(i-j)|^{1+p\theta}} \leq \sum_{|i|>N}h|u_i|^p \sum_{j=1}^{2N+1}h\frac{1}{|hj|^{1+p\theta}}
\leq h^{-p\theta}\sum_{|i|>N}h|u_i|^p\sum_{j=1}^\infty \frac{1}{j^{1+p\theta}} \underset{N\rightarrow\infty}{\longrightarrow} 0.
$$
Moreover, it is also clear that for any $q\in[p,\infty)$
$$
\|\bu-\bu^N\|_{l^q_h(\R)}\underset{N\rightarrow\infty}{\longrightarrow} 0.
$$}

\noindent
\emph{\underline{Case $\theta p<1$:}} Define
$$
A_k:=\{i\in\Z:\:|u_i|>2^k\},\quad k\in\Z,
$$
and
$$
D_k=A_k\setminus A_{k+1},\quad k\in\Z,\qquad D_{-\infty}=\{i\in\Z:\:u_i=0\}.
$$
Note that $A_{k+1}\subset A_k$ for any $k\in\Z$, and that $A_k=\emptyset$ for $k$ large enough. We have moreover,
\begin{equation}\label{union1}
\bigcup_{-\infty\leq l \leq k-1}D_l=A_k^c,\quad \bigcup_{l\geq k}D_l=A_k.
\end{equation}
Finally we introduce
$$
a_k=|A_k|h,\quad d_k=|D_k|h,\quad\text{for $k\in\Z$}.
$$
In, particular $(d_k)_{k\in\Z}$ and $(a_k)_{k\in\Z}$ vanish for $k$ large enough, $(a_k)_{k\in\Z}$ is a non-increasing sequence, and \eqref{union1} yields
\begin{equation}\label{union2}
\sum_{l\geq k} d_l=a_k.
\end{equation}
First $[\bu]_{W^{\theta,p}_h(\R)}$ is estimated from below by considering only indices $k\in D_i,l\in D_j$ with $|i-j|\geq 2$ in the definition of $[\bu]_{W^{\theta,p}_h(\R)}$ \eqref{Eqn:DefFracSobolevNorm_p},
$$
\begin{aligned}
[\bu]_{W^{\theta,p}_h(\R)}^p=&\sum_{\substack{-\infty\leq i,<+\infty\\ -\infty\leq j<+\infty}}\sum_{\substack{(k,l)\in D_i\times D_j\\k\neq l}}h^2\frac{|u_k-u_l|^p}{|h(k-l)|^{1+p\theta}},\\
\geq& 2\sum_{i\in\Z}\sum_{-\infty\leq j\leq i-2}\sum_{\substack{(k,l)\in D_i\times D_j\\k\neq l}}h^2\frac{|u_k-u_l|^p}{|h(k-l)|^{1+p\theta}}.
\end{aligned}
$$
Then, for $i\in\Z$ and $-\infty\leq j\leq i-2$, and $(k,l)\in D_i\times D_j$ it holds $|u_i-u_j|\geq 2^i-2^{i-1}=2^{i-1}$. Hence,
$$
\begin{aligned}
[\bu]_{W^{\theta,p}_h(\R)}^p\geq& 2\sum_{i\in\Z}2^{p(i-1)}\sum_{-\infty\leq j\leq i-2}\sum_{\substack{(k,l)\in D_i\times D_j\\k\neq l}}h^2\frac{1}{|h(k-l)|^{1+p\theta}}\\
=& 2\sum_{i\in\Z}2^{p(i-1)}\sum_{\substack{(k,l)\in D_i\times A_{i-1}^c\\k\neq l}}h^2\frac{1}{|h(k-l)|^{1+p\theta}}.
\end{aligned}
$$
Using Lemma \ref{Lem:EstimateKernel} to estimate the inner sum, and denoting $C=C(\theta,p)$ the constant in \eqref{Eqn:EstimateKernelConst} yields
\begin{equation}\label{calcul}
[\bu]_{W^{\theta,p}_h(\R)}^p\geq 2C\sum_{\substack{i\in\Z\\ a_{i-1}\neq 0}}2^{p(i-1)}a_{i-1}^{-p\theta}d_i.
\end{equation}
Next we use \eqref{union2} to obtain
\begin{equation}\label{final}
[\bu]_{W^{\theta,p}_h(\R)}^p\geq 2C\sum_{\substack{i\in\Z\\ a_{i-1}\neq 0}}2^{p(i-1)} a_{i-1}^{-p\theta}a_i  -  2C\sum_{\substack{i\in\Z\\ a_{i-1}\neq 0}}\sum_{l\geq i+1}2^{p(i-1)}a_{i-1}^{-p\theta}d_l.
\end{equation}
By definition of the sequences $(a_i)_{i\in\Z}$ and $(d_l)_{l\in\Z}$, it holds for any $i\in\Z$
$$
a_{i-1}\neq 0 \quad \Leftrightarrow \quad \exists l\geq i-1,\  d_l\neq 0.
$$
Therefore, adopting the convention $a_{i-1}^{-p\theta}d_l = 0$ whenever $d_l=0$, the last term can be written
$$
\sum_{\substack{i\in\Z\\ a_{i-1}\neq 0}}\sum_{l\geq i+1}2^{p(i-1)}a_{i-1}^{-p\theta}d_l = \sum_{i\in\Z}\sum_{l\geq i+1}2^{p(i-1)}a_{i-1}^{-p\theta}d_l.
$$
Using the decay of the sequence $(a_i)_{i\in\Z}$, exchanging the order of summation, and using \eqref{calcul} yields
$$
\begin{aligned}
\sum_{i\in\Z}\sum_{l\geq i+1}2^{p(i-1)}a_{i-1}^{-p\theta}d_l
\leq& 
\sum_{i\in\Z}\sum_{l\geq i+1}2^{p(i-1)}a_{l-1}^{-p\theta}d_l\\
=&
\sum_{l\in\Z}\sum_{i\leq l-1}2^{p(i-1)}a_{l-1}^{-p\theta}d_l\\
=&\frac{1}{2^p-1}\sum_{\substack{l\in\Z\\ a_{l-1}\neq 0}} 2^{p(l-1)}a_{l-1}^{-p\theta}d_l \leq \frac{1}{2C}\frac{1}{2^p-1}[\bu]_{W^{\theta,p}_h(\R)},
\end{aligned}
$$
where $C$ is the constant in \eqref{calcul}. Therefore,
\begin{equation}\label{estimation1}
[\bu]_{W^{\theta,p}_h(\R)}^p\geq 2C\frac{2^p-1}{2^p} \sum_{\substack{i\in\Z\\ a_{i-1}\neq 0}}2^{pi} a_{i-1}^{-p\theta}a_i.
\end{equation}
Let us on the other hand estimate from above $\|\bu\|_{l^{p^*_\theta}_h(\R)}$. We have
$$
\begin{aligned}
\|\bu\|_{l^{p^*_\theta}_h(\R)}^{p^*_\theta}=&\sum_{k\in\Z}\sum_{i\in D_k} h|u_i|^{p^*_\theta}\leq \sum_{k\in\Z} 2^{(k+1)p^*_\theta}d_k\\
\leq& \sum_{k\in\Z} 2^{(k+1)p^*_\theta}a_k.
\end{aligned}
$$
Hence,
$$
\|\bu\|_{l^{p^*_\theta}_h(\R)}^{p}\leq 2^p\left(\sum_{k\in\Z} 2^{kp^*_\theta}a_k\right)^{p/p^*_\theta}.
$$
Since $p/p^*_\theta=1-\theta p\in(0,1)$, it holds
\begin{equation}\label{estimate2}
\|\bu\|_{l^{p^*_\theta}_h(\R)}^{p}\leq 2^p\sum_{k\in\Z} 2^{kp}a_k^{1-\theta p}.
\end{equation}

The link between \eqref{estimation1} and \ref{estimate2} is done using Lemma \ref{algebre}. Indeed, taking $T=2^p$ and $s=\th p\in(0,1)$ in Lemma \ref{algebre} yields, for some constant $C_2=C_2(\theta,p)$, 
\begin{equation}\label{estimate3}
\sum_{k\in\Z} 2^{kp}a_k^{1-\theta p}\leq C_2\sum_{\substack{i\in\Z\\ a_{i-1}\neq 0}}2^{pi} a_{i-1}^{-p\theta}a_i.
\end{equation}
Combining \eqref{estimation1}, \eqref{estimate2}, \eqref{estimate3}, yields, up to relabelling the constant $C=C(\theta,p)$,
$$
\|\bu\|_{l^{p^*_\theta}_h(\R)}\leq C[\bu]_{W^{\theta,p}_h(\R)}.
$$
This conclude the proof of (1) of Proposition \ref{Prop:discreteSobolev}.

\noindent
\emph{\underline{Case $\theta p\geq 1$:}} Let $q\in [p,\infty)$. Let $\theta_0\in((1/p-1/q)_+,1/p)$ so that $p_{\theta_0}^* := p/(1-p\theta_0)_+ \geq q$ and $\theta_0 p<1$. Then by the previous case and \cite[Proposition 4.8]{ayi_herda_hivert_tristani}, there exists a positive constant $C(\theta,p)$ such that
$$
\|\bu\|_{l^{q}_h(\R)}\leq C(\theta,p)[\bu]_{W^{\theta_0,p}_h(\R)}\leq C(\theta,p) \|\bu\|_{W^{\theta,p}_h(\R)}.
$$
This conclude the proof of (2) of Proposition \ref{Prop:discreteSobolev}.
\end{prf}

\printbibliography
\end{document}